\documentclass{article}

\usepackage[english]{babel}
 \usepackage[letterpaper,top=2cm,bottom=2cm,left=3cm,right=3cm,marginparwidth=1.75cm]{geometry}
\usepackage{authblk}
\usepackage{amsmath}
\usepackage{amssymb}
\usepackage{amsthm}
\usepackage{color}
\usepackage{bm}
\usepackage[colorlinks,linkcolor=blue,anchorcolor=green,citecolor=red]{hyperref}
\usepackage{appendix}
\usepackage{graphicx}
\usepackage[utf8]{inputenc}
\usepackage{subfigure}
\usepackage{booktabs}
\usepackage{caption}
\usepackage{cancel}
\usepackage{float}
\usepackage[T1]{fontenc}  
\usepackage{indentfirst} 
\usepackage{mathrsfs}

 \usepackage{algorithm}
\usepackage{algpseudocode}

\allowdisplaybreaks[4]
\numberwithin{equation}{section}
\usepackage{epstopdf}

\newcommand\sH{{\mathcal H}}

\newcommand\sL{{\mathcal L}}

\newcommand\vel{\mathbf{u}}

\newcommand\bn{\mathbf{n}}

\newcommand{\ipd}[1]{\bigl(#1\bigr)}
\newcommand{\Gt}{{\Gamma(t)}}

\newtheorem{thm}{Theorem}[section]
\newtheorem{cor}[thm]{Corollary}
\newtheorem{lma}[thm]{Lemma}

\newtheorem{rmk}{Remark}[section]
\usepackage{authblk}
\title{Thermodynamically Consistent Modeling and Stable ALE Approximations  of Reactive Semi-Permeable Interfaces}
\author[1]{Weidong Shi}
\author[2]{Shixin Xu\footnote{Corresponding author: shixin.xu@dukekunshan.edu.cn}}
\author[3,4]{Zhen Zhang}
\author[5]{Quan Zhao}

\affil[1]{School of Applied Mathematics, Shanxi University of Finance and Economics, Taiyuan 030006, China }
\affil[2]{Zu Chongzhi Center, Duke Kunshan University, 8 Duke Ave, Kunshan, Jiangsu, China}
\affil[3]{Department of Mathematics, Southern University of Science and Technology (SUSTech), Shenzhen 518055, China}
\affil[4]{National Center for Applied Mathematics (Shenzhen), Southern University of Science and Technology (SUSTech), Shenzhen 518055, China}
\affil[5]{School of Mathematical Sciences, University of Science and Technology of China, Hefei 230026, China}
\begin{document}
\date{ }

\maketitle

\begin{abstract}
Reactive, semi-permeable interfaces play important roles in key biological processes such as targeted drug delivery, lipid metabolism, and signal transduction. These systems involve coupled surface reactions, transmembrane transport, and interfacial deformation, often triggered by local biochemical signals. The strong mechanochemical couplings complicate the modeling of such interfacial dynamics. We propose a thermodynamically consistent continuum framework that integrates bulk fluid motion, interfacial dynamics, surface chemistry, and selective solute exchange, derived via an energy variation approach to ensure mass conservation and energy dissipation. 
To efficiently solve the resulting coupled system, we develop a finite element scheme within an Arbitrary Lagrangian–Eulerian (ALE) framework, incorporating the Barrett–Garcke–Nürnberg (BGN) strategy to maintain mesh regularity and preserve conservation laws. Numerical experiments verify the convergence and conservation properties of the scheme and  demonstrate its ability in capturing complex interfacial dynamics.
Two biologically inspired examples showcase the model's versatility: cholesterol efflux via the ABCG1 pathway, involving multistage interfacial reactions and HDL uptake; and a self-propelled droplet system with reaction-activated permeability, mimicking drug release in pathological environments. 
This work provides a unified computational platform for studying strongly coupled biochemical and mechanical interactions at interfaces, offering new insights into reactive transport processes in both biological and industrial contexts.

\end{abstract}


%%%%%%%%%%%%%%%%%%%%%%%%%%%%%%%
\section{Introduction}\label{sec:Introduction}

Interfaces in biological systems often exhibit a complex interplay of fluid dynamics, membrane deformation, chemical reactivity, and selective permeability. A compelling example arises in lipid-based vesicular platforms for targeted drug delivery. These carriers, such as enzyme-responsive liposomes or redox-sensitive polymersomes, encapsulate therapeutic agents and remain stable in circulation. Upon reaching specific pathological environments, such as tumor tissues characterized by abnormal enzymatic activity, oxidative stress, or acidic pH, the vesicle surface undergoes chemical modifications that dramatically alter its mechanical and transport properties. Such modifications can activate membrane permeability, trigger morphological changes, and induce Marangoni-driven flows due to surface tension gradients, ultimately facilitating the controlled release of cargo. This highly coupled process involves fluid-structure interaction, bulk-interface communication,  mass transportation, and interfacial reaction dynamics.

For fluid-structure interaction, a wide range of mathematical models and numerical methods have been developed, including the immersed boundary method~\cite{peskin2002immersed}, immersed interface method~\cite{li1997immersed}, level set methods~\cite{osher2001level}, diffuse interface methods~\cite{anderson1998diffuse,yue2004diffuse}, and lattice Boltzmann methods~\cite{inamuro2004lattice}. These methods have advanced the modeling of fluid–structure interaction and multiphase flows, and comprehensive overviews can be found in~\cite{hou2012numerical,dowell2001modeling,kleinstreuer2017two}. 
The second feature of this problem is \textit{ the presence of bulk interface communication}, like surfactants or amphiphilic molecules that adsorb preferentially at interfaces and locally reduce surface tension. Surfactants are commonly classified as \textit{insoluble} (confined to the interface) or \textit{soluble} (capable of exchanging with the bulk via adsorption and desorption). The resulting gradients in surface tension generate Marangoni stresses, which can drive significant fluid motion and interfacial deformation. These effects have been shown to play important roles in interfacial instabilities, droplet migration, and film rupture. Several continuum models and numerical schemes have been proposed to account for surfactant dynamics in multiphase flows~\cite{pawar1996marangoni,james2004surfactant,lai2008immersed,teigen2011diffuse,zhao2021thermodynamically,zhu2019thermodynamically}. The foundational evolution equations for surfactants on deforming interfaces were derived by Stone~\cite{stone1990simple} and Wong et al.~\cite{wong1996surfactant}.

In parallel, modeling mass transport across \textit{semi-permeable membranes} introduces additional challenges. Such interfaces allow selective solute exchange and are crucial to processes including ion regulation in biological cells, osmotic pressure regulation in vesicles, and drug delivery across lipid bilayers~\cite{stillwell2013membrane,metheny2012fluid,gong2014immersed,zeng2023mathematical,qin2022phase}. Mathematically, these problems involve coupling flux boundary conditions for solute transport with interfacial dynamics and fluid flow. 
Further complexity arises when \textit{chemical reactions occur on the interface}.  Chemical modifications of membrane proteins, such as phosphorylation, acetylation, or oxidation, can significantly influence both their permeability to ions and their mechanical properties. This coupling between biochemical reactions and transport processes is fundamental to the regulation of cellular excitability and homeostatic regulation. Wang et al. \cite{wang2020field}  proposed a framework that incorporates bulk reactions with diffusion via an energy variation method. Later on Xu et al. \cite{xu2023coupled} generalized the derivation for interfacial reaction, offering a rigorous thermodynamic formulation for coupling reaction kinetics with fluid-interface interactions.

While the individual components have been extensively studied, integrating these mechanisms into a unified, thermodynamically consistent framework remains a substantial challenge. A central difficulty lies in accurately capturing the coupled interplay among surface forces, bulk and interfacial transport, reaction kinetics, and fluid-structure interactions, while simultaneously ensuring mass conservation and energy dissipation.
The first objective of this work is to develop a continuum model that fully couples fluid-structure interaction with reaction-modulated, semi-permeable interfaces in a thermodynamically consistent manner. We adopt an energy variation framework~\cite{liu2019energetic, qin2022phase} to derive a set of governing equations that rigorously satisfy mass and energy principles. The resulting system comprises a nonlinear coupled diffusion–convection–reaction problem on evolving domains with dynamic interface conditions and selective permeability.

% By systematically integrating bulk fluid flow, interfacial deformation, surface chemistry, and transmembrane transport, this model generalizes and unifies existing approaches that have previously addressed only isolated mechanisms. In later sections, we develop a finite element method within an Arbitrary Lagrangian–Eulerian (ALE) framework to efficiently solve the system, and we validate the model through simulations of biologically inspired scenarios such as droplet migration, bubble dynamics, and cholesterol transport via ABCG1-mediated efflux.

To accurately resolve the interfacial dynamics and preserve fundamental conservation laws, the ALE framework is adopted. ALE methods have been one of the numerical approaches to address free boundary problems by deforming the computational mesh not only to follow the interface motion but also to allow for flexibility in choosing the reference velocity for the interior points. The original ALE methods have been introduced for hydrodynamic problems \cite{Noh64, Franck1964mixed,Hirt1974arbitrary} and then generalized to other problems such as free surface flows, fluid-structure interaction, Stokes interface problems, see e.g.,\cite{Belytschko78, donea82, Soulaimani1998}. For the two-phase flow, an energy-stable ALE method was proposed in \cite{Duan2022energy} based the divergence-free velocity of the ALE frame. More recently, Garcke et al.~\cite{garcke2023structure, GKZ24ALE} devised ALE approximations that preserve the structure for the two-phase incompressible flow, where unconditional stability and exact volume preservation can be shown for fully discrete solutions. We note that both of the works in \cite{Duan2022energy, garcke2023structure} rely on a novel BGN formulation for the interfacial equations, which allows for tangential velocity to improve the interface mesh quality. In particular,  the interface is advected only in the normal direction according to the normal component of the fluid velocity, while the tangential degrees of freedom are used to optimize the mesh distribution,  see \cite{barrett2020parametric, BGN2013eliminating}. 

In the presence of tangential velocity during interface motion, the ALE finite element method can also be employed for solving PDEs on evolving surfaces, see \cite{Elliott12ale}. There have been extensions of ALE approaches to biointerface problems. For instance, MacDonald et al.~\cite{macdonald2016computational} proposed a moving-mesh finite element method to solve coupled bulk–surface reaction–diffusion systems in evolving cell domains. Building upon this framework, Mackenzie et al.~\cite{mackenzie2021conservative} developed a conservative ALE finite element scheme that rigorously enforces mass conservation across bulk–surface transport processes, ensuring that the total chemical content (membrane-bound and cytosolic) remains invariant up to numerical solver tolerance.  Besides, we note that the BGN approach has also been extended to simulate the two-phase flow with insoluble or soluble surfactants \cite{Barrett15insol, Barrett15sol}.

In the current work, we would like to combine the stable ALE approximation in \cite{garcke2023structure} with the ALE evolving surface finite element method in \cite{Elliott12ale},  and develop an ALE-based scheme for simulating semi-permeable interfaces with interfacial chemical reactions. The introduced ALE approach ensures discrete conservation laws for both bulk and interfacial species and provides accurate resolution of the coupled fluid-interface dynamics.

This paper is organized as follows. In Section~\ref{sec:Model_Derivation}, we present the mathematical model for semi-permeable interfaces with interfacial chemical reactions. Section~\ref{sec:Weak_Formulation} introduces the corresponding weak formulation of the system, followed by the finite element discretization detailed in Section~\ref{sec:Finite_Element_Approximations}. Section~\ref{sec:Numerical_Results} provides numerical experiments to validate the proposed model and demonstrate the effectiveness of the numerical method. Finally, concluding remarks are given in Section~\ref{sec:concluding}.

\section{Model Derivation}\label{sec:Model_Derivation}
\begin{figure}[htbp!]
    \centering
    \includegraphics[width=0.5\textwidth]{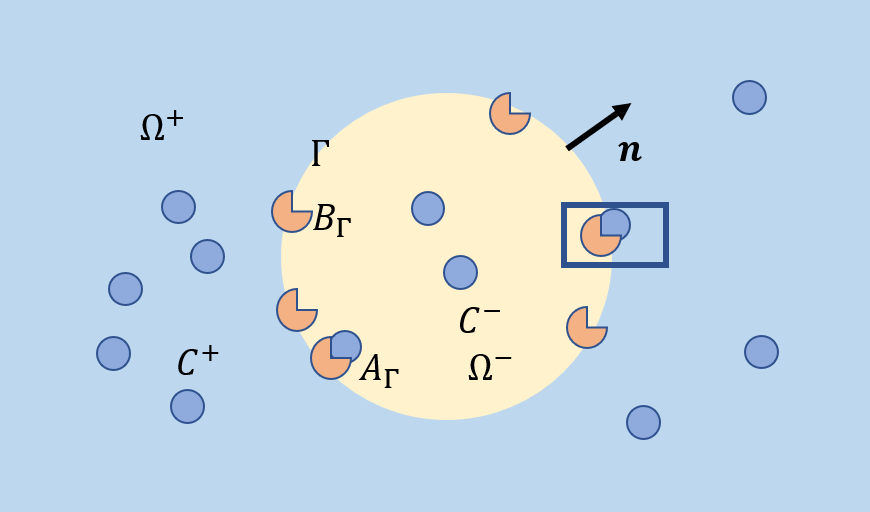}
    \caption{Schematic of mass transport  across the reactive interface.}
    \label{fig:Close_Mem}
\end{figure}

In this work, we consider a computational domain $\Omega \subset \mathbb{R}^d$ ($d = 2$ or $3$), as illustrated in Fig.~\ref{fig:Close_Mem}. The domain consists of two subdomains: the intracellular region $\Omega^-$ and the extracellular region $\Omega^+$, which are separated by a semi-permeable interface denoted by $\Gamma$. The interface $\Gamma$ allows selective leakage of substances between $\Omega^+$ and $\Omega^-$. Consequently, the domain can be decomposed as
\[
\Omega = \Omega^+ \cup \Gamma \cup \Omega^-.
\]
Let $\bn$ denote the unit normal vector on the interface $\Gamma$, oriented from the intracellular region $\Omega^-$ toward the extracellular region $\Omega^+$. For notational simplicity, we also use $\bn$ to represent the outward unit normal vector on the boundary $\partial\Omega$.

Let $C$ denote the concentration of a solute that is present in both bulk regions $\Omega^+$ and $\Omega^-$, and is capable of diffusing across the membrane interface $\Gamma$. In addition, $C$ can be adsorbed onto the interface, where it participates in surface chemical reactions; its surface concentration is denoted by $C_\Gamma$. Let $B_\Gamma$ represent the surface concentration of another species that exists exclusively on the interface $\Gamma$. These two interfacial species undergo a reversible chemical reaction described by
\begin{equation*}
      C_{\Gamma}+ B_{\Gamma} \overset{k_f}{\underset{k_r}{\rightleftharpoons}} A_{\Gamma}, 
\end{equation*}
where $A_\Gamma$ denotes the concentration of the reaction product, which also remains confined to the interface. The parameters $k_f$ and $k_r$ are the forward and reverse reaction rate constants, respectively.

The concentration field $C({\bm x}, t)$ is generally discontinuous across the interface and is defined piecewise as
\[
C(\bm{x}, t) = 
\begin{cases}
C^+(\bm{x}, t), & \text{for } (\bm{x}, t) \in \Omega^+ \times [0, T], \\
C^-(\bm{x}, t), & \text{for } (\bm{x}, t) \in \Omega^- \times [0, T].
\end{cases}
\]
The jump of $C$ across the interface $\Gamma$ is denoted by
\[
[\![ C(\bm{x}, t) ]\!] = \lim_{\Omega^+ \ni \bm{y} \to \bm{x}} C^+(\bm{y}, t) - \lim_{\Omega^- \ni \bm{y} \to \bm{x}} C^-(\bm{y}, t), \quad \forall (\bm{x}, t) \in \Gamma.
\]
For notational simplicity, this is abbreviated as
\[
[\![ C ]\!] = C^+ - C^-.
\]
The above notation and definitions also apply to other physical quantities, including scalars, vectors, and tensors, unless otherwise specified.
% In fact, we need to notice that $C$ is a piecewise continuous function 
% \begin{equation}
%   C({\bm x},t) = 
%   \begin{cases}
%     C^+({\bm x},t), & \mbox{if } ({\bm x},t)\in \Omega^+\times[0,T], \\
%     C^-({\bm x},t), & \mbox{if } ({\bm x},t)\in \Omega^-\times[0,T].
%   \end{cases}
% \end{equation}
% The jump of $C$ across the interface $\Gamma$ is defined on $\Gamma$ and denoted by 
% \begin{equation}
% [\![C({\bm x},t)]\!]=\lim_{\Omega^+\ni{\bm y}\rightarrow{\bm x}} C^+({\bm y},t)-\lim_{\Omega^-\ni{\bm y}\rightarrow{\bm x}}C^-({\bm y},t),\quad \forall ({\bm x},t) \in \Gamma.
% \end{equation}
% Here, the above equation is also simplified as 
% \begin{equation}
%   [\![C]\!] = C^+ - C^-. 
% \end{equation}
% Without loss of generality, the description about $C$ is also available for other variables, including scaler, vector and tensor.

The total free energy of the system is given by
\begin{equation}\label{eq:Total_Free_Energy}
    E_{\text{tot}} = E_{\text{kin}} + E_{\text{mix}} + E_{\Gamma},
\end{equation}
where the individual components represent the kinetic energy of the fluids, the bulk mixing energy, and the interfacial energy, respectively:
\begin{subequations}\label{eq:Different_Energy}
\begin{align}
    E_{\text{kin}} &= \frac{1}{2} \sum_{\pm} \int_{\Omega^\pm} \rho^\pm |\vel^\pm|^2 \, d\sL^{d}, \label{eq:Kinetic_Energy} \\
    E_{\text{mix}} &= \sum_{\pm} \int_{\Omega^\pm} f(C^\pm) \, d\sL^{d}, \label{eq:Mix_Energy} \\
    E_{\Gamma} &= \int_{\Gamma} \left( \gamma_0 + \sum_{K} g(K) \right) \, d\mathcal{H}^{d-1}.\label{eq:Interface_Energy}
\end{align}
\end{subequations}
Here, $\rho^\pm$ and ${\bm u}^\pm$ denote the density and velocity of the fluid in $\Omega^\pm$, $d\sL^{d}$ represents the Lebesgue measure in $\mathbb{R}^d$, and $d\mathcal{H}^l$ denotes the $l$-dimensional Hausdorff measure in $\mathbb{R}^d$. The term $E_{\text{mix}}$ accounts for the chemical energy of the solute in the bulk.
The bulk free energy density $f(C)$ is defined in $\Omega^\pm$ as
\begin{equation}\label{eq:Bulk_Energy_Density}
    f(C) = U_C C + RT C_\infty \left( \frac{C}{C_\infty} \ln \left( \frac{C}{C_\infty} \right) - \frac{C}{C_\infty} \right),
\end{equation}
where $U_C$ is the standard chemical potential in the bulk, $R = k_B N_A$ is the universal gas constant with $k_B$ and $N_A$ denoting the Boltzmann and Avogadro constants, respectively, $T$ is the absolute temperature, and $C_\infty$ is the maximum molar concentration of the solute in the bulk. It is straightforward to verify that
\[
f'(C) = U_C + RT \ln \left( \frac{C}{C_\infty} \right), \quad f''(C) = \frac{RT}{C} > 0.
\]

 $E_\Gamma$ captures the energy on the interface $\Gamma$, including a constant intrinsic surface energy density $\gamma_0$ and additional energy contributions $g(K)$ from the interfacial species $K \in \{A_\Gamma, B_\Gamma, C_\Gamma\}$  defined on $\Gamma$ by
\begin{equation}\label{eq:Interfacial_Energy_Density}
    g(K) = U_K K + \omega_K RTK_\infty \left( \frac{K}{K_\infty} \ln \left( \frac{K}{K_\infty} \right) - \frac{K}{K_\infty} \right),
\end{equation}
where $U_K$ is the standard chemical potential of species $K$ on the interface, $K_\infty$ denotes the maximum interfacial packing concentration and $\omega_K$ measures the adsorption intensity. For $A_\Gamma, B_\Gamma$ and $C_\Gamma$, we use $\omega_a, \omega_b$ and $\omega_c$ to represent the corresponding values of $\omega_K$, respectively. In this content, we assume that the following fact holds 
\begin{equation*}
    A_{\Gamma,\infty}=B_{\Gamma,\infty}=C_{\Gamma,\infty}.
\end{equation*}
Similarly, one has
\[
g'(K) = U_K + RT\omega_K \ln \left( \frac{K}{K_\infty} \right), \quad g''(K) = \frac{RT\omega_K}{K} > 0.
\]

We assume that the fluids are   incompressible, and their dynamics within the bulk domains $\Omega^\pm$ are governed by the incompressible Navier–Stokes equations:
\begin{subequations}\label{eq:NS_1}
\begin{align}
  \rho^\pm \partial_t^\bullet {\bm u}^\pm &= -\nabla p^\pm + \nabla \cdot {\bm \sigma}^\pm, \label{eq:NS_11} \\
  \nabla \cdot {\bm u}^\pm &= 0, \label{eq:NS_12}
\end{align}
\end{subequations}
where $p^\pm$ denote the pressure fields, and ${\bm \sigma}^\pm = \eta^\pm (\nabla {\bm u}^\pm + (\nabla {\bm u}^\pm)^T)$ are the viscous stress tensors with $\eta^\pm$ being the dynamic viscosity in $\Omega^\pm$. Here the material derivative is defined by
\begin{equation}\label{eq:Material_Derivative}
  \partial_t^\bullet = \partial_t + {\bm u} \cdot \nabla.
\end{equation}

We further assume that the velocity field is continuous across the  interface $\Gamma$, i.e.,
\begin{equation}\label{eq:Velocity_Continuity}
  [\![ {\bm u} ]\!] = {\bm 0}.
\end{equation}

To describe the interface evolution, we introduce a parameterization of $\Gamma(t)$ over a fixed reference surface $\Upsilon$:
\begin{equation}
  {\bm X}({\bm q}, t): \Upsilon \times [0, T] \to \mathbb{R}^d, \quad \text{for all } {\bm q} \in \Upsilon,
\end{equation}
and  the interface obeys the following kinematic condition:
\begin{equation}\label{eq:Interface_Evolution}
  {\bm V} \cdot {\bm n} = {\bm u} \cdot {\bm n},
\end{equation}
where the interface velocity ${\bm V}$ is defined as
\begin{equation}
  {\bm V}({\bm X}({\bm q}, t), t) = \partial_t {\bm X}({\bm q}, t), \quad \text{for all } {\bm q} \in \Upsilon.
\end{equation}

The mean curvature of the interface $\Gamma(t)$ is defined by $\kappa = \nabla_s \cdot \bn$, where $\nabla_s$ denotes the surface gradient operator. The curvature vector of the evolving interface $\Gamma(t)$ satisfies
\begin{equation}\label{eq:Curvature_Relation}
  \kappa {\bm n} = -\Delta_s {\bm X}.
\end{equation}
 Here $\Delta_s = \nabla_s \cdot \nabla_s$ is the Laplace–Beltrami operator on $\Gamma(t)$. %, and $\nabla_s$ denotes the surface gradient operator   defined as the tangential projection of the standard gradient onto the interface \(\Gamma\):
% \[
% \nabla_s = (\mathbb{I}- \bn \otimes \bn)\nabla,
% \]
% where \(\mathbb{I}\) is the identity tensor.  

According to the conservation law in domain \eqref{eq:Bulk_Mass_Conservation_Diff2}, the concentration  of the substance $C$ in the bulk domains $\Omega^\pm$ satisfies the following equations
\begin{equation}
  \label{eq:Bulk_Mass_Conservation_DiffC}
  \partial^\bullet_t C^\pm + \nabla\cdot{\bm J}_C^\pm = 0,
\end{equation} 
where ${\bm J}_C$ is the diffusion flux. 
% On the interface $\Gamma$, we assume 
% \begin{subequations}\label{eq:Source_and_Flux}
% \begin{align}
%   &S^+ + J_s+{\bm n}\cdot{\bm J}_C^+=0, \label{eq:Bulk_Flux_Condition_Plus}
%   \\
%   &S^- - J_s-{\bm n}\cdot{\bm J}_C^-=0. \label{eq:Bulk_Flux_Condition_Minus}
% \end{align}
% \end{subequations}

Similarly,   the conservation law on the interface \eqref{eq:Surface_Mass_Conservation_Diff1},  yields 
\begin{subequations}\label{eq:Surface_Mass_Conservation}
\begin{align} 
{\bm n}\cdot{\bm J}_C^\pm=&-J_s\mp S^{\pm}\label{eq:Bulk_Flux_Condition}\\
  \partial^\bullet_t A_\Gamma +\nabla_s\cdot{\bm u}A_\Gamma + \nabla_s\cdot{\bm J}_{A_\Gamma} &=S_{A_\Gamma}=~~ \mathcal{R},\label{eq:Surface_Mass_Conservation_DiffA}
  \\
  \partial^\bullet_t B_\Gamma +\nabla_s\cdot{\bm u}B_\Gamma + \nabla_s\cdot{\bm J}_{B_\Gamma} &=S_{B_\Gamma}= -\mathcal{R},\label{eq:Surface_Mass_Conservation_DiffB}
  \\
  \partial^\bullet_t C_\Gamma +\nabla_s\cdot{\bm u}C_\Gamma + \nabla_s\cdot{\bm J}_{C_\Gamma} &=S_{C_\Gamma}= -\mathcal{R} + \sum_\pm S^\pm,\label{eq:Surface_Mass_Conservation_DiffC}  
\end{align}
\end{subequations}
where $J_s$ is the transmembrane flux, \( S^\pm \) stem from the bulk–interface exchange, \( \mathcal{R} \) arises from the interfacial reaction, and \( {\bm J}_K \) represents the surface diffusion flux tangential to the interface \( \Gamma \). For convenience, we define the surface material derivative \( \partial_t^\bullet \) as  
\begin{equation}
\label{eq:Surface_Material_Derivative}
\partial_t^\bullet = \partial_t + {\bm u} \cdot \nabla_s.
\end{equation}

The above formulation establishes the governing equations of the coupled system. To complete the model, it remains to specify the constitutive expressions for the fluxes \( {\bm J}_C \) and \( {\bm J}_K \), the source terms \(J_s, S^\pm \) and \( \mathcal{R} \), appropriate interfacial conditions on \( \Gamma \), and boundary conditions on \( \partial \Omega \).

\subsection{Constitutive relations, interfacial conditions and boundary conditions}

In this subsection, we will make use of the principles of nonequilibrium thermodynamics to examine the expressions of ${\bm J}_C,{\bm J}_K,\mathcal{R},S^\pm$, the interfacial conditions on $\Gamma$ and the boundary conditions on $\partial\Omega$, which are consistent with the second law of thermodynamics.

Firstly, taking the derivative of the kinetic energy $E_{kin}$ with respect to time, we obtain
\begin{align}
  \frac{d E_{kin}}{dt} & = \sum_{\pm}\int_{\Omega^\pm} \rho^\pm \partial^\bullet_t{\bm u}^\pm \cdot\vel^\pm d\sL^{d}\nonumber\\
&= \sum_{\pm}\int_{\Omega^\pm} \vel^\pm\cdot(\nabla\cdot{\bm \sigma}^\pm) d\sL^{d}-\sum_{\pm}\int_{\Omega^\pm} \vel^\pm \cdot \nabla p^\pm d\sL^{d}\nonumber\\
& = \sum_{\pm}\int_{\Omega^\pm} (\nabla\cdot({\bm\sigma}^\pm\cdot\vel^\pm)-{\bm\sigma}^\pm:\nabla\vel^\pm)  d\sL^{d}- \sum_{\pm}\int_{\Omega^\pm} (\nabla\cdot(p^\pm\vel^\pm) - p^\pm\nabla\cdot\vel^\pm) d\sL^{d}\nonumber\\
& = \sum_{\pm}\int_{\Omega^\pm} \big[\nabla\cdot(\mathbb{T}^\pm \cdot \vel^\pm) -{\bm\sigma}^\pm:\nabla\vel^\pm\big]  d\sL^{d}\nonumber\\
&= -\sum_{\pm}\int_{\Omega^\pm} \frac{1}{2\eta^\pm} ||{\bm\sigma}^\pm||_F^2 d\sL^{d} - \int_{\Gamma}\bn\cdot [\![\mathbb{T}]\!]\cdot \vel d\mathcal{H}^{d-1} +\int_{\partial\Omega} \bn\cdot \mathbb{T}\cdot \vel d\mathcal{H}^{d-1}, \label{eq:Kinetic_Energy_Derivative}
\end{align}
where $\mathbb{T}^\pm = (-p^\pm\mathbb{I}+{\bm\sigma}^\pm)$ is the stress tensor with $\mathbb{I}_d\in\mathbb{R}^{d\times d}$ being the identity matrix and $||\cdot||_F$ denotes the Frobenius norm. In the above derivation, we have used the Reynolds transport formula and Eq. \eqref{eq:NS_1}.

Secondly, for the chemical energy $E_{mix}$, its time derivative satisfies
\begin{align}
  \frac{dE_{mix}}{dt} &= \int_{\Omega^+} f^{\prime\prime}(C^+)\nabla C^+\cdot{\bm J}_C^+ d\sL^{d} +\int_{\Gamma} f^{\prime}(C^+){\bm n}\cdot{\bm J}_C^+ d\mathcal{H}^{d-1}- \int_{\partial\Omega} f^{\prime}(C^+){\bm n}\cdot{\bm J}_C^+ d\mathcal{H}^{d-1} \nonumber\\
&+ \int_{\Omega^-} f^{\prime\prime}(C^-)\nabla C^-\cdot{\bm J}_C^+ d\sL^{d} -\int_{\Gamma} f^{\prime}(C^-){\bm n}\cdot{\bm J}_C^- d\mathcal{H}^{d-1},\label{eq:Mix_Energy_Derivative}
\end{align}
where we have used Eq. \eqref{eq:Bulk_Energy_Derivative2}.

Finally, for the interfacial energy $E_\Gamma$, we can compute its time derivative as follows 
\begin{align}
  \frac{dE_{\Gamma}}{dt} & = \int_{\Gamma}\gamma_0 \kappa{\bm n}\cdot{\bm u} d\mathcal{H}^{d-1} \nonumber\\
                         &+ \sum_{K}\int_{\Gamma} \big\{{\bm u}\cdot(\tilde{\gamma}(K)\kappa{\bm n}-\nabla_s\tilde{\gamma}(K)) + g''(K){\bm J}_{K}\cdot\nabla_s K + g^{\prime}(K)S_K \big\}d\mathcal{H}^{d-1}\nonumber\\
                         &=\int_{\Gamma}{\bm u}\cdot(\gamma\kappa{\bm n}-\nabla_s\gamma) d\mathcal{H}^{d-1}+\sum_{K}\int_{\Gamma} \big\{ g^{\prime\prime}(K){\bm J}_{K}\cdot\nabla_s K + g^{\prime}(K)S_K \big\}d\mathcal{H}^{d-1},\label{eq:Interface_Energy_Derivative}
\end{align}
where $\tilde{\gamma}(K) = -\omega_K RTK$, $\gamma=\gamma_0+\sum_{K}\tilde{\gamma}(K)$ and we have used Eq. \eqref{eq:Surface_Energy_Derivative2}.

Combining Eqs. \eqref{eq:Kinetic_Energy_Derivative}, \eqref{eq:Mix_Energy_Derivative} and \eqref{eq:Interface_Energy_Derivative}, yields
\begin{align}
  \frac{dE_{tot}}{dt} &= \frac{dE_{kin}}{dt} + \frac{dE_{mix}}{dt} + \frac{dE_{\Gamma}}{dt} \nonumber\\
                      &= -\sum_{\pm}\int_{\Omega^\pm} \frac{1}{2\eta^\pm} ||{\bm\sigma}^\pm||_F^2 d\sL^{d} + \sum_{\pm}\int_{\Omega^\pm} f^{\prime\prime}(C^\pm)\nabla C^\pm\cdot{\bm J}_C^\pm d\sL^{d} \nonumber\\
                      &+ \int_{\Gamma}\big[(\gamma\kappa{\bm n}-\nabla_s\gamma)-{\bn}\cdot [\![\mathbb{T}]\!]\big]\cdot{\bm u} d\mathcal{H}^{d-1} + \sum_{K}\int_{\Gamma} g^{\prime\prime}(K){\bm J}_{K}\cdot\nabla_s K d\mathcal{H}^{d-1}\nonumber \\
                      &+ \sum_{\pm}\int_{\Gamma} (g^{\prime}(C_\Gamma)-f^{\prime}(C^\pm)) S^\pm d\mathcal{H}^{d-1} + \int_{\Gamma}[g^{\prime}(A_\Gamma)-g^{\prime}(B_\Gamma)-g^{\prime}(C_\Gamma)]\mathcal{R}d\mathcal{H}^{d-1} \nonumber \\
                    %  &+\int_{\Gamma}[S^+ + J_s+{\bm n}\cdot{\bm J}_C^+]f(C^+)d\mathcal{H}^{d-1} + \int_{\Gamma}[S^- - J_s-{\bm n}\cdot{\bm J}_C^-]f(C^-)d\mathcal{H}^{d-1}\nonumber \\
                      &-\int_{\Gamma}[\![f^{\prime}(C)]\!]J_sd\mathcal{H}^{d-1}+\int_{\partial\Omega} {\bn}\cdot \mathbb{T}\cdot \vel d\mathcal{H}^{d-1}- \int_{\partial\Omega} f^{\prime}(C^+){\bm n}\cdot{\bm J}_C^+ d\mathcal{H}^{d-1}.\label{eq:Total_Free_Energy_Derivative}
\end{align}

Based on the second law of thermodynamics, we choose the constitutive relations, the interfacial conditions, and the boundary conditions such that the total free energy has non-positive dissipation. Here, we require that each term on the right-hand side of Eq. \eqref{eq:Total_Free_Energy_Derivative} is non-positive.  

\textbf{\textsl{ Bulk region}}: It is obvious that the first term is the viscous dissipation in the bulk and non-positive. For the second term, which stems from the bulk diffusion, we set
\begin{equation}
  \label{eq:Bulk_Flux}
  {\bm J}^\pm_C = -D_C^\pm \nabla C^\pm,
\end{equation}
where $D_C^\pm$ is the diffusion coefficient for the bulk substance.

\textbf{\textsl{Interface}}: By vanishing the third term, the Laplace-Young condition, which describes the balance of the stress jump of the fluids and the surface force, is derived
\begin{equation}
  \label{eq:Laplace_Young_Condition}
  {\bn}\cdot [\![\mathbb{T}]\!] = (\gamma\kappa{\bm n}-\nabla_s\gamma).
\end{equation}
Similarly to the second term, the fourth term describes the dissipation that arises from the surface diffusion, so the following expression can be given
\begin{equation}
  \label{eq:Surface_Flux}
  {\bm J}_K = -D_K \nabla_s K,
\end{equation} 
where $D_K$ is the diffusion coefficient for the surface substance. 

By applying the result of Lemma~\ref{lma:Adsorption_Desorption_Source}, we obtain the following expression for the adsorption-desorption source term from the fifth term:
\begin{equation}
  \label{eq:Adsorption_Desorption_Source}
  S^\pm = k_{\text{ad}}^\pm \frac{C^\pm}{C_\infty} - k_d^\pm \left(\frac{C_\Gamma}{C_{\Gamma,\infty}}\right)^{\omega_c},
\end{equation}
where \(S^\pm\) represents the net mass flux between the bulk region \(\Omega^\pm\) and the interface \(\Gamma\), arising from adsorption and desorption processes. The first term corresponds to the adsorption from the bulk onto the interface, while the second term accounts for the desorption back into the bulk.

Similarly, from the sixth term, it indicates that the source term in Eq. \eqref{eq:Surface_Mass_Conservation_DiffA} satisfies
\begin{equation}
  \label{eq:Chemical_Reaction_Source}
  \mathcal{R} = k_f\left(\frac{B_{\Gamma}}{B_{\Gamma,\infty}}\right)^{\omega_b}\left(\frac{C_{\Gamma}}{C_{\Gamma,\infty}}\right)^{\omega_c} - k_r\left(\frac{A_{\Gamma}}{A_{\Gamma,\infty}}\right)^{\omega_a},
\end{equation}
where the result of Lemma \ref{lma:Chemical_Reaction_Source} has been used.

% Next, we consider that the dissipation terms, which range from the seventh term to the eighth term, uniformly vanish. 
% Then, we directly have on the interface $\Gamma$
% \begin{subequations}\label{eq:Source_and_Flux}
% \begin{align}
%   &S^+ + J_s+{\bm n}\cdot{\bm J}_C^+=0, \label{eq:Bulk_Flux_Condition_Plus}
%   \\
%   &S^- - J_s-{\bm n}\cdot{\bm J}_C^-=0. \label{eq:Bulk_Flux_Condition_Minus}
% \end{align}
% \end{subequations}

For the seventh term, Lemma~\ref{lma:Bulk_Concentration_Jump} provides the following expression :
\begin{equation}
  \label{eq:Bulk_Concentration_Jump}
  J_s = k_c \frac{[\![ C ]\!]}{C_\infty},
\end{equation}
where  \(k_c\) is a transport coefficient characterizing the permeability of the membrane.

\textbf{\textsl{The outer boundary}}: Let $\partial\Omega_1$ denote the left and right boundaries of $\Omega$, and let $\partial\Omega_2$ be the bottom and top boundaries of $\Omega$. Consequently, the outer boundary of $\Omega$ can be expressed as $\partial\Omega=\partial\Omega_1\cup\partial\Omega_2$. Let us suppose that the velocity satisfies
\begin{subequations}\label{eq:Velocity_Outer_Condition}
\begin{align}
  \vel &= {\bm 0}, \quad\text{on}~\partial\Omega_1,
  \\
  {\bn}\cdot \mathbb{T} &= {\bm 0}, \quad\text{on}~\partial\Omega_2,
\end{align}
\end{subequations}
and the concentration of the bulk substance obeys
\begin{equation}
  \label{eq:Bulk_Concentration_Outer_Condition}
  {\bm n}\cdot{\bm J}_C^+ = 0,
\end{equation}
such that there is no dissipation on the outer boundary $\partial\Omega$. 

% \begin{rmk}
%   In fact, we can also directly apply the linear response to the fifth term, i.e., 
%   \begin{equation}
%     \label{eq:Adsorption_Desorption_Source_Non}
%     S^\pm = -k_d^\pm(g^{\prime}(C_\Gamma)-f^{\prime}(C^\pm)).
%   \end{equation}
%   If we carry out the Taylor expansion for $\ln(x)$ at $x=1$, Eq. \eqref{eq:Adsorption_Desorption_Source_Non} can result in  
%   \begin{equation}
%     S^\pm \approx -k_d^\pm\bigg(\lambda_a\frac{C^\pm}{C_\infty} - \left(\frac{C_\Gamma}{C_{\Gamma,\infty}}\right)^{\omega_c}\bigg).
%   \end{equation}
%   This indicates that Eq. \eqref{eq:Adsorption_Desorption_Source} is the linear approximation of Eq. \eqref{eq:Adsorption_Desorption_Source_Non}. Similar results hold for the seventh and eighth terms.
% \end{rmk}

Using the fluxs \eqref{eq:Bulk_Flux} and \eqref{eq:Surface_Flux}, the Laplace-Young condition \eqref{eq:Laplace_Young_Condition}, the source terms \eqref{eq:Adsorption_Desorption_Source} and \eqref{eq:Chemical_Reaction_Source}, the interface conditions \eqref{eq:Bulk_Concentration_Jump} and \eqref{eq:Bulk_Flux_Condition}, and the outer boundary conditions \eqref{eq:Velocity_Outer_Condition} and \eqref{eq:Bulk_Concentration_Outer_Condition}, the energy law \eqref{eq:Kinetic_Energy_Derivative} reduces to 
\begin{align*}
  \frac{dE_{tot}}{dt} &= -\sum_{\pm}\int_{\Omega^\pm} \frac{1}{2\eta^\pm} ||{\bm\sigma}^\pm||_F^2 d\sL^{d} 
                       - \sum_{\pm}\int_{\Omega^\pm} D_C^\pm f^{\prime\prime}(C^\pm)|\nabla C^\pm|^2 d\sL^{d} \nonumber\\
                      &- \sum_{K}\int_{\Gamma} D_{K}g^{\prime\prime}(K)|\nabla_s K|^2 d\mathcal{H}^{d-1}
                       + \sum_{\pm}\int_{\Gamma} (g^{\prime}(C_\Gamma)-f^{\prime}(C^\pm)) S^\pm d\mathcal{H}^{d-1}\nonumber \\
                      &+ \int_{\Gamma}[g^{\prime}(A_\Gamma)-g^{\prime}(B_\Gamma)-g^{\prime}(C_\Gamma)]\mathcal{R}d\mathcal{H}^{d-1}
                       -\int_{\Gamma}[\![f^{\prime}(C)]\!]J_s d\mathcal{H}^{d-1} \leq 0,
\end{align*}
where the following facts are used
\begin{align*}
  &f^{\prime\prime}(C^\pm) >0,\quad g^{\prime\prime}(K)> 0,\quad (g^{\prime}(C_\Gamma)-f^{\prime}(C^\pm)) S^\pm \leq 0,  
  \\
  &[g^{\prime}(A_\Gamma)-g^{\prime}(B_\Gamma)-g^{\prime}(C_\Gamma)]\mathcal{R} \leq 0,\quad [\![f^{\prime}(C)]\!]J_s \geq 0.
\end{align*}

\subsection{Dimensionless}\label{sec:Dimensionless_Model}
The physical variables are scaled by the corresponding characteristic quantities
\begin{align*}
  & \hat{\rho} = \frac{\rho}{\rho^-},\quad \hat{\eta} = \frac{\eta}{\eta^-},\quad \hat{k_c} = \frac{k_c}{k_c}, \quad \hat{k}_r = \frac{k_r}{k_c}, \quad \hat{k}_d = \frac{k_d}{k_c}, \quad \hat{\bm x} = \frac{{\bm x}}{L},
  \\
  &\hat{\bm u} = \frac{{\bm u}}{U}, \quad \hat{C} = \frac{C}{C_\infty}, \quad \hat{K} =\frac{K}{K_{\Gamma,\infty}},\quad \hat{U}_C = \frac{U_C}{RT},\quad \hat{U}_K := \frac{U_K}{RT}, 
  \\
  &\hat{D}_C := \frac{D_C}{D_C^-}, \quad \hat{D}_K := \frac{D_K}{D_{C_\Gamma}},\quad \hat{t} := \frac{Ut}{L},\quad \hat{p} := \frac{p}{\rho^-U^2},\quad \hat{\gamma}:=\frac{\gamma}{\gamma_0}, 
\end{align*}
where $L$ and $U$ are the characteristic length and velocity, respectively. We still use the same notations for the same variables after the non-dimensionalization. 

Therefore, the relevant dimensionless numbers are the Reynolds number $Re$, the capillary number $Ca$, the Weber number $We$, the Biot number $Bi$, the adsorption depth $Da$, the bulk Peclet number $Pe$, the surface Peclet number $Pe_\Gamma$, the elasticity $E$, the equilibrium constant $\lambda_a$ that corresponds to the adsorption, and the equilibrium constant $\lambda_c$ that corresponds to the chemical reaction. These quantities are defined as follows.
\begin{align*}
  & Re = \frac{\rho^-UL}{\eta^-}, \quad Ca = \frac{\eta^-U}{\gamma_0},\quad We = \frac{\rho^- U^2 L}{\gamma_0},\quad Bi = \frac{k_c L}{U C_{\Gamma,\infty}}, \quad  Da = \frac{C_{\Gamma,\infty}}{L C_\infty}, \quad Pe = \frac{UL}{D_C},
  \\
  & Pe_\Gamma = \frac{UL}{D_{C_\Gamma}}, \quad E = \frac{RTC_{\Gamma,\infty}}{\gamma_0}, \quad  \lambda_a = \exp\bigg(\frac{U_C-U_{C_\Gamma}}{RT}\bigg),\quad \lambda_c = \exp\bigg(\frac{U_{B_\Gamma}+U_{C_\Gamma}-U_{A_\Gamma}}{RT}\bigg).
\end{align*}

The dimensionless system is summarized as follows: in the bulk region $\Omega^{\pm}$
\begin{subequations} \label{eq:bulkeq}
\begin{align}
  \rho^\pm\partial^\bullet_t {\bm u}^\pm &= -\nabla p^\pm + \frac{1}{Re}\nabla\cdot {\bm \sigma}^\pm, \label{eq:NS_11_Ndim}
  \\
  \nabla\cdot {\bm u}^\pm &= 0,\label{eq:NS_12_Ndim}\\
  \partial^\bullet_t C^\pm &= \frac{1}{Pe}\nabla\cdot(D_C^\pm\nabla C^\pm),\label{eq:Bulk_Mass_Conservation_DiffC_Ndim}
\end{align}
\end{subequations}
% The dimensionless Navier–Stokes equations become
% \begin{subequations}\label{eq:NS_1_Ndim}
% \begin{align}
%   \rho^\pm\partial^\bullet_t {\bm u}^\pm &= -\nabla p^\pm + \frac{1}{Re}\nabla\cdot {\bm \sigma}^\pm, \label{eq:NS_11_Ndim}
%   \\
%   \nabla\cdot {\bm u}^\pm &= 0,\label{eq:NS_12_Ndim}
% \end{align}
% \end{subequations}
on the interface $\Gamma$
\begin{subequations}\label{eq:Interfacial_Condition_Ndim}
\begin{align}
  {\bn}\cdot [\![\mathbb{T}]\!] &= \frac{1}{We}(\gamma\kappa{\bm n}-\nabla_s\gamma),\label{eq:Laplace_Young_Condition_Ndim}
  \\
  [\![{\bm u}]\!] &= {\bm 0},\label{eq:Velocity_Continuity_Ndim}
  \\
  {\bm X}_t\cdot{\bm n} &= {\bm u}\cdot{\bm n},\label{eq:Interface_Evolution_Ndim}
  \\ 
  \kappa{\bm n} &= -\Delta_s{\bm X},\label{eq:Curvature_Relation_Ndim}\\
  \frac{D_C^\pm}{Da\cdot Pe}{\bm n}\cdot \nabla C^\pm &=   J_s \pm S^\pm, \label{eq:Bulk_Flux_Condition_Ndim}\\
    \partial^\bullet_t A_\Gamma +\nabla_s\cdot{\bm u}A_\Gamma &= \frac{1}{Pe_\Gamma}\nabla_s\cdot(D_{A_\Gamma}\nabla_s A_\Gamma) + \mathcal{R},\label{eq:Surface_Mass_Conservation_DiffA_Ndim}
  \\
  \partial^\bullet_t B_\Gamma +\nabla_s\cdot{\bm u}B_\Gamma &= \frac{1}{Pe_\Gamma}\nabla_s\cdot(D_{B_\Gamma}\nabla_s B_\Gamma) - \mathcal{R},\label{eq:Surface_Mass_Conservation_DiffB_Ndim}
  \\
  \partial^\bullet_t C_\Gamma +\nabla_s\cdot{\bm u}C_\Gamma &= \frac{1}{Pe_\Gamma}\nabla_s\cdot(D_{C_\Gamma}\nabla_s C_\Gamma) +  \sum_{\pm}S^\pm - \mathcal{R},\label{eq:Surface_Mass_Conservation_DiffC_Ndim}
\end{align}
\end{subequations}
and the outer boundary $\partial\Omega$
\begin{subequations}\label{eq:outbd}
    \begin{align}
        &{\bm u} = {\bm 0}\quad\text{on}~\partial\Omega_1,\quad \mathbb{T}\cdot{\bm n} = 0\quad\text{on}~\partial\Omega_2,\label{eq:Velocity_Outer_Condition_Ndim}
        \\
        &{\bm n}\cdot \nabla C^+ = 0\quad\mbox{on}~\partial\Omega,\label{eq:Bulk_Concentration_Outer_Condition_Ndim}
    \end{align}
\end{subequations}
 where the interfacial tension reads
\begin{equation}
  \label{eq:gamma}
  \gamma = 1 - E\sum_{K}\omega_K K,  
\end{equation}
reaction rate function and leaky flux are defined as
\begin{subequations}\label{eq:Source_Ndim}
\begin{align}
& J_s = Bi [\![C]\!]\\
  &\mathcal{R} = Bi(k_f B_\Gamma^{\omega_b} C_\Gamma^{\omega_c}- k_r A_\Gamma^{\omega_a}),\label{eq:Chemical_Reaction_Source_Ndim}
  \\
  &S^\pm = Bi(k_{ad}^\pm C^\pm - k_d^\pm C_\Gamma^{\omega_c}).\label{eq:Adsorption_Desorption_Source_Ndim}
\end{align}
\end{subequations}

\begin{thm}
Consider the system \eqref{eq:bulkeq} with the corresponding interfacial and boundary conditions \eqref{eq:Interfacial_Condition_Ndim}-\eqref{eq:outbd}. Then the following properties hold:
   \begin{itemize}

       \item  \textbf{Mass Conservation Laws:}
\begin{subequations}\label{eq:Total_Mass_Conservation}
\begin{align}
  &  \mathrm{\frac{d}{dt}}m_s(A_\Gamma,B_\Gamma,0,0;t) = 0, \label{eq:Total_Mass_Conservation1}
  \\
  &  \mathrm{\frac{d}{dt}}m_s(A_\Gamma,0,C_\Gamma,C;t) = 0, \label{eq:Total_Mass_Conservation2}
\end{align} 
\end{subequations}
where 
\begin{equation}\label{eq:massdefinition}
  m_s(A_\Gamma,B_\Gamma,C_\Gamma,C;t) = \sum_K\int_{\Gamma} K d\mathcal{H}^{d-1} + \frac{1}{Da}\sum_\pm\int_{\Omega^\pm} C^\pm d\sL^{d}.
\end{equation}
       \item \textbf{Energy Dissipation Law:} The system satisfies the following energy dissipation identity:
     \begin{equation}\label{eq:Energy_Dissipation_Law_Ndim}
         \frac{d E_{tot}}{dt} = -(\Delta_k + \Delta_C + \Delta_K + \Delta_a + \Delta_c +\Delta_m)\le 0,
     \end{equation}  
where 
\begin{equation*}
  E_{tot} = \frac 1 2 \sum_{\pm}\int_{\Omega^\pm}\rho^\pm |\vel^\pm|^2 d\sL^{d} + \frac{1}{We\cdot Da}\sum_{\pm}\int_{\Omega^\pm} f(C^\pm) d\sL^{d} + \frac{1}{We}\int_{\Gamma} \big(1 + \sum_{K}g(K)\big)d\mathcal{H}^{d-1},
\end{equation*}
and 
\begin{align*}
  &\Delta_k =\frac{1}{Re}\sum_{\pm}\int_{\Omega^\pm}\frac{1}{2\eta^\pm}||{\bm \sigma}^\pm||_F^2 d\sL^{d},~~\Delta_C = \frac{1}{We\cdot Da\cdot Pe} \sum_\pm \int_{\Omega^\pm} D_C^\pm f^{\prime\prime}(C^\pm) |\nabla C^\pm|^2 d\sL^{d},
  \\
  &\Delta_K = \frac{1}{We\cdot Pe_\Gamma} \sum_{K}\int_\Gamma D_K g^{\prime\prime}(K) |\nabla_s K|^2 d\mathcal{H}^{d-1},~~\Delta_a = \frac{1}{We}\sum_\pm\int_\Gamma (f^{\prime}(C^\pm) - g^{\prime}(C_\Gamma)) S^\pm d\mathcal{H}^{d-1},
  \\
  &\Delta_c = \frac{1}{We}\int_{\Gamma}[g^{\prime}(B_\Gamma)+g^{\prime}(C_\Gamma)-g^{\prime}(A_\Gamma)]\mathcal{R} d\mathcal{H}^{d-1},~~  \Delta_m = \frac{1}{We}\int_{\Gamma}[\![f^{\prime}(C)]\!]J_s d\mathcal{H}^{d-1},
\end{align*}
Here the six terms are the viscous dissipation in the bulk fluid $\Delta_k$, the dissipation arising from the diffusion of the bulk substance $\Delta_C$, the dissipation stemming from the diffusion of the interfacial substances $\Delta_K$, the dissipation due to substance adsorption and desorption $\Delta_a$, the dissipation originating from the chemical reaction of the interfacial substances $\Delta_c$, and the dissipation coming from the mass transfer across the interface $\Delta_m$, respectively. $f$ and $g$ are given respectively as 
\begin{align*}
  f(C) = E[U_C  + (\ln C - 1)]C,\quad g(K) = E[U_K + \omega_K(\ln (K) - 1)]K.
\end{align*}
Their derivatives are calculated as 
\begin{align*}
  &f^{\prime}(C) = E[U_C  + \ln C], \qquad\qquad~  f^{\prime\prime}(C) = \frac{E}{C},
  \\
  &g^{\prime}(K) = E[U_K + \omega_K\ln (K)], \quad  g^{\prime\prime}(C) = \frac{\omega_K E}{K}.
\end{align*}
   \end{itemize}
\end{thm}

\begin{rmk}
When there is no chemical reaction on the interface, i.e., \( \mathcal{R} = 0 \), the proposed model reduces to the classical soluble surfactant formulation with Henry-type surface energy, as studied in~\cite{zhao2021thermodynamically}. In the absence of a bulk solvent or interfacial absorption, the model becomes consistent with the classical formulation for insoluble surfactants~\cite{zhang2014derivation}. Furthermore, if both the chemical reaction \( \mathcal{R} = 0 \) and interfacial absorption \( S^\pm = 0 \), the model simplifies to the formulation for mass transport across a permeable interface as described in~\cite{qin2022phase}. Finally, by setting \( k_c = 0 \), the model degenerates to the standard fluid–structure interaction problem with a non-permeable membrane.

\end{rmk}

\begin{rmk}
In this work, we only consider the passive transmembrane flux \( J_s \) driven by the difference in chemical potential across the interface, corresponding to passive channels in biological membranes. In addition to passive transport, biological membranes also contain active pumps, such as \(\mathrm{Na/K}\) -ATPase, which transport ions against their concentration gradients, from regions of lower to higher concentration, by consuming metabolic energy (e.g. ATP).

To incorporate active pumping into the model, one may extend the interfacial flux condition as
\[
\frac{D_C^\pm}{Da \cdot Pe} {\bm n} \cdot \nabla C^\pm = J_s - J_{\mathrm{pump}} \pm S^\pm,
\]
where \( J_{\mathrm{pump}} \) represents an oriented active flux from \( \Omega^- \) to \( \Omega^+ \). A commonly used expression for \( J_{\mathrm{pump}} \) is the modified Michaelis–Menten form \cite{zaheri2020comprehensive}:
\[
J_{\mathrm{pump}} = J_{\max} \left( \frac{C^-}{K_M + C^-} \right)^{\beta},
\]
where \( J_{\max} \) is the maximal pumping rate, \( K_M \) is the Michaelis constant, and \( \beta \) characterizes the binding cooperativity. This extension allows for the modeling of energy-consuming transport processes essential for maintaining physiological membrane potential and ionic homeostasis.
\end{rmk}

\section{Weak Formulation}\label{sec:Weak_Formulation}
In this section, we propose an ALE weak formulation for the whole system by considering it into two interacting subsystem. The first subsystem is for the fluid dynamics, which is given by  Eqs.~\eqref{eq:NS_11_Ndim}--\eqref{eq:NS_12_Ndim} with boundary conditions \eqref{eq:Laplace_Young_Condition_Ndim}--\eqref{eq:Curvature_Relation_Ndim} and \eqref{eq:Velocity_Outer_Condition_Ndim}. The second subsystem describes the dynamics of the solute species in the bulk and on the interface. These include~\eqref{eq:Bulk_Mass_Conservation_DiffC_Ndim}, \eqref{eq:Bulk_Flux_Condition_Ndim}--\eqref{eq:Surface_Mass_Conservation_DiffC_Ndim} together with the associated interfacial and outer boundary conditions. 

Let $\mathcal{O}\subset\mathbb{R}^d$ be the ALE reference domain of $\Omega$ and let $\{\mathcal{\bm A}[t]\}_{t\in[0,T]}$ be a family of ALE mappings 
\begin{equation}\label{eq:ALEmap}
  \mathcal{\bm A}[t]:\mathcal{O}\rightarrow\Omega, \quad {\bm y}\mapsto \mathcal{\bm A}[t]({\bm y})={\bm x}({\bm y},t)\quad \text{for all} \quad t\in[0,T] \quad \text{and} \quad {\bm y}\in\mathcal{O},
\end{equation}
with $\Upsilon\subset\mathcal{O}$ and $\mathcal{O}^\pm\subset\mathcal{O}$ such that $\Gamma(t)={\bm X}(\Upsilon,t)$ and $\mathcal{\bm A}[t](\mathcal{O}^\pm)=\Omega^\pm(t)$, respectively. We further assume that  $\mathcal{\bm A}[t]\in[W^{1,\infty}(\mathcal{O})]^d$ and $\mathcal{\bm A}[t]^{-1}\in[W^{1,\infty}(\Omega)]^d$, and introduce the domain mesh velocity
\begin{equation}
  {\bm w}({\bm x},t):=\frac{\partial{\bm x}({\bm y},t)}{\partial t}\bigg{|}_{{\bm y}=\mathcal{\bm A}[t]^{-1}({\bm x})}\quad \text{for all}\quad t\in[0,T] \quad \text{and}  \quad {\bm x}\in\Omega,
\end{equation}
where we allow $\mathcal{\bm A}[t]$ to be somehow arbitrary, except that on the boundary it satisfies
\begin{equation}
  \label{eq:ALE_Restraint}
  [{\bm w}({\bm x},t)-{\bm u}({\bm x},t)]\cdot{\bm n} = 0\quad\text{on}~\Gamma(t)\cup\partial\Omega.
\end{equation}
We notice that this condition for the mesh velocity is crucial for geometric consistency so that the interface or boundary mesh follows the evolving interface or boundary.  The precise construction of the ALE mappings will be left later in Section \ref{sec:Finite_Element_Approximations}.

For a given scalar field $\phi : \Omega \times [0, T] \rightarrow \mathbb{R}$, we define the time derivative with respect to the ALE  reference frame as
\begin{equation}\label{eq:ALE_Derivative}
\partial_t^\circ \phi = \partial_t \phi + (\bm w \cdot \nabla) \phi.
\end{equation}
 This leads to the relation with the material derivative:
\begin{equation}\label{eq:Relation_Material_ALE_Bulk}
  \partial_t^\bullet \phi - \partial_t^\circ \phi = ((\bm{u} - \bm{w}) \cdot \nabla) \phi,
\end{equation}
on recalling \eqref{eq:Material_Derivative}.  This relation is essential for converting governing equations written in a material (Lagrangian) frame into the ALE frame. It plays a key role in deriving consistent weak formulations for both fluid dynamics and substance transport within deforming domains.

Analogously, on the moving interface $\Gamma(t)$, for a surface quantity $\phi_\Gamma : \Gamma(t) \times [0, T] \rightarrow \mathbb{R}$, the material derivative relates to the ALE derivative via
\begin{equation} \label{eq:Relation_Material_ALE_Surface}
\partial_t^\bullet \phi_\Gamma - \partial_t^\circ \phi_\Gamma = \left( (\mathbf{u} - \bm w_\Gamma) \cdot \nabla_s \right) \phi_\Gamma,
\end{equation}
where $\bm w_\Gamma = \bm w|{\Gamma}$ is the restriction of the mesh velocity to the interface, and $\nabla_s$ denotes the surface gradient operator.

To simplify notation in what follows, we introduce the piecewise-defined density and viscosity fields as
\begin{equation*}
\rho := \rho^+ \chi_{\Omega^+} + \rho^- \chi_{\Omega^-}, \quad \eta := \eta^+ \chi_{\Omega^+} + \eta^- \chi_{\Omega^-},\label{eq:rhoetaf}
\end{equation*}
where $\chi_{\Omega^\pm}$ denotes the characteristic function of the domain $\Omega^\pm$. We next aim to derive the weak formulation for the fluid dynamics and substance concentrations, in which we need to transfer the material derivative to the derivative with respect to the ALE frame.

\subsection{ALE weak formulation for fluid dynamics}
We follow the work in \cite{garcke2023structure} and define the following function spaces 
\begin{align*}
  & \mathbb{U}:= \big\{ {\bm\chi}\in[H^1(\Omega)]^d: {\bm\chi}=0~\text{on}~\partial\Omega_1 \big\},
  \\
  & \mathbb{P}:= \big\{ \chi\in L^2(\Omega): (\chi,1) = 0 \big\},\quad \mathbb{V} := H^1(0,T;[L^2(\Omega)]^d)\cap L^2(0,T;\mathbb{U}),
\end{align*}
where we denote by $(\cdot,\cdot)$ the $L^2$-inner product on $\Omega$.

It follows from the divergence free condition \eqref{eq:NS_12_Ndim}, the constraint on the mesh velocity \eqref{eq:ALE_Restraint} and \eqref{eq:Relation_Material_ALE_Bulk} that  \cite[(4.7)]{garcke2023structure}
\begin{equation}
  \label{eq:Inertia_Term_Weak_Form3}
  \ipd{\rho\partial_t^\bullet{\bm u},{\bm v}} = \ipd{\rho\partial_t^\circ{\bm u},{\bm v}} + \mathscr{A}(\rho,{\bm u}-{\bm w};{\bm u},{\bm v}) + \frac{1}{2}\ipd{\rho\nabla\cdot{\bm w},{\bm u}\cdot{\bm v}},\quad\forall{\bf v}\in [H^1(\Omega)]^d,
\end{equation}
where the antisymmetric operator $\mathscr{A}$ is defined as
\begin{equation}
  \mathscr{A}(\rho,{\bm g};{\bm u},{\bm v}) = \frac{1}{2}\left[\ipd{\rho({\bm g}\cdot\nabla){\bm u},{\bm v}} - \ipd{\rho({\bm g}\cdot\nabla){\bm v},{\bm u}}\right].
\end{equation}

Let $\ipd{\cdot,\cdot}_\Gt$ be the $L^2$-inner product over $\Gt$. Then for the viscous term in Eq.~\eqref{eq:NS_11_Ndim}, we take the inner product with ${\bm v}\in\mathbb{V}$, integrate by parts to get
\begin{align}
  \label{eq:Viscous_Term_Weak_Form}
  \ipd{\nabla\cdot\mathbb{T},{\bm v}} &= \int_{\Omega^-(t)} {\bm v}\cdot(\nabla\cdot\mathbb{T}^-)d\sL^{d} + \int_{\Omega^+(t)} {\bm v}\cdot(\nabla\cdot\mathbb{T}^+)d\sL^{d} \nonumber
  \\
  &= \ipd{p,\nabla\cdot{\bm v}} -\frac{2}{Re}\ipd{\eta D({\bm u}), D({\bm  v})} - \frac{1}{We}\ipd{\gamma\kappa{\bm n}-\nabla_s\gamma,{\bm v}}_{\Gamma(t)}, 
\end{align}
where $D(\bm u) = \frac{1}{2}\left(\nabla \bm u+(\nabla\bm u)^T\right)$, and we used \eqref{eq:Laplace_Young_Condition_Ndim}, \eqref{eq:Velocity_Continuity_Ndim} and \eqref{eq:Velocity_Outer_Condition_Ndim}.  

For the curvature term in \eqref{eq:Curvature_Relation_Ndim}, multiplying it by a test function ${\bm g}\in[H^1(\Gamma(t))]^d$ and integrating it over $\Gamma(t)$ followed by integration by parts yields
\begin{equation}
  \label{eq:Curvature_Relation_Ndim_Weak_Form1}
\ipd{\kappa{\bm n},~{\bm g}}_{\Gamma(t)} -\ipd{\nabla_s{\bm X},~\nabla_s{\bm g}}_{\Gamma(t)}=0.
\end{equation}

Collecting these results, it is natural to consider the following ALE weak formulation for the fluid dynamics. Given $\Gamma(0) = \Gamma_0$ and ${\bm u}(\cdot,0)={\bm u}_0$, we find $\Gamma(t)={\bm X}(\Upsilon,t)$ with interface velocity ${\bf V}\in [H^1(\Gamma(t))]^d,{\bm u}\in\mathbb{V},~p\in L^2(0,T;\mathbb{P})$ and $\kappa(\cdot,t)\in L^2(\Gamma(t))$ such that for all $t\in (0,T]$
\begin{subequations}
\label{eq:Two_Phase_Flow_WK}
\begin{align}
  &\ipd{\rho\partial_t^\circ{\bm u},~{\bm v}} + \frac{1}{2}\ipd{\rho\,\nabla\cdot{\bm w},{\bm u}\cdot{\bm v}} + \mathscr{A}(\rho,{\bm u}-{\bm w};{\bm u},{\bm v}) - \ipd{p,\nabla\cdot{\bm v}} \nonumber 
  \\
  &\qquad + \frac{2}{Re}\ipd{\eta D({\bm u}), D({\bm v})} + \frac{1}{We}\ipd{\gamma\,\kappa\,{\bm n}-\nabla_s\gamma,~{\bm v}}_{\Gamma(t)}=0\quad \forall {\bm v}\in\mathbb{V},\label{eq:NS_11_Ndim_WK}
  \\[0.6em]
  &\ipd{\nabla\cdot{\bm u},~q} =0 \qquad \forall q\in\mathbb{P}, \label{eq:NS_12_Ndim_WK}
  \\[0.6em]
  &\ipd{{\bm V}\cdot{\bm n},~\psi}_{\Gamma(t)} =\ipd{{\bm u}\cdot{\bm n},~\psi}_{\Gamma(t)}\qquad \forall \psi\in L^2(\Gamma(t)), \label{eq:Interface_Evolution_Ndim_WK}
  \\[0.6em]
  &\ipd{\kappa{\bm n},~{\bm g}}_{\Gamma(t)} =\ipd{\nabla_s{\bm X},~\nabla_s{\bm g}}_{\Gamma(t)}\qquad \forall {\bm g}\in [H^1(\Gamma(t))]^d,\label{eq:Curvature_Relation_Ndim_WK}
\end{align}
\end{subequations}
where we assume that the mesh velocity ${\bf  w}$ is provided with ${\bf w}_\Gamma = {\bf V}$.  Here we notice that \eqref{eq:NS_11_Ndim_WK} results from \eqref{eq:Inertia_Term_Weak_Form3} and \eqref{eq:Viscous_Term_Weak_Form}. \eqref{eq:NS_12_Ndim_WK} is due to the divergence-free condition \eqref{eq:NS_12_Ndim}.  \eqref{eq:Interface_Evolution_Ndim_WK} comes from \eqref{eq:Interface_Evolution_Ndim} and \eqref{eq:Curvature_Relation_Ndim_WK} is from \eqref{eq:Curvature_Relation_Ndim_Weak_Form1}.   
In the case when $\gamma(\cdot)$ is a constant, \eqref{eq:Two_Phase_Flow_WK} collapses to the weak formulation that was introduced in \cite[(4.9)]{garcke2023structure}.

\subsection{ALE weak formulation for substance concentrations}
We define the space-time bulk domain 
\begin{equation}
\mathcal{Q}^\pm_T:=\bigcup_{t\in[0,T]}\Omega^\pm(t)\times\{t\},\nonumber 
\end{equation} 
and denote by $\ipd{\cdot,\cdot}_{\Omega_\pm(t)}$ the $L^2$-inner product on $\Omega_\pm(t)$, respectively.  It follows from the Reynolds transport formula and the divergence free condition \eqref{eq:NS_12_Ndim} that
\begin{equation}
  \label{eq:Bulk_Integration_Time_Derivative}
    \mathrm{\frac{d}{dt}}\ipd{C,~\Phi}_{\Omega^\pm(t)}=\ipd{\partial_t^\bullet C,~\Phi}_{\Omega^\pm(t)}+\ipd{\partial_t^\bullet\Phi,~C}_{\Omega^\pm(t)}\quad\forall \Phi\in H^1(\mathcal{Q}^\pm_T),
\end{equation}
We next take the inner product of  \eqref{eq:Bulk_Mass_Conservation_DiffC_Ndim} with $\Phi\in H^1(\mathcal{Q}^\pm_T)$ and apply integration by parts to get
\begin{equation}
  \mathrm{\frac{d}{dt}}\ipd{C,~\Phi}_{\Omega^\pm(t)} - \ipd{\partial_t^\bullet\Phi,~C}_{\Omega^\pm(t)} + \frac{D_C^\pm}{Pe}\ipd{\nabla C,\nabla\Phi}_{\Omega^\pm(t)} + Da\ipd{J_s \pm  S^\pm,~\Phi}_{\Gamma(t)}=0,\label{eq:dcphi}
\end{equation}
%\footnote{\color{red} It seems that $J_s$ is not defined. }
where we use \eqref{eq:Bulk_Integration_Time_Derivative} and recall the boundary conditions \eqref{eq:Bulk_Flux_Condition_Ndim}, \eqref{eq:Bulk_Concentration_Outer_Condition_Ndim}. 

We also define the evolving surface
\begin{equation}
\mathcal{G}_T:=\bigcup_{t\in[0,T]}\Gamma(t)\times\{t\}.\nonumber 
\end{equation} 
Then it follows from the surface transport formula that
\begin{equation}
  \label{eq:Surface_Integration_Time_Derivative}
  \mathrm{  \mathrm{\frac{d}{dt}}}\ipd{ C_\Gamma,~\Psi_C}_{\Gamma(t)}=\ipd{\partial_t^\bullet C_\Gamma,\Psi_C }_{\Gamma(t)}+\ipd{C_\Gamma,\partial_t^\bullet\Psi_C}_{\Gamma(t)}+\ipd{C_\Gamma\nabla_s\cdot{\bm u},\Psi_C}_{\Gamma(t)}\quad\forall\Psi_C\in H^1(\mathcal{G}_T).
\end{equation}
Multiplying \eqref{eq:Surface_Mass_Conservation_DiffC_Ndim} with $\Psi_C\in H^1(\mathcal{G}_T)$ followed by an integration on $\Gamma(t)$, and using \eqref{eq:Surface_Integration_Time_Derivative} lead us to
\begin{align}
   \mathrm{  \mathrm{\frac{d}{dt}}}\ipd{C_\Gamma,\Psi_C}_{\Gamma(t)}-\ipd{ C_\Gamma,\partial_t^\bullet\Psi_C }_{\Gamma(t)}+\frac{D_{C_\Gamma}}{Pe_\Gamma}\ipd{\nabla_s C_\Gamma,\nabla_s\Psi_C}_{\Gamma(t)} = \ipd{\sum_\pm S^\pm-\mathcal{R},\Psi_C}_{\Gamma(t)}.
\end{align}
Similarly, taking inner product of \eqref{eq:Surface_Mass_Conservation_DiffA_Ndim} and \eqref{eq:Surface_Mass_Conservation_DiffB_Ndim} with $\Psi_A, \Psi_B\in H^1(\mathcal{G}_T)$ gives us 
\begin{align}
  \mathrm{\frac{d}{dt}}\ipd{ A_\Gamma,\Psi_A}_{\Gamma(t)} &- \ipd{ A_\Gamma,\partial_t^\bullet\Psi_A }_{\Gamma(t)}+ \frac{D_{A_\Gamma}}{Pe_\Gamma}\ipd{\nabla_s A_\Gamma,\nabla_s\Psi_A}_{\Gamma(t)} - \ipd{\mathcal{R},\Psi_A}_{\Gamma(t)}=0,\label{eq:Surface_Mass_Conservation_DiffA_Ndim_WK}
  \\ 
  \mathrm{  \mathrm{\frac{d}{dt}}}\ipd{ B_\Gamma,\Psi_B}_{\Gamma(t)} &- \ipd{ B_\Gamma,\partial_t^\bullet\Psi_B}_{\Gamma(t)}+ \frac{D_{B_\Gamma}}{Pe_\Gamma}\ipd{ \nabla_s B_\Gamma,\nabla_s\Psi_B}_{\Gamma(t)} +\ipd{\mathcal{R},\Psi_B}_{\Gamma(t)}=0.\label{eq:Surface_Mass_Conservation_DiffB_Ndim_WK}
\end{align}

We further define the function spaces that are compatible with the ALE mappings by 
\begin{subequations}
  \begin{align}
    &\mathbb{C}^\pm = \{\chi: \mathcal{Q}^\pm_T\rightarrow\mathbb{R},\quad \chi = \hat{\chi}\circ\mathcal{A}[t]^{-1},\quad\hat{\chi}\in H^1(\mathcal{O}^\pm)\}\subset H^1(\mathcal{Q}^\pm_T),\nonumber 
    \\
    &\mathbb{C}_\Gamma = \{\chi: \mathcal{G}_T\rightarrow\mathbb{R},\quad \chi = \hat{\chi}\circ\mathcal{A}[t]^{-1},\quad\hat{\chi}\in H^1(\Upsilon)\}\subset H^1(\mathcal{G}_T).\nonumber 
  \end{align}
\end{subequations}
Then the following identities hold 
\begin{equation}
  \partial_t^\circ \Phi = 0, \quad \forall \Phi\in\mathbb{C}^\pm;\quad \partial_t^\circ \Psi = 0, \quad \forall \Psi\in\mathbb{C}_\Gamma, \nonumber 
\end{equation}
since the reference frame is time independent. Then we can transfer the material derivative as follows
\begin{equation}
  \label{eq:material_derivative_alternative}
  \partial_t^\bullet \Phi = ({\bm u}-{\bm w})\cdot\nabla \Phi, \quad \forall \Phi\in\mathbb{C}^\pm;\quad \partial_t^\bullet \Psi_C = ({\bm u}-{\bm w}_\Gamma)\cdot\nabla_s \Psi_C, \quad \forall \Psi_C\in\mathbb{C}_\Gamma.
\end{equation}
on recalling \eqref{eq:Relation_Material_ALE_Bulk} and \eqref{eq:Relation_Material_ALE_Surface}. 
 
Collecting these results, we are ready to present the weak form of the substance concentrations. Given the substance concentrations $A_{\Gamma,0}:=A_\Gamma(\cdot,0),B_{\Gamma,0}:=B_\Gamma(\cdot,0),C_{\Gamma,0}:=C_\Gamma(\cdot,0)$ and $C_0^\pm=C^\pm(\cdot,0)$, we find $A_\Gamma(\cdot,t),B_\Gamma(\cdot,t),C_\Gamma(\cdot,t)\in H^1(\mathcal{G}_T)$ and $C^\pm(\cdot,t)\in H^1(\mathcal{Q}_T^\pm)$ such that for any $\Phi\in\mathbb{C}^\pm$ and $(\Psi_A,\Psi_B,\Psi_C)\in [\mathbb{C}_\Gamma]^3$
\begin{subequations}
\label{eq:Substance_WK_ALE}
\begin{align}  
  \mathrm{  \mathrm{\frac{d}{dt}}}\ipd{C,\Phi}_{\Omega^\pm(t)} &- (C,\ipd{{\bm u}-{\bm w})\cdot\nabla \Phi}_{\Omega^\pm(t)} + \frac{D_C^\pm}{Pe}\ipd{\nabla C,\nabla\Phi}_{\Omega^\pm(t)} + Da\ipd{ S^\pm \pm J_s,\Phi}_{\Gamma(t)}=0,\label{eq:Bulk_Mass_Conservation_DiffC_Ndim_WK_ALE}
  \\ 
  \mathrm{  \mathrm{\frac{d}{dt}}}\ipd{ A_\Gamma,\Psi_A}_{\Gamma(t)} &- \ipd{ A_\Gamma,({\bm u}-{\bm w}_\Gamma)\cdot\nabla_s \Psi_A }_{\Gamma(t)}+ \frac{D_{A_\Gamma}}{Pe_\Gamma}\ipd{\nabla_s A_\Gamma,\nabla_s\Psi_A}_{\Gamma(t)} - \ipd{\mathcal{R},\Psi_A}_{\Gamma(t)}=0,\label{eq:Surface_Mass_Conservation_DiffA_Ndim_WK_ALE}
  \\ 
  \mathrm{  \mathrm{\frac{d}{dt}}}\ipd{B_\Gamma,\Psi_B}_{\Gamma(t)} &- \ipd{ B_\Gamma,({\bm u}-{\bm w}_\Gamma)\cdot\nabla_s \Psi_B }_{\Gamma(t)}+ \frac{D_{B_\Gamma}}{Pe_\Gamma}\ipd{ \nabla_s B_\Gamma,\nabla_s\Psi_B}_{\Gamma(t)} + \ipd{\mathcal{R},\Psi_B}_{\Gamma(t)}=0,\label{eq:Surface_Mass_Conservation_DiffB_Ndim_WK_ALE}
  \\ 
  \mathrm{  \mathrm{\frac{d}{dt}}}\ipd{C_\Gamma,\Psi_C}_{\Gamma(t)} &- \ipd{ C_\Gamma,({\bm u}-{\bm w}_\Gamma)\cdot\nabla_s \Psi_C }_{\Gamma(t)} + \frac{D_{C_\Gamma}}{Pe_\Gamma}\ipd{\nabla_s C_\Gamma,\nabla_s\Psi_C}_{\Gamma(t)} + \ipd{\mathcal{R},\Psi_C}_{\Gamma(t)}\nonumber
  \\
  & \hspace{5cm}- \ipd{\sum_\pm S^\pm,\Psi_C}_{\Gamma(t)}=0.\label{eq:Surface_Mass_Conservation_DiffC_Ndim_WK_ALE}
\end{align}
\end{subequations}

%\subsection{The mass conservation law and the energy law}
%\subsection{The mass conservation law}
\subsection{Mass conservation within the weak formulation}

The ALE weak formulation for the whole system consists of \eqref{eq:Two_Phase_Flow_WK} and \eqref{eq:Substance_WK_ALE}. We have the following properties for the weak solution
\begin{itemize}
    \item Volume preservation of the fluid:  \begin{equation}
  \label{eq:Bulk_Volume_Evolution}
    \mathrm{\frac{d}{dt}}\, \text{vol}(\Omega^-(t)) =   \mathrm{\frac{d}{dt}} \int_{\Omega^-(t)} d\sL^{d} = \int_{\Gamma(t)}{\bf  V}\cdot{\bf n}\, d\sH^{d-1}=0.
\end{equation}
This could be proved by  setting \(q = \chi_{\Omega^-(t)} - \omega(t)\) in \eqref{eq:NS_12_Ndim_WK} and $\psi=1$ in \eqref{eq:Interface_Evolution_Ndim_WK},
  where $\omega(t)$ is  the normalized volume fraction of $\Omega^-(t)$
\[
  \omega(t) = \frac{\text{vol}(\Omega^-(t))}{\text{vol}(\Omega)}, \quad t \in (0,T].
\]
% and it is not difficult to show $(\chi_{\Omega^-(t)} - \omega(t))\in\mathbb{P}$. Then setting \(q = \chi_{\Omega^-(t)} - \omega(t)\) in \eqref{eq:NS_12_Ndim_WK} and $\psi=1$ in \eqref{eq:Interface_Evolution_Ndim_WK}, and recalling \eqref{eq:Bulk_Volume_Evolution}, we get
% \begin{align}
%   \mathrm{  \mathrm{\frac{d}{dt}}} \, \text{vol}(\Omega^-(t)) 
%   &=\int_{\Gamma(t)}{\bf V}\cdot{\bf n}\,d\sH^{d-1} = \int_{\Gamma(t)}{\bf u}\cdot{\bf n}\,d\sH^{d-1}= \int_{\Omega^-(t)} \nabla \cdot {\bm u} \, d\sL^{d} 
%   = \int_{\Omega} \chi_{\Omega^-(t)} \nabla \cdot {\bm u} \, d\sL^{d}  = 0,\nonumber 
% \end{align}
% which implies the volume preservation of the fluid within the weak formulation.
\item Conservation of mass: 
\begin{subequations}
\begin{align*}
  \mathrm{  \mathrm{\frac{d}{dt}}} \, m_s(A_\Gamma, B_\Gamma, 0, 0; t) 
  &=   \mathrm{\frac{d}{dt}} \ipd{ A_\Gamma, 1}_{\Gamma(t)} 
   +   \mathrm{\frac{d}{dt}} \ipd{ B_\Gamma, 1}_{\Gamma(t)} = 0, \\
  \mathrm{  \mathrm{\frac{d}{dt}}} \, m_s(A_\Gamma, 0, C_\Gamma, C; t) 
  &=   \mathrm{\frac{d}{dt}} \ipd{ A_\Gamma, 1 }_{\Gamma(t)} 
   +   \mathrm{\frac{d}{dt}} \ipd{ C_\Gamma, 1 }_{\Gamma(t)} 
   + \frac{1}{Da}   \mathrm{\frac{d}{dt}}\ipd{C, 1} = 0.
\end{align*}
\end{subequations}
which  could be proved by choosing \(\Phi = \Psi_A = \Psi_B = \Psi_C = 1\) in \eqref{eq:Substance_WK_ALE}.  

\end{itemize}

\section{Finite Element Approximations}\label{sec:Finite_Element_Approximations}
In this section, we first present the temporal-spatial discretization and then propose the ALE-FEM approximations of the two weak formulations \eqref{eq:Two_Phase_Flow_WK} and \eqref{eq:Substance_WK_ALE}.  

\subsection{Discretization}

We partition the time domain uniformly as $[0,T]=\bigcup_{m=1}^{M}[t_{m-1},t_m]$ with $t_m = m\Delta t$ and time step size $\Delta t = T/M$. For any variable $\zeta(\cdot, t)$, we denote by $\zeta^m$ the numerical approximation of $\zeta$ at time $t_m$.

\paragraph{\textsl{Interface discretization}:}
We approximate the interface $\Gamma(t_m)$ by a $(d-1)$-dimensional polyhedral surface which is given by 
\[\Gamma^m = \bigcup_{j=1}^{J_\Gamma} \overline{\sigma_j^m}\quad\mbox{with}\quad \mathcal{T}_\Gamma^m = \{\sigma_j^m: j = 1,\cdots, J_\Gamma\},\quad Q_\Gamma^m = \{{\bm q}_k^m: k = 1,\cdots, K_\Gamma\},\]
where  $\mathcal{T}_\Gamma^m$ is a collection of mutually disjoint $(d-1)$-simplices, and $Q_\Gamma^m$ is a collection of the vertices. Let $\{{\bm q}_{j_k}^{m}\}_{k=0}^{d-1}$ be the vertices of $\sigma_j^m$, ordered with the same orientation for all $\sigma_j^m\in\mathcal{T}_\Gamma^m$. For simplicity,  we denote $\sigma_j^m=\Delta \{{\bm q}_{j_k}^m\}_{k=0}^{d-1}$. 
 
 Then we introduce the unit normal ${\bm n}$ to $\Gamma^m$; that is, 
\begin{equation*}
  {\bm n}_j^m:={\bm n}_j^m|_{\sigma_j^m}:=\frac{{\bm A}\{\sigma_j^m\}}{|{\bm A}\{\sigma_j^m\}|}\quad \text{with}\quad {\bm A}\{\sigma_j^m\}=({\bm q}_{j_1}^m-{\bm q}_{j_0}^m)\wedge\cdots\wedge({\bm q}_{j_{d-1}}^m-{\bm q}_{j_0}^m),
\end{equation*}
where $\wedge$ is the wedge product and ${\bm A}\{\sigma_j^m\}$ is the orientation vector of $\sigma_j^m$. To approximate the inner product $\ipd{\cdot,\cdot}_{\Gamma(t_m)}$, we introduce mass-lumped products over the current polyhedral surface $\Gamma^m$ via
\begin{align*}
    \ipd{ u, v}_{\Gamma^m}^h := \frac{1}{d} \sum_{j=1}^{J_\Gamma}|\sigma_j^m|\sum_{k=0}^{d-1}\lim_{\sigma_j^m\ni {\bm p} \rightarrow {\bm q}_{j_k}^m}(u\cdot v)({\bm p}), 
\end{align*}
where $u,v$ are piecewise continuous, with possible jumps across the edges of $\{\sigma_j^m\}_{j=1}^{J_\Gamma}$, and $|\sigma_j^m|=\frac{1}{(d-1)!}|{\bm A}\{\sigma_j^m\}|$ is the measure of $\sigma_j^m$.

We next introduce the surface finite element space
\begin{equation}
  V^h(\Gamma^m):=\{\chi\in C(\Gamma^m): \chi|_{\sigma}~\text{is affine} \quad \forall \sigma\in\mathcal{T}_\Gamma^m\}.\nonumber 
\end{equation}

Following the work in \cite{Bao2021SJNA,Bao2023NMPDE}, define the interpolated polyhedral surfaces $\Gamma^h(t)$ via a linear interpolation of nodes between $\Gamma^{m}$ and $\Gamma^{m+1}$:
\begin{equation*}
  \Gamma^h(t):=\frac{t_{m+1}-t}{\Delta t}\Gamma^m + \frac{t-t_{m}}{\Delta t}\Gamma^{m+1} \quad \text{and} \quad \Gamma^h(t)=\bigcup_{j=1}^{J_\Gamma}\overline{\sigma_j^h}(t),\quad t\in [t_m,t_{m+1}],
\end{equation*}
where $\{\sigma_j^h(t)\}_{j=1}^{J_\Gamma}$ are mutually disjoint $(d-1)-$simplices and the vertices $\{{\bm q}_k^h(t)\}_{k=1}^{K_\Gamma}$ of $\Gamma^h(t)$ are given by
\begin{equation*}
  {\bm q}^h_k(t) = \frac{t_{m+1}-t}{\Delta t}{\bm q}_k^{m} + \frac{t-t_{m}}{\Delta t}{\bm q}_k^{m+1}, \quad t\in [t_m,t_{m+1}], \quad k=1,\cdots,K_\Gamma.
\end{equation*}
We then define the time-weighted normals ${\bm n}^{m+\frac{1}{2}}\in [L^\infty(\Gamma^m)]^d$
\begin{equation}\label{eq:twnormal}
  {\bm n}^{m+\frac{1}{2}}|_{\sigma_j^m} = {\bm n}^{m+\frac{1}{2}}_j:=\frac{1}{\Delta t|{\bm A}\{\sigma_j^m\}|}\int_{t_m}^{t_{m+1}} {\bm A}\{\sigma_j^h(t)\} dt \quad \forall 1\leq j\leq J_\Gamma.
\end{equation}
Let $\Omega^{-,m}$ be the interior of $\Gamma^m$ and $\Omega^{+,m}$ be the exterior of $\Gamma^m$ in $\Omega$. Then we have the following lemma.
\begin{lma}
  Let ${(\bm X}^{m+1}\in [V^h(\Gamma^m)]^d$. Then it holds
  \begin{equation}
  \label{eq:Vol_Conservation_Dis_Zhao}
  \ipd{{\bm X}^{m+1}-{\bm X}^m)\cdot{\bm n}^{m+\frac{1}{2}},1}_{m}^h = {\rm vol}(\Omega^{-,m+1}) - {\rm vol}(\Omega^{-,m}).
  \end{equation}
\end{lma}
\begin{proof}
  The details of the proof can be found in \cite{Bao2021SJNA, Bao2023NMPDE}.
\end{proof}

 \paragraph{\textsl{Bulk discretization}:} 
We consider the discretization of the bulk domain $\Omega$. At time $t_m$, a regular partition of $\Omega$ is given by 
\begin{equation}
  \overline{\Omega} = \cup_{j=1}^{J_\Omega^m} \overline{o_j^m}\quad\mbox{with}\quad \mathcal{T}^m = \{o_j^m: j =1,\cdots, J_\Omega\},\quad Q^m=\{{\bm a}_k^m: k =1,\cdots, K_\Omega\},\nonumber 
\end{equation}
where $\mathcal{T}^m$ is a collection of mutually disjoint open simplices, and $Q^m$ is a collection of the vertices of the mesh. We employ the fitted mesh strategy so that the interface triangulation $\Gamma^m$ aligns with the bulk partition $\mathcal{T}^m$, i.e., $\Gamma^m \subsetneq \{ \partial o : o\in\mathcal{T}^m \}.$
Let $\Omega_-^m$  and $\Omega_+^m$ denote the interiors and exterior of $\Gamma^m$. In our fitted mesh approach, we can divide the elements of $\mathcal{T}^m$ into interior and exterior elements as 
\begin{equation}  
\mathcal{T}^{m}_\pm:=\left\{o\in\mathcal{T}^m: o\subset\Omega_\pm^m\right\}.\nonumber 
\end{equation}  
 
We define the bulk finite element spaces:
\begin{subequations}\label{eq:Bulk_FE_Space}
\begin{align*}
  &S_k^m:=\{\chi\in C(\overline{\Omega}): \chi|_{o}\in P_k(o)\quad\forall o\in\mathcal{T}^m\},\quad  k\in \mathbb{N}^+,
  \\
  &S_0^m:=\{\chi\in L^2(\Omega): \chi|_{o}\in P_0(o) \quad\forall o\in\mathcal{T}^m\},\\
  &S_k^{\pm,m}:=\{\chi\in C(\overline{\Omega_\pm^m}): \chi|_{o}\in P_k(o)\quad\forall o\in\mathcal{T}_\pm^m\}, \quad k\in \mathbb{N}^+,\label{eq:Sub_Bulk_FE_Space}
\end{align*}
\end{subequations}
where $P_k(o)$ denotes the space of polynomials with degree at most $k$ on the element $o$. In this work,  we employ a stable Taylor-Hood element:
\begin{equation}\label{eq:UPspace}
  (\mathbb{U}^m,\mathbb{P}^m) = ([S_2^m]^d\cap\mathbb{U},(S_1^m+S_0^m)\cap\mathbb{P}),
\end{equation}
for the finite element spaces for the discrete fluid velocity and pressure, which satisfies the LBB inf-sup stability condition.  Furthermore, the density $\rho({\bf x},t)$ and viscosity function $\eta({\bf x}, t)$ in \eqref{eq:rhoetaf} are approximated by $\rho^m, \eta^m\in S_0^m$ such that
\begin{equation}
  \rho^m= \rho^+\chi_{\Omega_+^{m}} + \rho^-\chi_{\Omega_-^{m}},\quad\eta^m= \eta^+\chi_{\Omega_+^{m}} + \eta^-\chi_{\Omega_-^{m}}.\nonumber 
\end{equation}

\subsection{Discrete ALE mappings}

We now consider the moving mesh approach. In general, we construct $\mathcal{T}^m$ base on $\mathcal{T}^{m-1}$. In particular, we keep the mesh connectivity and topology unchanged and update the vertices of the mesh according to 
\begin{equation}\label{eq:verticesupdate}
  {\bm a}_k^{m} = {\bm a}_k^{m-1} + {\bm\Phi}^m({\bm a}_k^{m-1}), \quad 1\leq k\leq K_\Omega,\quad 1\leq m\leq M,
\end{equation}
where ${\bm\Phi}^m\in [S_1^{m-1}]^d$ is the displacement of the bulk mesh which can be obtained as follows. On introducing 
\begin{align*}
    &\mathbb{Y}^{m-1} = \{{\bm \chi}\in[S_1^{m-1}]^d: {\bm\chi}\cdot{\bm n}=0~\text{on}~\partial\Omega;~{\bm\chi}={\bm X}^{m}-{\bm X}^{m-1}~\text{on}~\Gamma^{m-1}\}, 
    \\
    &\mathbb{Y}^{m-1}_0 = \{{\bm \chi}\in[S_1^{m-1}]^d: {\bm\chi}\cdot{\bm n}=0~\text{on}~\partial\Omega;~{\bm\chi}={\bm 0}~\text{on}~\Gamma^{m-1}\},
\end{align*}
then  we find  ${\bm\Phi}^{m}\in\mathbb{Y}^{m-1}$  by solving the following elastic equation
\begin{equation}
  \label{eq:Elastic_Equation_WK_Dis}
  2(\lambda^{m-1} D({\bm \Phi}^{m}),D({\bm \psi})) + (\lambda^{m-1}\nabla\cdot{\bm\Phi}^{m},\nabla{\bm\psi})=0\quad\forall{\bm\psi}\in \mathbb{Y}_0^{m-1},
\end{equation}
where $\lambda^{m-1}\in S_0^{m-1}$ is used to limit the distortion of small elements \cite{MASUD1997CMAME,ZHAO2020JCP}, and defined as
\begin{equation*}
  \lambda^{m-1}|_{o_j^{m-1}} = 1 + \frac{\max\limits_{o\in\mathcal{T}^{m-1}}|o|-\min\limits_{o\in\mathcal{T}^{m-1}}|o|}{|o_j^{m-1}|},\quad j = 1,\cdots,J_\Omega. 
\end{equation*}

At time $t_m$, it is natural to define the discrete mesh velocity $\bm{w}^m \in S_1^m$ with
\begin{equation*}
  \bm{w}^m(\vec a_k^m) :=  \frac{\bm{a}_k^m - \bm{a}_k^{m-1}}{\Delta t} \quad\forall k = 1,\cdots, K_\Omega.
\end{equation*}
The corresponding discrete ALE mapping $\mathcal{\bm{A}}^m[t]({\bf x}) \in [S_1^m]^d$ for $t \in [t_{m-1}, t_m]$ is given by the linear interpolation
\begin{equation*}
  \mathcal{\bm{A}}^m[t]({\bf x}):= \sum_{k=1}^{K_\Omega} \left( \frac{t_m - t}{\Delta t} \bm{a}_k^{m-1} + \frac{t - t_{m-1}}{\Delta t} \bm{a}_k^m \right) \phi_k^m({\bf x}),
\end{equation*}
where $\phi_k^m$ is the nodal basis function of $S_1^m$ at ${\bf a}_k^m$.  It should be noted that $\mathcal{\bm{A}}^m[t_m] = {\bf id}|_{\Omega}$ is the identity map, and  
\begin{equation*}
  \mathcal{{\bm A}}^{m}[t_{m-1}] = {\bf id}|_{\Omega}-\Delta t\,{\bm w}^m\quad\text{and}\quad\mathcal{{\bm A}}^{m}[t_{m-1}]({\bm a}_k^{m}) = {\bm a}_k^{m-1},\quad k=1,\cdots,K_\Omega,
\end{equation*} 
which implies that $\mathcal{{\bm A}}^{m}[t_{m-1}](\Omega^{m}_\pm) = \Omega^{m-1}_\pm$. Let $\mathcal{J}^m$ be the Jacobian determinant of the element-wise linear map $\mathcal{\bm{A}}^m[t_{m-1}]$. It is not difficult to show that
\begin{equation}
  \label{eq:Jacobian_Dis}
  \mathcal{J}^m:={\rm det}(\nabla\mathcal{\bm A}^m[t_{m-1}]) ={\rm det}(\mathbb{I}_d-\Delta t\nabla{\bm w}^m) = 1 - \Delta t\nabla\cdot{\bm w}^m + O(\Delta t^2),
\end{equation}
Here, \eqref{eq:Jacobian_Dis} will provide a consistent approximation of the first two terms in \eqref{eq:NS_12_Ndim_WK}. 

\subsection{ALE-FEM approximations for Navier-Stokes equations}
We are now ready to present the finite element approximation for the fluid-structure interface part \eqref{eq:Two_Phase_Flow_WK}, which couple the standard ALE-FEM for the Navier-Stokes with the parametric FEM for the interface evolution. For simplicity, we define
\begin{equation*}
\zeta_{\mathcal{A}}=\zeta^\circ\mathcal{A}^m[t_{m-1}]\in\mathbb{U}^{m}\quad\forall\zeta\in\mathbb{U}^{m-1}.
\end{equation*}
Let $\Gamma^0$ be a polyhedral approximation of the fluid interface and $\mathcal{T}^0$ be a regular fitted partition of $\Omega$. We set ${\bf X}^0={\bf id}|_{\Gamma^0}$, $\Gamma^{-1}=\Gamma^0$, $\Omega_\pm^{-1}=\Omega_\pm^0$ and $\mathcal{J}^0({\bf x})=1$. Given ${\bf u}^0\in \mathbb{U}^0$, for $m\geq 0$, we seek ${\bm u}^{m+1}\in\mathbb{U}^m$, $p^{m+1}\in\mathbb{P}^m$, ${\bm X}^{m+1}\in[V^h(\Gamma^m)]^d$ and $\kappa^{m+1}\in V^h(\Gamma^m)$ such that
\begin{subequations}\label{eq:Two_Phase_Flow_WK_Dis}
  \begin{align}
    &\ipd{\rho^m\frac{{\bm u}^{m+1}-{\bm u}^m_{\mathcal{\bm A}}\sqrt{\mathcal{J}^m}}{\Delta t},{\bm v}^h}+ \mathscr{A}\ipd{\rho^m,{\bm u}^m_{\mathcal{\bm A}}-{\bm w}^m;{\bm u}^{m+1},{\bm v}^h} - \ipd{p^{m+1},\nabla\cdot{\bm v}^h}\nonumber
    \\
    &\qquad+\frac{2}{Re}\ipd{\eta^m D({\bm u}^{m+1}),D({\bm v}^h)}+\frac{1}{We}\ipd{\gamma^m\kappa^{m+1}{\bm n}^m-\nabla_s\gamma^m,{\bm v}^h}_{\Gamma^m}=0\quad\forall{\bm v}^h\in\mathbb{U}^m,\label{eq:NS_11_Ndim_WK_Dis}
    \\[0.5em]
    &\ipd{\nabla\cdot{\bm u}^{m+1},q^h} = 0\quad \forall q^h\in\mathbb{P}^m,\label{eq:NS_12_Ndim_WK_Dis}
    \\[0.5em]
    &\ipd{\frac{{\bm X}^{m+1}-{\bm X}^m}{\Delta t}\cdot{\bm n}^{m+\frac{1}{2}},\psi^h}_{m}^h - \ipd{{\bm u}^{m+1}\cdot{\bm n}^m,\psi^h}_{\Gamma^m}=0\quad\forall\psi^h\in V^h(\Gamma^m),\label{eq:Interface_Evolution_Ndim_WK_Dis}
    \\[0.5em]
    &\ipd{\kappa^{m+1}{\bm n}^{m+\frac{1}{2}},{\bm g}^h}_{\Gamma^m}^h + \ipd{\nabla_s{\bm X}^{m+1},\nabla_s{\bm g}^h}_{\Gamma^m}=0\quad\forall {\bm g}^h\in[V^h(\Gamma^m)]^h,\label{eq:Curvature_Relation_Ndim_WK_Dis}
  \end{align}
\end{subequations}
and we then set $\Gamma^{m+1}=\vec X^{m+1}(\Gamma^m)$ to construct the new bulk mesh $\mathcal{T}^{m+1}$ through \eqref{eq:verticesupdate}, \eqref{eq:Elastic_Equation_WK_Dis}. 

Here $\gamma^m$ is an explicit approximation of $\gamma(\cdot)$ defined through \eqref{eq:gamma} such that
 \[\gamma^m = 1 - E\sum_{K}\omega_K K^m, \]
 with $K^m$ being the approximation of the substance concentration $K$ at the time $t_m$. Computations of $K^m$ are based on discretization of \eqref{eq:Substance_WK_ALE}, which will be discussed in the next subsection. On the other hand,  ${\bf n}^{m+\frac{1}{2}}$ is a time-weighted interface normals defined in \eqref{eq:twnormal} to enable an exact volume preservation (see \cite[Theorem 4.4]{garcke2013diffuse}). 
 
 Following \cite{garcke2013diffuse}, on recalling \eqref{eq:Jacobian_Dis},  the first term in \eqref{eq:NS_11_Ndim_WK_Dis} can be rewritten as 
\begin{align*}
  &\bigg(\rho^m\frac{{\bm u}^{m+1}-{\bm u}^m_{\mathcal{\bm A}}\sqrt{\mathcal{J}^m}}{\Delta t},{\bm v}^h\bigg)\nonumber
  \\
 =&\bigg(\rho^m\frac{{\bm u}^{m+1}-{\bm u}^m_{\mathcal{\bm A}}}{\Delta t},{\bm v}^h\bigg)+\frac{1}{2}(\rho^m\nabla\cdot{\bm w}^m,{\bm u}^m_{\mathcal{\bm A}}\cdot{\bm v}^h) + O(\Delta t),
\end{align*}
which is hence a consistent temporal discretization of the first two terms in \eqref{eq:NS_11_Ndim_WK}. We notice that in the case when $\gamma^m$ is a constant, \eqref{eq:Two_Phase_Flow_WK_Dis} collapse to the ALE structure-preserving method in \cite[(4.20)]{garcke2023structure}, where an unconditional stability in terms of the discrete energy can be proved for \eqref{eq:Two_Phase_Flow_WK_Dis}.

\subsection{ALE-FEM approximations for concentrations}
We next consider numerical approximations for the substance concentrations \eqref{eq:Substance_WK_ALE}. Given the initial the initial bulk substance concentrations $C^0$ and the interfacial surfactant concentration $A_\Gamma^0,B_\Gamma^0,C_\Gamma^0$, for each $m\geq 0$ and we are given $\Gamma^k$ with $k\leq m+1$,  we then find  $C^{m+1}\in S_1^{\pm,m+1}$ and $(A_\Gamma^{m+1},B_\Gamma^{m+1},C_\Gamma^{m+1})\in [V^h(\Gamma^{m+1})]^3$ such that
\begin{subequations}
\label{eq:Substance_WK_Dis}
\begin{align}  
  \frac{1}{\Delta t}&\ipd{C^{m+1},\Phi^h}_{\Omega_\pm^{m+1}} - \ipd{C^{m+1}({\bm u}^{m+1}-{\bm w}^{m+1}),\nabla\Phi^h}_{\Omega_{\pm}^{m+1}} + \frac{D_C^\pm}{Pe}(\nabla C^{m+1},\nabla\Phi^h)_{\Omega_\pm^{m+1}} \nonumber
  \\
  &+ Da\ipd{S^{\pm,m+1} \pm J_s^{m+1},\Phi^h}_{\Gamma^{m+1}}=\frac{1}{\Delta t}\ipd{C^{m},\Phi^h\circ\mathcal{A}^{m+1}[t_m])^{-1}}_{\Omega_\pm^{m+1}},\label{eq:Bulk_Mass_Conservation_DiffC_Ndim_WK_Dis}
  \\[0.5em]
  \frac{1}{\Delta t}&\ipd{A_\Gamma^{m+1},\Psi_A^h}_{\Gamma^{m+1}} - \ipd{ A_\Gamma^{m+1}({\bm u}^{m+1}-{\bm w}^{m+1}_\Gamma), \nabla_s\Psi_A^h}_{\Gamma^{m+1}}+ \frac{D_{A_\Gamma}}{Pe_\Gamma}\ipd{\nabla_s A_\Gamma^{m+1},\nabla_s\Psi_A^h}_{\Gamma^{m+1}} \nonumber
  \\
  &-\ipd{\mathcal{R}^{m+\frac{1}{2}},\Psi_A^h}{\Gamma^{m+1}}=\frac{1}{\Delta t}\ipd{ A_\Gamma^m,\Psi_A^h\circ\mathcal{A}^{m+1}[t_m])^{-1}}_{\Gamma^m},\label{eq:Surface_Mass_Conservation_DiffA_Ndim_WK_Dis}
  \\[0.5em]
  \frac{1}{\Delta t}&\ipd{B_\Gamma^{m+1},\Psi_B^h}_{\Gamma^{m+1}} - \ipd{ B_\Gamma^{m+1}({\bm u}^{m+1}-{\bm w}^{m+1}_\Gamma),\nabla_s\Psi_B^h}_{\Gamma^{m+1}}+ \frac{D_{B_\Gamma}}{Pe_\Gamma}\ipd{ \nabla_s B_\Gamma^{m+1},\nabla_s\Psi_B^h}_{\Gamma^{m+1}} \nonumber 
  \\
  &+ \ipd{\mathcal{R}^{m+\frac{1}{2}},\Psi_B^h}_{\Gamma^{m+1}}=\frac{1}{\Delta t}\ipd{ B_\Gamma^m,\Psi_B^h\circ(\mathcal{A}^{m+1}[t_m])^{-1}}_{\Gamma^m},\label{eq:Surface_Mass_Conservation_DiffB_Ndim_WK_Dis}
  \\[0.5em]
  \frac{1}{\Delta t}&\ipd{C_\Gamma^{m+1},\Psi_C^h}_{\Gamma^{m+1}} - \ipd{C_\Gamma^{m+1}({\bm u}^{m+1}-{\bm w}^{m+1}_\Gamma),\nabla_s\Psi_C^h }_{\Gamma^{m+1}}+ \frac{D_{C_\Gamma}}{Pe_\Gamma}\ipd{\nabla_s C_\Gamma^{m+1},\nabla_s\Psi_C^h}_{\Gamma^{m+1}} \nonumber 
  \\
  &+ \ipd{\mathcal{R}^{m+\frac{1}{2}},\Psi_C^h}_{\Gamma^{m+1}}-\ipd{\sum_\pm S^{\pm,m+\frac{1}{2}},\Psi_C^h}_{\Gamma^{m+1}}=\frac{1}{\Delta t}\ipd{ C_\Gamma^m,\Psi_C^h\circ(\mathcal{A}^{m+1}[t_m])^{-1}}_{\Gamma^m},\label{eq:Surface_Mass_Conservation_DiffC_Ndim_WK_Dis}
\end{align}
\end{subequations}
for all $(\Phi^h, \Psi_A^h, \Psi_B^h, \Psi_C^h)\in S_1^{\pm,m+1}\times [V^h(\Gamma^{m+1})]^3$, 
where $\mathcal{R}^{m+\frac{1}{2}}$ and $S^{\pm,m+\frac{1}{2}}$ are linear approximations to the reaction term $\mathcal{R}(A_\Gamma,B_\Gamma,C_\Gamma)$ and the source term $S^\pm(C^\pm,C_
\Gamma)$, respectively. In particular, they are given by
\begin{align*}
&\mathcal{R}^{m+\frac{1}{2}} = Bi\bigg(\frac{1}{2} k_f \left(B^{m+1}_\Gamma (\tilde{B}_\Gamma^m)^{\omega_b-1} (\tilde{C}_\Gamma^m)^{\omega_c}+C^{m+1}_\Gamma (\tilde{B}_\Gamma^m)^{\omega_b} (\tilde{C}_\Gamma^m)^{\omega_c-1}\right)
 -k_rA^{m+1}_\Gamma (\tilde{A}_\Gamma^m)^{\omega_a-1}\bigg),\nonumber 
 \\
 &S^{\pm,m+\frac{1}{2}} = Bi \left(k_{ad}^\pm C^{\pm,m+1}-k_d^\pm C^{m+1}_\Gamma (\tilde{C}_\Gamma^m)^{\omega_c-1}\right),\nonumber 
\end{align*}
where 
\begin{equation*}
 \tilde{A}_\Gamma^m = A^m_\Gamma\circ(\mathcal{A}^{m+1}[t_m])^{-1},\quad \tilde{B}_\Gamma^m = B^m_\Gamma\circ(\mathcal{A}^{m+1}[t_m])^{-1}, \quad \tilde{C}_\Gamma^m = C^m_\Gamma\circ(\mathcal{A}^{m+1}[t_m])^{-1}. 
\end{equation*}

\subsection{The overall ALE scheme}
The approximation of the whole scheme consists of \eqref{eq:Two_Phase_Flow_WK_Dis} and \eqref{eq:Substance_WK_Dis}. Here \eqref{eq:Two_Phase_Flow_WK_Dis} leads to a ``weakly'' nonlinear system due to the presence of the time-weighted interface normals. We solve the nonlinear system efficiently through a Picard-type iterative method. This leads to a linear system at each iterative step, which can be solved via the preconditioned GMRES method. The details of the solver can be found in \cite{garcke2023structure}. Scheme \eqref{eq:Substance_WK_Dis} leads to a lienar system which can be solved via th sparse LU factorization.   The overall procedure of the proposed numerical scheme can be summarized as follows: 
\begin{algorithm}
\caption{Stable ALE method for reactive semi-permeable interfaces}
\begin{algorithmic}[1]
\State Initialize: set $m \gets 0$; prepare the initial data $(\Omega^0, \Gamma^0)$ and $({\bf u}^0, C^0, A_\Gamma^0, B_\Gamma^0, C_\Gamma^0)$; 
\While{$m \leq M$}
    \State Solve  \eqref{eq:Two_Phase_Flow_WK_Dis} on the mesh  $\mathcal{T}^m$ to get $({\bf u}^{m+1}, p^{m+1}, {\bf  X}^{m+1}, \kappa^{m+1})$; 
    \State Based on $\mathcal{T}^m$ and ${\bf  X}^{m+1}$, we update the mesh $\mathcal{T}^{m+1}$ through \eqref{eq:verticesupdate}, \eqref{eq:Elastic_Equation_WK_Dis}; 
    \State Compute  \eqref{eq:Substance_WK_Dis} for the substance concentrations $(C^{m+1}, A_\Gamma^{m+1}, B_\Gamma^{m+1}, C_\Gamma^{m+1})$; 
    \State $m \gets m + 1$
\EndWhile
\end{algorithmic}
\end{algorithm}

We have the following theorem for the volume preservation of the fluid and the mass conservation of the substance on fully discrete level. 
\begin{thm} The introduced scheme satisfies the volume conservation
\begin{equation}
  \label{eq:Volume_Conservation_Dis}
  {\rm vol}(\Omega_-^{m+1}) = {\rm vol}(\Omega_-^{m}),
\end{equation}
and the mass conservations that 
\begin{subequations} \label{eq:Mass_Conservation_Dis}
\begin{align}
  m_s(A_\Gamma, B_\Gamma, 0, 0; t_{m+1}) &= m_s(A_\Gamma, B_\Gamma, 0, 0; t_m), \\
  m_s(A_\Gamma, 0, C_\Gamma, C; t_{m+1}) &= m_s(A_\Gamma, 0, C_\Gamma, C; t_m),
\end{align}
\end{subequations}
\end{thm}

\begin{proof}
\eqref{eq:Volume_Conservation_Dis} is a direct consequence of the time-weighted approximation \eqref{eq:twnormal} and the chosen finite element space \eqref{eq:UPspace}. Its proof can be found in \cite[Theorem 4.4]{garcke2023structure}.

For the substance concentrations, it is simple to choose $\Phi = \Psi_A = \Psi_B = \Psi_C = 1$ in Eq.~\eqref{eq:Substance_WK_Dis} and combine these equations, we obtain the discrete mass conservation laws \eqref{eq:Mass_Conservation_Dis}, where the total substance mass is defined as
\begin{equation*}
  m_s(A_\Gamma, B_\Gamma, C_\Gamma, C; t_m) = \sum_K \langle K^m, 1 \rangle_{m} + \frac{1}{Da} (C^m, 1)_m.
\end{equation*}

\end{proof}

% \begin{rmk}
% It doesn't  appear possible to obtain a discrete energy law for the introduced scheme. Nevertheless, in the case of a constant surface tension force, we are able to obtain an unconditional energy stability for the introduced scheme, see \cite[Theorem 4.4]{garcke2023structure}.     
% \end{rmk}

\section{Numerical Results}\label{sec:Numerical_Results}

This section presents a sequence of numerical experiments to validate the proposed model and demonstrate its capability to capture reactive interfacial dynamics with selective permeability. We begin in Subsection~\ref{subsec:conv_tese} with a convergence test to assess the accuracy of the numerical scheme and verify mass and volume conservation. Next, we consider a model problem of cholesterol efflux mediated by the ABCG1 transporter. This test, in the absence of fluid motion, illustrates the model’s ability to capture multistage surface reactions and transmembrane solute exchange relevant to lipid metabolism. We then proceed to simulations of fluid-structure interactions driven by interfacial chemistry, including Marangoni-induced droplet migration, shear deformation, and buoyancy-driven bubble dynamics. These examples underscore the intricate coupling between flow, interfacial dynamics, and reactive transport, as uniquely captured by the proposed framework.

\subsection{Convergence test and conservation law} \label{subsec:conv_tese}

In this subsection, we assess the convergence order of the proposed numerical method by conducting simulations with varying mesh resolutions and time steps. The computational domain is set to $\Omega = [-0.5, 0.5]^2$, and the initial interface is defined by the ellipse
\[
\Gamma(0) = \left\{(x, y) \,\bigg|\, \left(\frac{x}{1.25}\right)^2 + \left(\frac{y}{0.8}\right)^2 = 0.25^2 \right\}.
\]
We initialize the velocity field as ${\bm u}({\bm x},0) = \bm{0}$, and set the initial concentrations as $C({\bm x},0) = 0.8$ in the bulk, and $A_\Gamma({\bm x},0) = B_\Gamma({\bm x},0) = C_\Gamma({\bm x},0) = 0.8$ on the interface. Neumann boundary conditions are imposed for the bulk concentration, i.e., ${\bm n} \cdot \nabla C^+ = 0$ on $\partial \Omega$. For the velocity field, we prescribe ${\bm u} = \bm{0}$ on the no-slip boundary $\partial\Omega_1$, and free-stress conditions $\mathbb{T} \cdot \bm{n} = 0$ on $\partial\Omega_2$.

The physical and numerical parameters are chosen as follows:
\[
\begin{aligned}
&\rho^+ = 10, \quad \rho^- = 1, \quad \mu^+ = 10, \quad \mu^- = 1, \quad Re = 10, \quad Ca = 0.1, \quad Da = 1, \\
&Bi = 0.4, \quad E = 0.1, \quad Pe = 1, \quad Pe_\Gamma = 1, \quad \omega_a = \omega_b = \omega_c= \lambda_a = \lambda_c = 1, \\
&k_r = k_d^\pm = 1, \quad D_{A_\Gamma} = D_{B_\Gamma} = D_{C_\Gamma} = 1, \quad D^+ = 0.5, \quad D^- = 1.
\end{aligned}
\]

We consider three levels of spatial discretization, corresponding to interface meshes with $J_\Gamma = 16$, $32$, and $64$ elements. The associated time step sizes are chosen as $\Delta t = 2.0 \times 10^{-2}$, $5.0 \times 10^{-3}$, and $1.25 \times 10^{-3}$, respectively.

Table~\ref{tab:conv_test} reports the numerical errors at time $t = 10$ in the interfacial arc length $L$, the interfacial concentrations $A_\Gamma$, $B_\Gamma$, $C_\Gamma$, and the bulk concentration $C$, along with their convergence orders. The reference values $L^e$, $A_\Gamma^e$, $B_\Gamma^e$, $C_\Gamma^e$, and $C^e$ correspond to the equilibrium states. The results demonstrate that the errors decrease consistently with mesh refinement, and the observed convergence rates are approximately second-order, confirming the expected accuracy of the numerical scheme.

\begin{table}[htb!]
\centering
\caption[scale=1]{Errors of $L, A_\Gamma, B_\Gamma, C_\Gamma$ and $C$, and convergence orders at $t=10$.}
\scalebox{0.85}{
\begin{tabular}{ccccccccccc}
\hline
 $J_\Gamma$ & $|L-L^e|$  & order & $|A_\Gamma-A_\Gamma^e|$ & order & $|B_\Gamma-B_\Gamma^e|$ & order & $|C_\Gamma-C_\Gamma^e|$ & order & $|C-C^e|$ & order \\
\hline
     16     & $1.01\times 10^{-2}$ &   -  & $5.34\times 10^{-5}$ &   -  & $5.31\times 10^{-5}$ &  -   & $1.11\times 10^{-4}$ &  -   & $1.11\times 10^{-4}$ & -    \\
     32     & $2.52\times 10^{-3}$ & 2.00 & $1.34\times 10^{-5}$ & 1.99 & $1.31\times 10^{-5}$ & 2.02 & $2.76\times 10^{-5}$ & 2.00 & $2.80\times 10^{-5}$ & 2.00 \\
     64     & $6.30\times 10^{-4}$ & 2.00 & $3.43\times 10^{-6}$ & 1.97 & $3.14\times 10^{-6}$ & 2.06 & $6.62\times 10^{-6}$ & 2.06 & $7.36\times 10^{-6}$ & 1.93 \\
\hline
\end{tabular}
}
\label{tab:conv_test}
\end{table} 

\begin{figure}[htb!]
    \centering
    \includegraphics[scale=0.55]{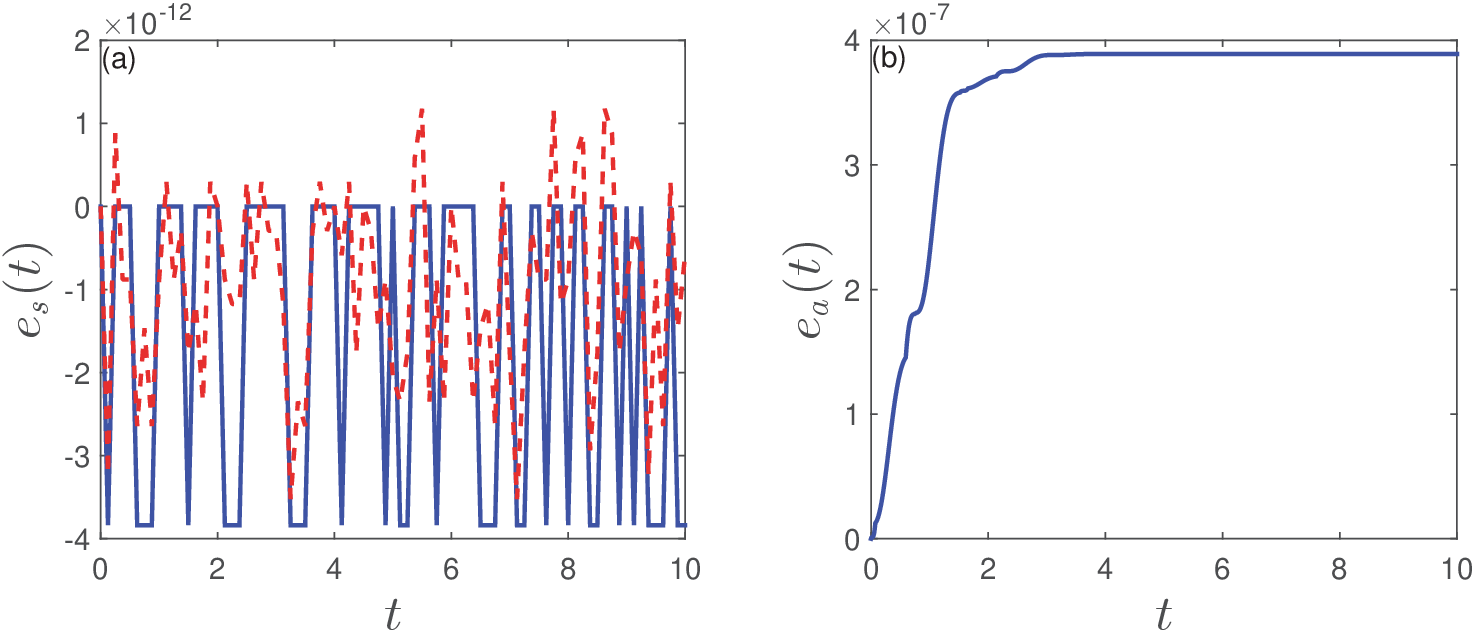}
    \caption{(a) Time evolution of the relative mass error. The blue solid line and red dashed line correspond to the errors in $m_s(A_\Gamma, B_\Gamma, 0, 0; t)$ and $m_s(A_\Gamma, 0, C_\Gamma, C; t)$, respectively. (b) Time evolution of the relative area error.}
    \label{fig:relaxation_area_mass}
\end{figure}

In Fig.~\ref{fig:relaxation_area_mass}, we plot the time evolution of the relative changes in total area and substance mass. The relative area error is defined as
\begin{equation*}
    e_a(t)\big|_{t=t_m} = \frac{\text{vol}(\Omega^{-,m}) - \text{vol}(\Omega^{-,0})}{\text{vol}(\Omega^{-,0})},
\end{equation*}
and the relative mass error is defined as
\begin{equation*}
    e_s(t)\big|_{t=t_m} = \frac{m_s(\cdot; t_m) - m_s(\cdot; 0)}{m_s(\cdot; 0)}.
\end{equation*}

The results confirm that both the mass and volume are well-conserved over time, consistent with the theoretical predictions given by Eqs.~\eqref{eq:Volume_Conservation_Dis} and~\eqref{eq:Mass_Conservation_Dis}.

\begin{figure}[htb!]
  \centering
  \includegraphics[scale=0.55]{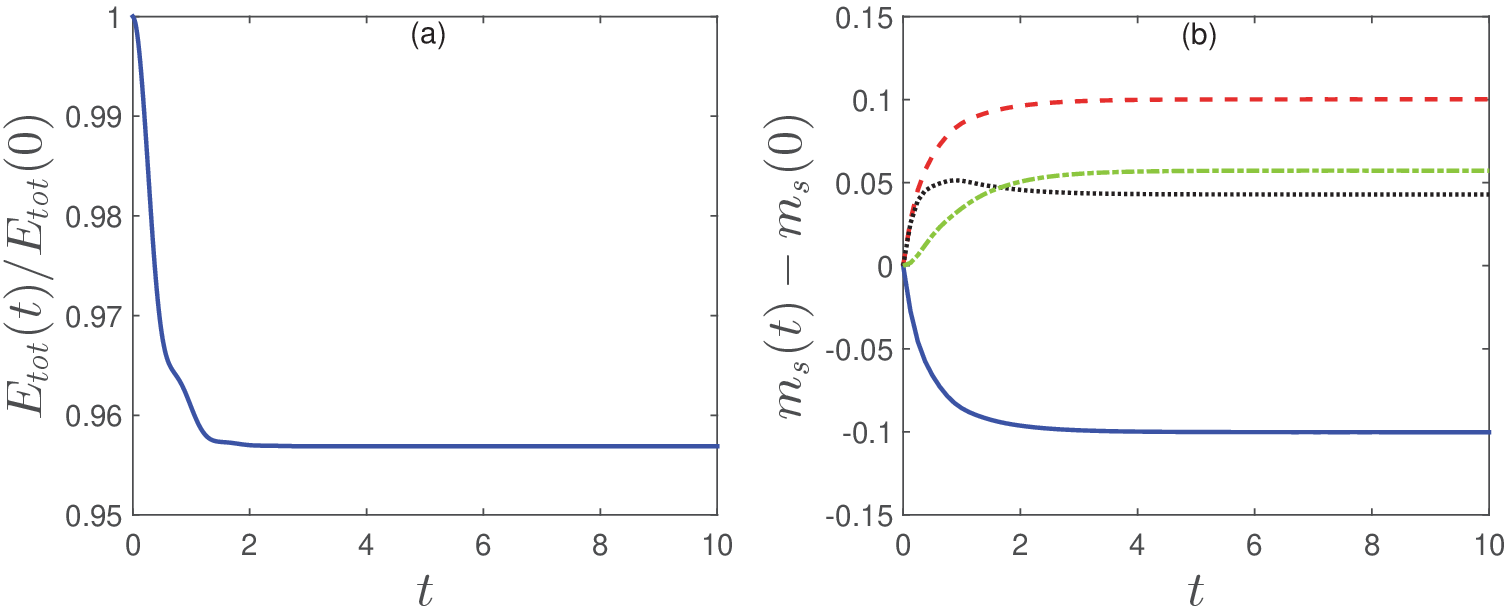}
  \caption{(a) Time evolution of the normalized total energy $E_{tot}(t)/E_{tot}(0)$; (b) Evolution of the individual substance masses: $m_s(A_\Gamma,0,0,0;t)$ (blue solid line), $m_s(0,B_\Gamma,0,0;t)$ (red dashed line), $m_s(0,0,C_\Gamma,0;t)$ (black dotted line), and $m_s(0,0,0,C;t)$ (green dash-dotted line).}
  \label{fig:relaxation_energy_mass}
\end{figure}

Fig.~\ref{fig:relaxation_energy_mass} shows the evolution of the normalized total energy $E_{tot}(t)/E_{tot}(0)$ and the individual mass components $m_s(A_\Gamma, 0, 0, 0; t)$, $m_s(0, B_\Gamma, 0, 0; t)$, $m_s(0, 0, C_\Gamma, 0; t)$, and $m_s(0, 0, 0, C; t)$. Although a discrete energy law has not been rigorously established for the fully discretized scheme, the left panel of Fig.~\ref{fig:relaxation_energy_mass} clearly demonstrates a monotonic decrease in total energy over time, which is consistent with the theoretical behavior predicted by the continuum model. The right panel of Fig.~\ref{fig:relaxation_energy_mass} indicates that when the system is sufficiently relaxed, different chemical species can reach their equilibrium states, with their mass gains or losses unchanged. This is also consistent with the theoretical prediction.

\begin{figure}[htbp!]
    \centering
    \includegraphics[scale=0.6]{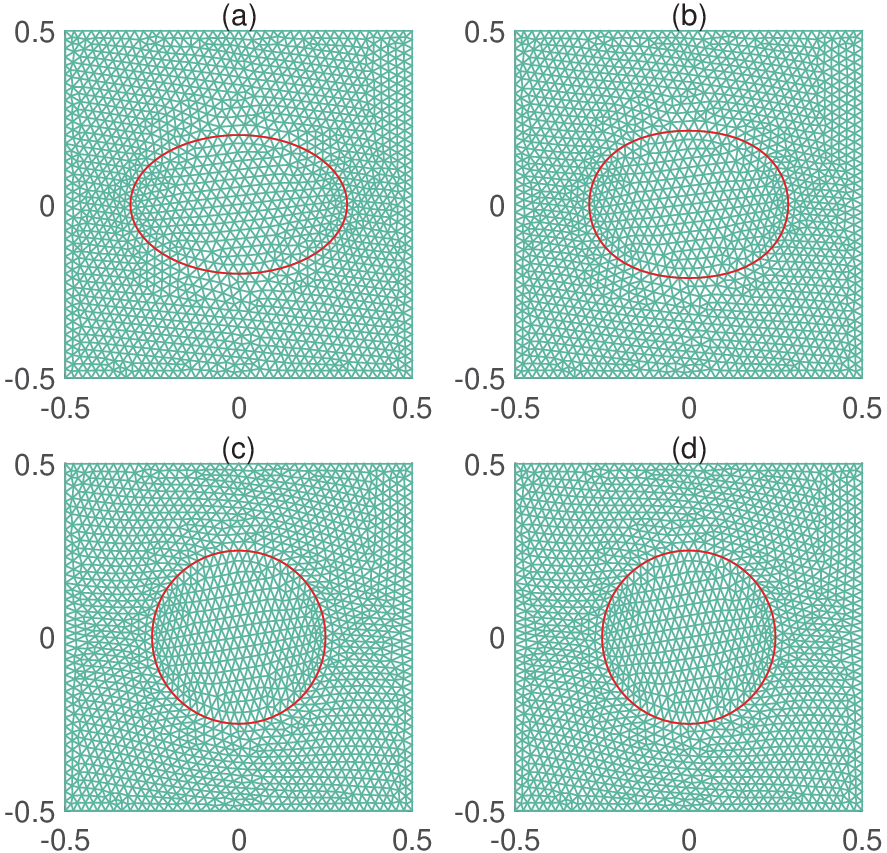}
    \caption{Snapshots of the fluid interface together with the mesh at different times, where $J_\Gamma=64$ and $J_\Omega=4010$. (a) $t=0$; (b) $t=0.2$; (c) $t=2$; (d) $t=10$.}
    \label{fig:relaxation_mesh}
\end{figure}

In Fig.~\ref{fig:relaxation_mesh}, we display snapshots of the fluid interface and mesh at selected time points: $t = 0$, $0.2$, $2$, and $10$. The interface is discretized using $J_\Gamma = 64$ markers, while the bulk mesh consists of $K_\Omega = 1940$ vertices and $J_\Omega = 4010$ triangles. The visualizations clearly illustrate the relaxation of the initially elliptical interface toward a circular shape, reflecting the system's tendency toward a lower-energy equilibrium configuration.

% \begin{figure}[htbp!]
%   \centering
%   \includegraphics[scale=0.55]{images/relaxation_mesh}
%   \caption{Snapshots of the fluid interface and bulk mesh at times $t = 0$, $0.2$, $2$, and $10$. The interface evolves from an initial ellipse toward a circular shape.}
%   \label{fig:relaxation_mesh}
% \end{figure}

\begin{rmk}
  The above problem exists an equilibrium state, and $L^e, A_\Gamma^e, B_\Gamma^e, C_\Gamma^e$ and $C^e$ can be obtained by solving the following equations
  \begin{align}
    &\mathcal{R}(A_\Gamma^e,B_\Gamma^e,C_\Gamma^e)=0,\quad S(C_\Gamma^e,C^e)=0,\nonumber
    \\
    &m_s(A_\Gamma^e,B_\Gamma^e,0,0;\infty) = m_s(A_\Gamma({\bm x},0),B_\Gamma({\bm x},0),0,0;0),\nonumber
    \\
    &m_s(A_\Gamma^e,0,C_\Gamma^e,C^e;\infty) = m_s(A_\Gamma({\bm x},0),0,C_\Gamma({\bm x},0),C({\bm x},0);0).\nonumber
  \end{align}
\end{rmk}
\subsection{Reaction-induced mass transfer: cholesterol efflux via the ABCG1 pathway}
Cholesterol homeostasis is essential for maintaining cellular function and membrane integrity, especially in macrophages, neurons, and vascular endothelial cells. Excess membrane cholesterol is removed through specialized transport mechanisms, forming a key component of reverse cholesterol transport. Among these, the ATP-binding cassette sub-family G member 1 (ABCG1) transporter plays a central role in facilitating the efflux of cholesterol and phospholipids to high-density lipoproteins (HDL)~\cite{kennedy2005abcg1}.

% In this subsection, we neglect fluid dynamics and focus on the concentration evolution governed by our proposed model, illustrating its application to cholesterol efflux mediated by interfacial reactions. This efflux pathway is schematically represented in Fig.~\ref{fig:ABCG1}. Initially, cholesterol with concentration \( C^- \) is confined within the intracellular region \( \Omega_1^- \)  and interface $\Gamma_1$. A portion of the cholesterol binds to the plasma membrane interface \( \Gamma_1 \), where it reacts with the transporter protein ABCG1, modeled as an interfacial species with concentration \( B_{\Gamma_1} \), to form a complex \( A_{\Gamma_1} \). This complex then dissociates into the extracellular space \( \Omega^+ \) as soluble species \( A \). Upon reaching the HDL surface \( \Gamma_2 \), the species \( A \) is absorbed and undergoes a second reaction with the enzyme LCAT (lecithin–cholesterol acyltransferase), represented by the interfacial species \( F_{\Gamma_2} \), producing \( G_{\Gamma_2} \). Finally, \( G_{\Gamma_2} \) dissociates into the HDL interior as the encapsulated product \( G \). This sequence mimics the biological process of ABCG1-facilitated cholesterol transport and its subsequent enzymatic processing by HDL-associated factors.

\begin{figure}[htbp!]
    \centering
    \includegraphics[width=0.85\linewidth]{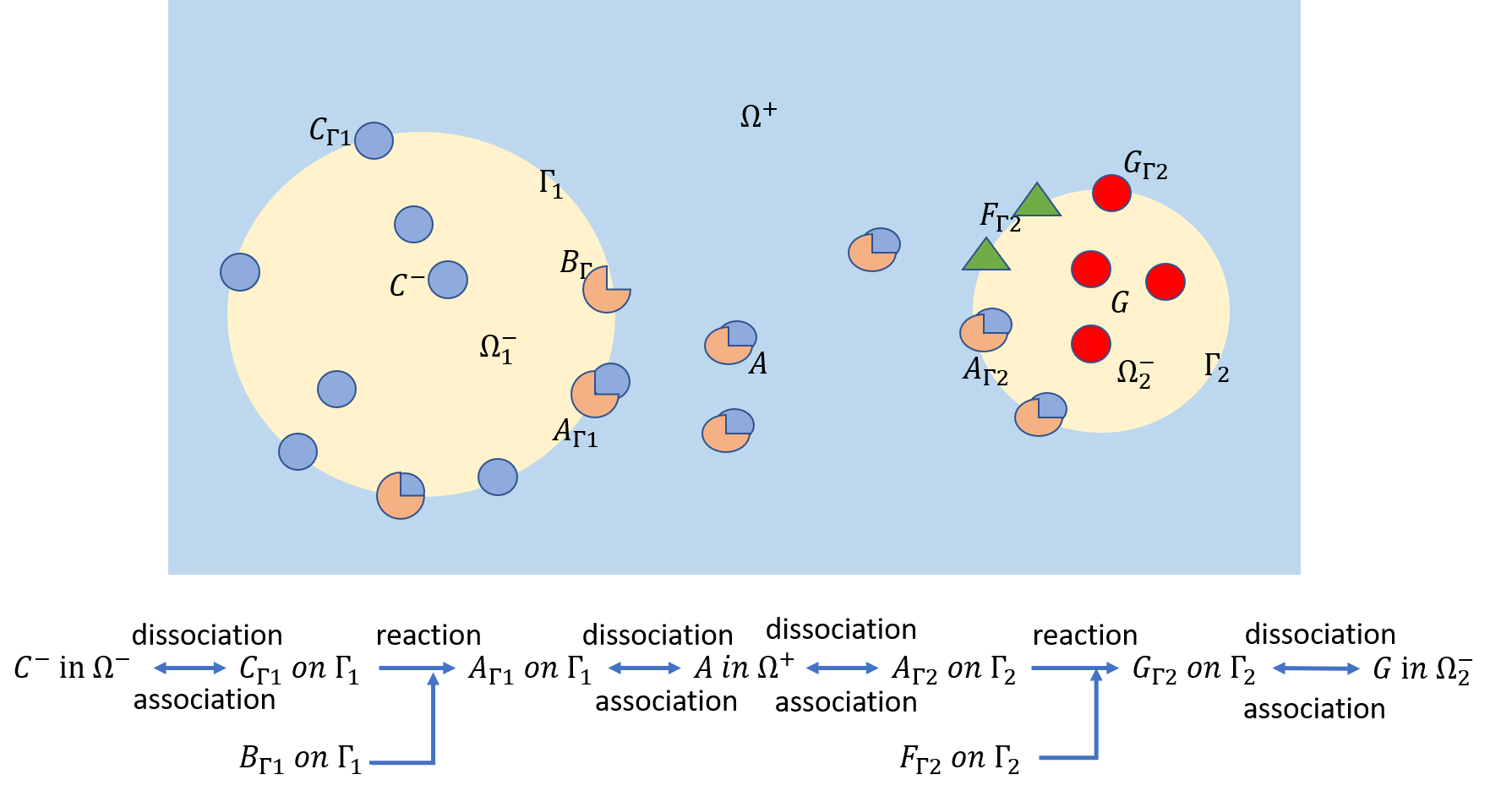}
    \caption{Schematic figure of cholesterol efflux via the ABCG1 pathway. The whole domain consists of the cell inner region $\Omega_1^-$, the HDL particle inner region $\Omega_2^-$, and the extracellular space $\Omega^+$. }
    \label{fig:ABCG1}
\end{figure}

In this subsection, we neglect fluid dynamics and focus on concentration evolution to illustrate the applicability of our proposed model to cholesterol efflux mediated by interfacial reactions. This efflux pathway, which is typically activated under conditions of elevated membrane cholesterol, is schematically represented in Fig.~\ref{fig:ABCG1}. 
To reflect biological reality, we assume that cholesterol is initially concentrated on the inner plasma membrane interface \( \Gamma_1 \), with a low bulk intracellular concentration \( C^- \), corresponding to the quiescent cytosolic state. 
The process begins with membrane-bound cholesterol \( C_{\Gamma_1} \) reacting with the transporter protein ABCG1, modeled as an interfacial species \( B_{\Gamma_1} \), to form a complex \( A_{\Gamma_1} \). This complex subsequently dissociates into the extracellular region \( \Omega^+ \) as a soluble intermediate species \( A \). The extracellular species \( A \) then diffuses to the surface of high-density lipoprotein (HDL), represented by interface \( \Gamma_2 \), where it binds to and reacts with the enzyme LCAT (lecithin–cholesterol acyltransferase), modeled as the interfacial species \( F_{\Gamma_2} \), producing a second complex \( G_{\Gamma_2} \). Finally, this complex dissociates into the HDL interior as the encapsulated product \( G \).
This multistage pathway captures key aspects of ABCG1-mediated cholesterol trafficking, including membrane-to-extracellular efflux, carrier-mediated uptake, and enzymatic processing at the lipoprotein interface. It also highlights how interfacial reactions and transport coupling govern the overall directionality and efficiency of lipid clearance.

Mathematically, we constrain the $C$ inside $\Omega_1^-$ and $G$ inside $\Omega_2^-$, the interface permeability is set to be zero, i.e. $J_s=0$. The communication of $C_{\Gamma1}$ (or \textbf{$G_{\Gamma2}$ }) with the bulk region is only from inner region $\Omega_1^-$ (or $\Omega_2^-$) and the flux from outer region $\Omega^+$   is zero. Similarly, the communication of $A_{\Gamma1}$ (or $A_{\Gamma2}$) with the bulk region is only between the region $\Omega^+$ and the interface $\Gamma_1$ (or $\Gamma_2$). 
So, the initial values  are set as follows
\[
C(\bm{x}, 0) =  
\begin{cases}
0.1, & \text{for } \bm{x} \in \Omega_1^-, \\
0, & \text{for } \bm{x} \in \Omega_2^-,\\
0, & \text{for } \bm{x} \in \Omega^+,\\
\end{cases}
~~~A(\bm{x}, 0)=G(\bm{x},0)=0, \bm{x} \in \Omega, 
\]
on interface $\Gamma_1$, 
\[A_{\Gamma1}=0, C_{\Gamma1}=1, B_{\Gamma_1}=1,\]
on interface $\Gamma_2$, 
\[A_{\Gamma2}=G_{\Gamma2}=0, F_{\Gamma_2}=1.\]
And the parameters are as follows
\[
\begin{aligned}
&\omega_a=\omega_b=\omega_c=1,\quad Bi = 0.1,\quad Da = 0.1, \quad Pe_\Gamma = 10, \quad Pe = 1,
\\
& D^- = 1,\quad D^+ = 5,\quad k_r =  k_d^\pm = 1, \quad \lambda_c = 10, \quad\lambda_{a,C,\Gamma_1} =\lambda_{a,A,\Gamma_2} = 10, \\
& \lambda_{a,A,\Gamma_1} = \lambda_{a,G,\Gamma_2} = 0.1,\quad D_{A_\Gamma} = D_{B_\Gamma} = D_{C_\Gamma} =  D_{F_\Gamma} = D_{G_\Gamma} = 5.
\end{aligned}
\]

\begin{figure}[htbp!]
    \centering
    \includegraphics[scale=0.6]{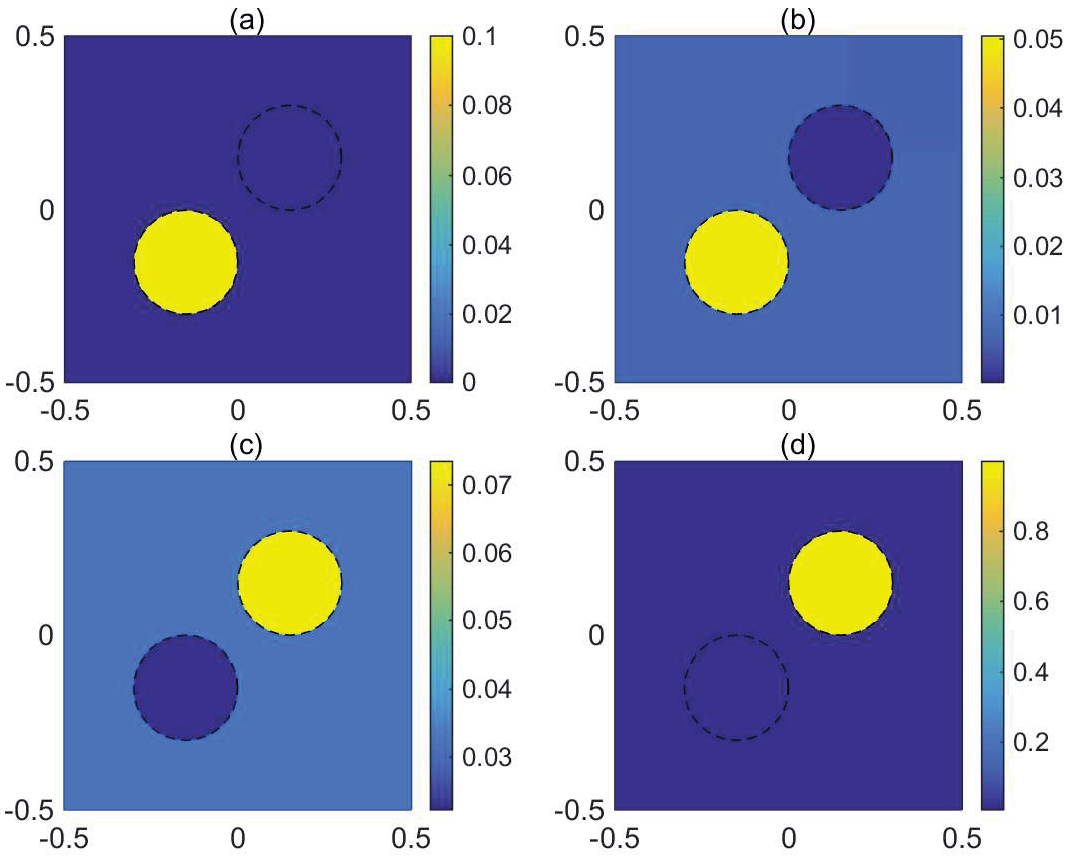}
    \caption{Snapshots of the bulk solute concentrations \( C^- \), \( A \), and \( G \) at selected times. (a) \( t = 0 \); (b) \( t = 1 \); (c) \( t = 10 \); (d) \( t = 100 \). The lower left circle denotes the cell, and the upper right circle denotes the HDL particle. The dash lines indicate the interfaces. }
    \label{fig:mass_transfer_sb}
\end{figure}

\begin{figure}[htbp!]
    \centering
    \includegraphics[scale=0.5]{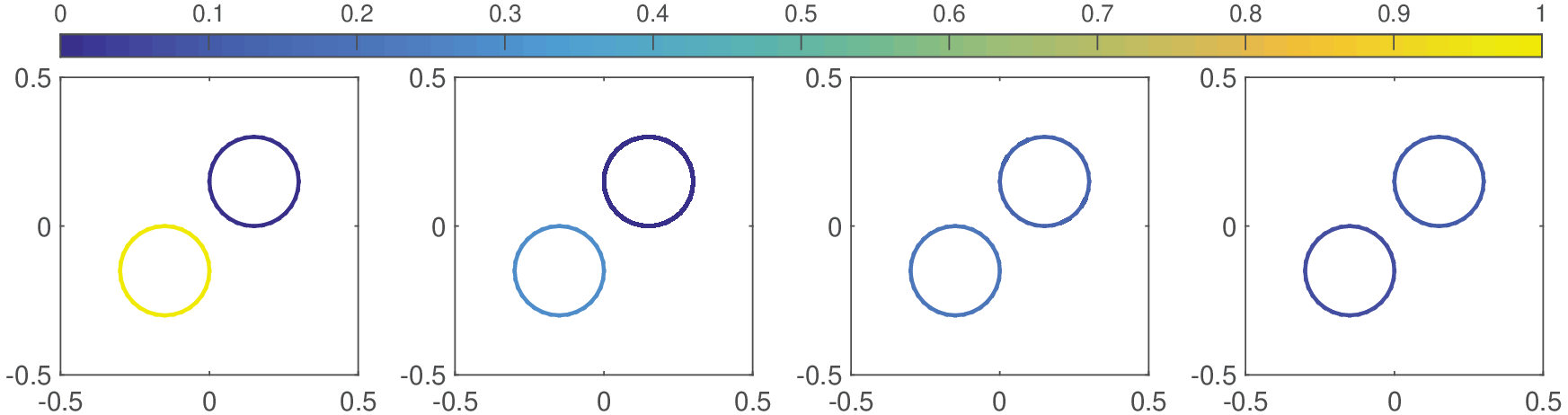}
    \caption{Snapshots of the interfacial solute concentrations \( C_{\Gamma_1} \) and \( G_{\Gamma_2} \) at $t=0,1,10,100$.  }
    \label{fig:mass_transfer_sf}
\end{figure}

\begin{figure}[htbp!]
    \centering
    \includegraphics[scale=0.6]{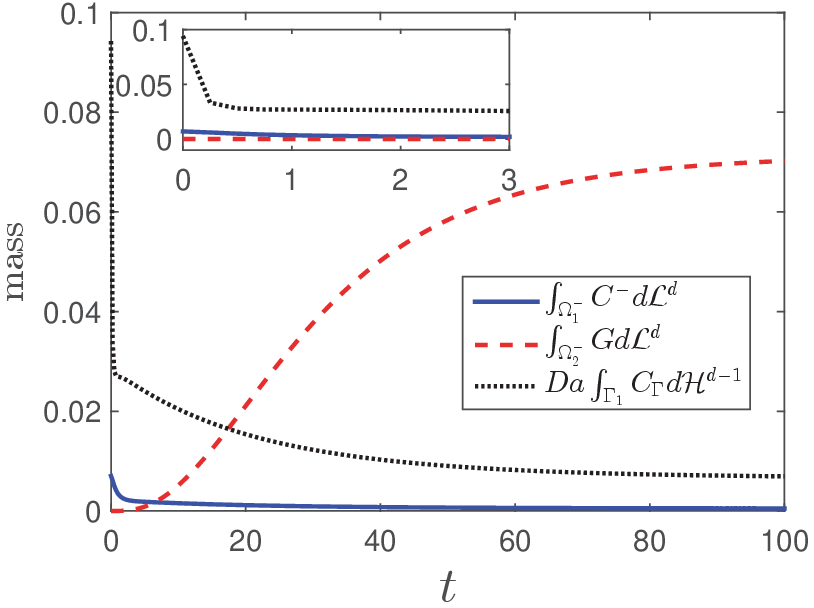}
    \caption{Temporal evolution of the total mass of intracellular cholesterol \( C^- \) inside the cell (blue line), cholesterol on the interface $\Gamma_1$  (black dot line) and  HDL-encapsulated cholesterol \( G \) (red dash line).  }
    \label{fig:mass_transfer_mass_evolution}
\end{figure}

Fig.~\ref{fig:mass_transfer_sb} illustrates: the intracellular concentration \( C^- \) within the cell domain \( \Omega_1^- \), the extracellular intermediate species \( A \) in \( \Omega^+ \), and the final captured form \( G \) within the HDL domain at different time slots. 
Initially, cholesterol is localized exclusively in the interface  \( \Gamma_1 \). Through the interfacial reaction on the plasma membrane, cholesterol at the interface (\( C_{\Gamma_1} \)) is rapidly converted and released into the extracellular region as \( A \), resulting in a decrease of \( C_{\Gamma_1}\) and an increase in \( A_{\Gamma} \). As time progresses, the accumulated extracellular \( A \) is gradually absorbed by the HDL interface \( \Gamma_2 \), where it undergoes a secondary reaction to form \( G_{\Gamma_2} \). This interfacial product is then transported into the HDL interior as the bulk species \( G \). 
Fig.~\ref{fig:mass_transfer_sf} depicts the different time slots of interfacial species \( C_{\Gamma_1} \) and \( G_{\Gamma_2} \). The corresponding total masses in the intracellular domain \( \Omega_1^- \), HDL region \( \Omega_2^- \), and on the membrane interface \( \Gamma_1 \) are summarized in Fig.~\ref{fig:mass_transfer_mass_evolution}, which confirms the cholesterol transformation from the cell membrane to the HDL particle. This multistage transformation, from cholesterol release on the cell membrane to HDL-mediated uptake, captures the essential characteristics of ABCG1-regulated cholesterol efflux. The simulation highlights the coordinated interplay of membrane-bound reactions, interfacial transport, and selective transmembrane exchange that underpins efficient and targeted cholesterol clearance from cells to circulating lipoproteins.

\subsection{Reaction-driven droplet motion} \label{subsec:drop_motion_dirichlet}

In this subsection, we investigate the motion of a droplet driven by interfacial chemical reactions. The droplet is placed in the computational domain $\Omega = [-0.5, 0.5]^2$. The initial concentration of the bulk substance is prescribed as
\[
C(\bm{x}, 0) = 
\begin{cases}
x + 0.5, & \text{for } \bm{x} \in \Omega^+, \\
0, & \text{for } \bm{x} \in \Omega^-,
\end{cases}
\]
which implies that no solute is initially present within the droplet.
The initial interfacial concentrations are set to $A_\Gamma(\bm{x}, 0) = 0$, $B_\Gamma(\bm{x}, 0) = 1$, and $C_\Gamma(\bm{x}, 0) = 0.5$. The initial fluid velocity is zero, i.e., ${\bm u}(\bm{x}, 0) = \bm{0}$.

For the boundary conditions, Dirichlet conditions are applied to the left and right boundaries to maintain the concentration gradient, while no-flux (Neumann) conditions are imposed on the top and bottom boundaries:
\begin{align*}
    {\bm u} = \bm{0}, \quad C = 0.5 + x, &\quad \text{on } \partial\Omega_1, \\
    \mathbb{T} \cdot \bm{n} = \bm{0}, \quad \bm{n} \cdot \nabla C^+ = 0, &\quad \text{on } \partial\Omega_2.
\end{align*}

The physical and numerical parameters used in this experiment are summarized as follows:
\[
\begin{aligned}
&\rho^+ = 0.1, \quad \rho^- = 1, \quad \mu^+ = 0.01, \quad \mu^- = 1, \quad Re = 10, \quad Ca = 0.1, \quad Bi = 0.8, \\
&Da = 1, \quad E = 1 \times 10^{-3}, \quad Pe_\Gamma = 10, \quad Pe = 1, \quad \lambda_a = 1, \quad \lambda_c = 1, \quad k_r = 0.25, \\
&k_d^\pm = 0.25, \quad D_{A_\Gamma} = 1, \quad D_{B_\Gamma} = 1, \quad D_{C_\Gamma} = 1, \quad D^- = 1, \quad D^+ = 0.1.
\end{aligned}
\]

To emphasize the dominant influence of species $A_\Gamma$ on surface tension, we set the weighting coefficients as $\omega_a = 100$, $\omega_b = 1$, and $\omega_c = 1$. The time step is uniformly set to $\Delta t = 10^{-3}$.

\begin{figure}[htb!]
    \centering
    \includegraphics[scale=0.6]{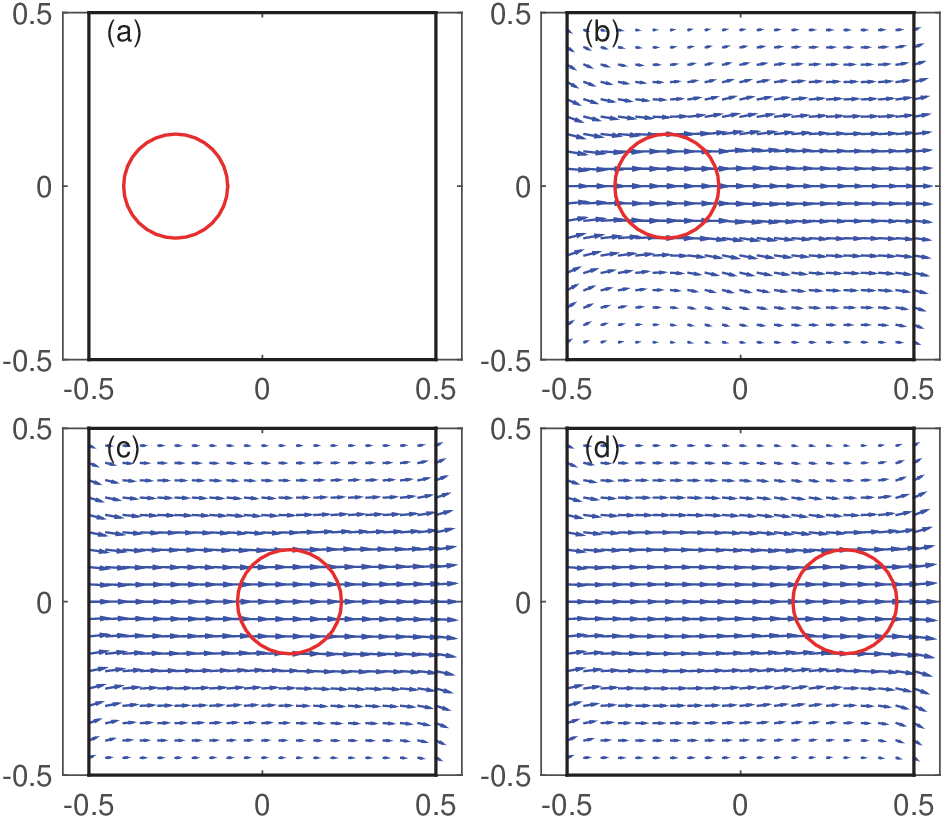}
    \caption{Snapshots of the fluid interface and velocity field at selected times. (a) $t = 0$, $\max_{\bm{x} \in \Omega} |\bm{u}| = 0$, $K_\Omega = 386$, $J_\Omega = 838$; (b) $t = 10$, $\max |\bm{u}| = 0.0064$, $K_\Omega = 386$, $J_\Omega = 838$; (c) $t = 50$, $\max |\bm{u}| = 0.0057$, $K_\Omega = 408$, $J_\Omega = 882$; (d) $t = 100$, $\max |\bm{u}| = 0.0040$, $K_\Omega = 386$, $J_\Omega = 838$.}
    \label{fig:var_weight_dirichlet_vel}
\end{figure}

\begin{figure}[htb!]
    \centering
    \includegraphics[scale=0.6]{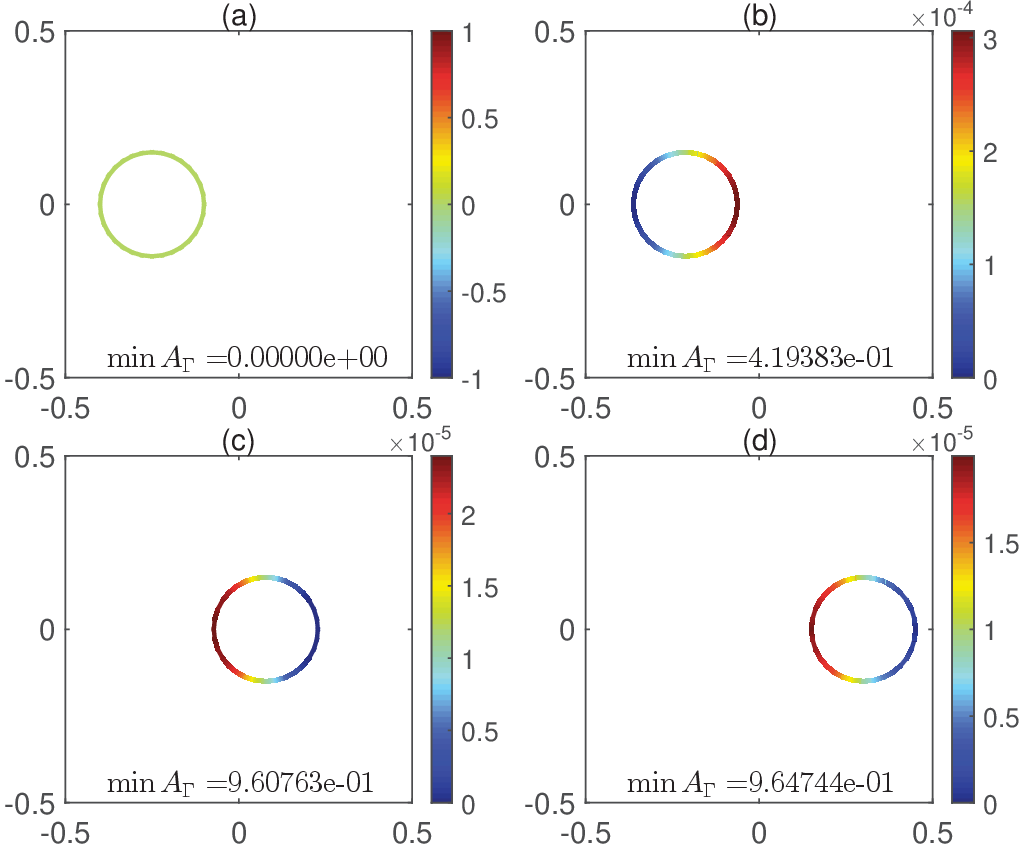}
    \caption{Snapshots of the interfacial substance concentration $A_\Gamma-\min_\Gamma A_\Gamma$ at selected times. (a) $t=0$; (b) $t=10$; (c) $t=50$; (d) $t=100$.}
    \label{fig:var_weight_dirichlet_sfa}
\end{figure}

\begin{figure}[htb!]
    \centering
    \includegraphics[scale=0.6]{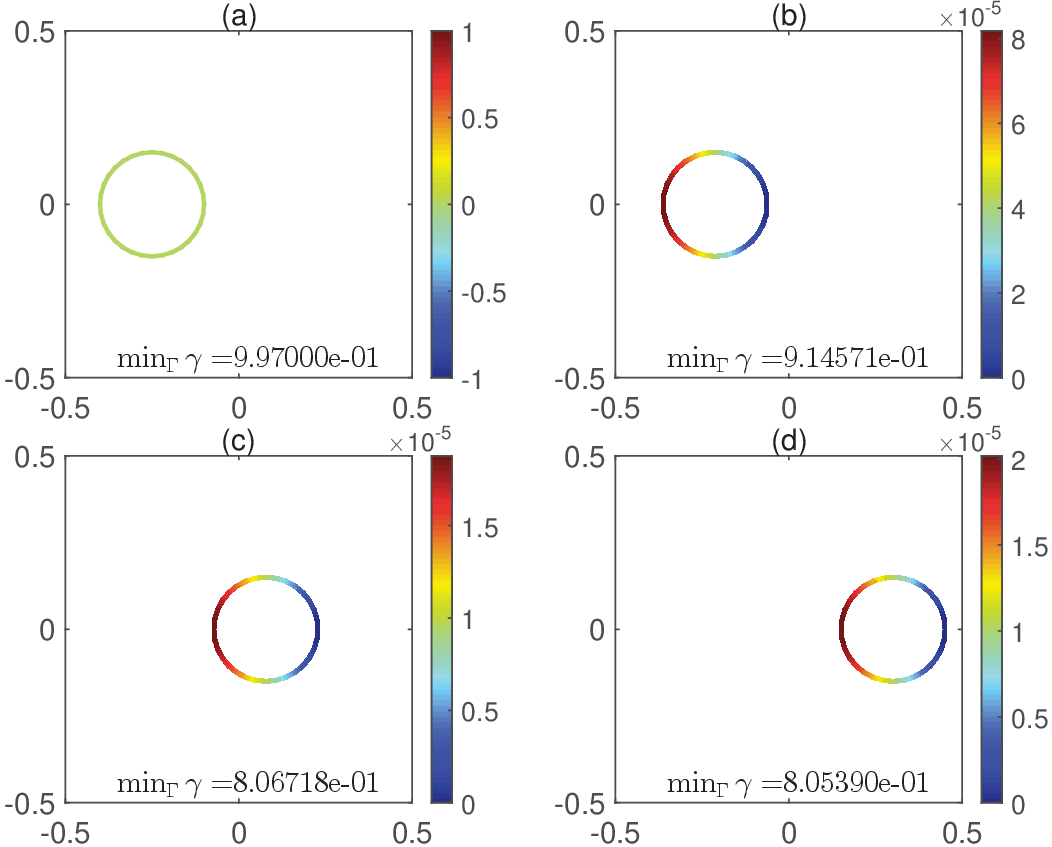}
    \caption{Snapshots of the surface tension $\gamma-\min_\Gamma \gamma$ at selected times. (a) $t=0$; (b) $t=10$; (c) $t=50$; (d) $t=100$.}
    \label{fig:var_weight_dirichlet_st}
\end{figure}

\begin{figure}[htb!]
    \centering
    \includegraphics[scale=0.45]{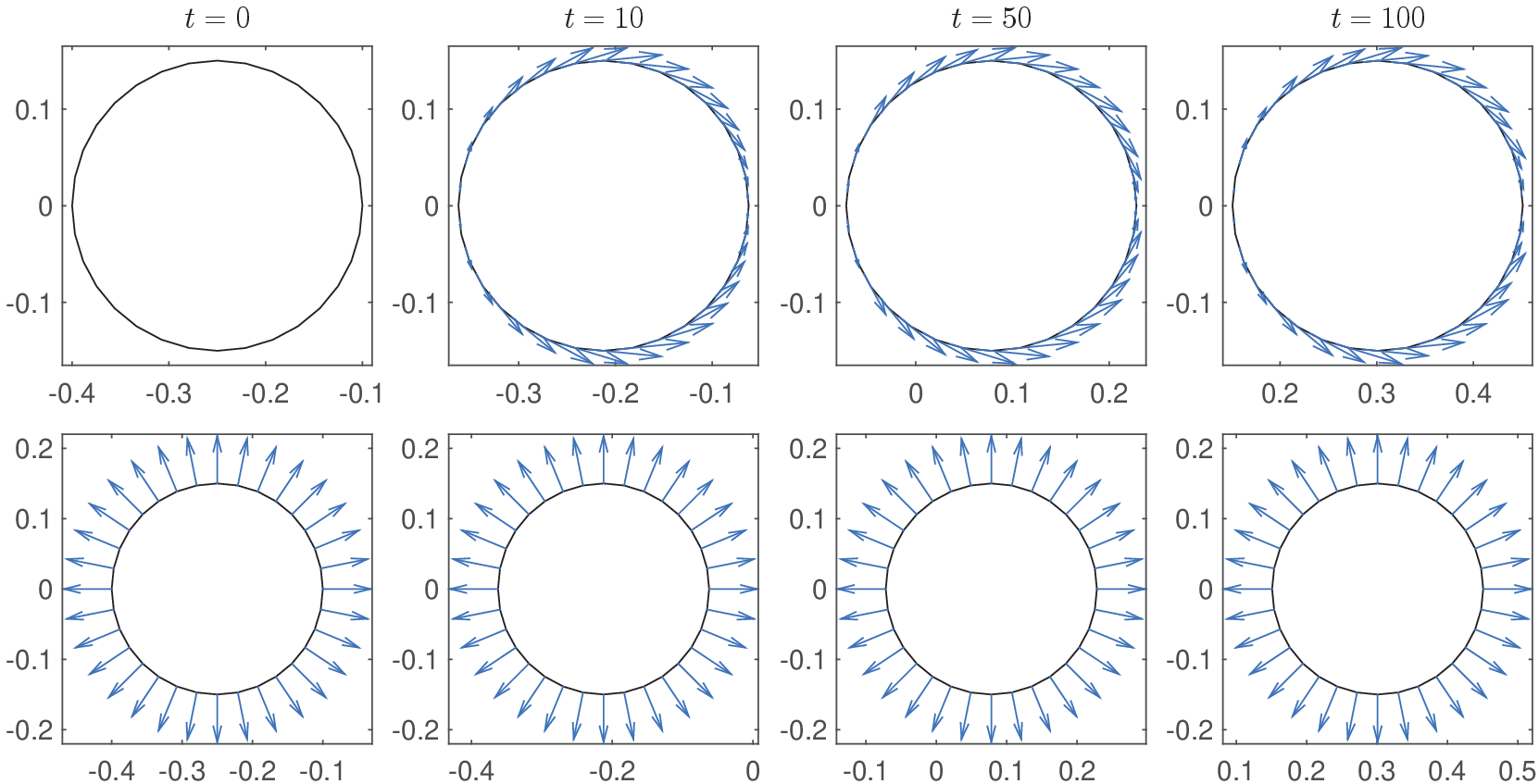}
    \caption{Snapshots of the interfacial force at selected times. Top panel: the Marangoni force; bottom panel: the Capillary force.}
    \label{fig:var_weight_dirichlet_force}
\end{figure}

Fig.~\ref{fig:var_weight_dirichlet_vel} presents the droplet interface and velocity field at selected time points. In the absence of interfacial chemical reactions, the droplet remains stationary. However, the presence of reactions generates interface-driven forces that induce droplet migration toward the region with higher bulk concentration \( C \).
The corresponding interfacial concentration of the reaction product \( A_\Gamma \) is shown in Fig.~\ref{fig:var_weight_dirichlet_sfa}. At early times, a larger amount of \( A_\Gamma \) accumulates on the right side of the droplet due to the higher availability of \( C \), leading to a local reduction in surface tension, as illustrated in Fig.~\ref{fig:var_weight_dirichlet_st}.
Although the capillary force \( \gamma \kappa \bm{n} \) remains symmetric along the interface, the tangential Marangoni force \( -\nabla_s \gamma \), shown in Fig.~\ref{fig:var_weight_dirichlet_force}, becomes asymmetric. This asymmetry drives the droplet to move rightward, toward regions of high bulk substance concentration.

\begin{figure}[htb!]
    \centering
    \includegraphics[scale=0.52]{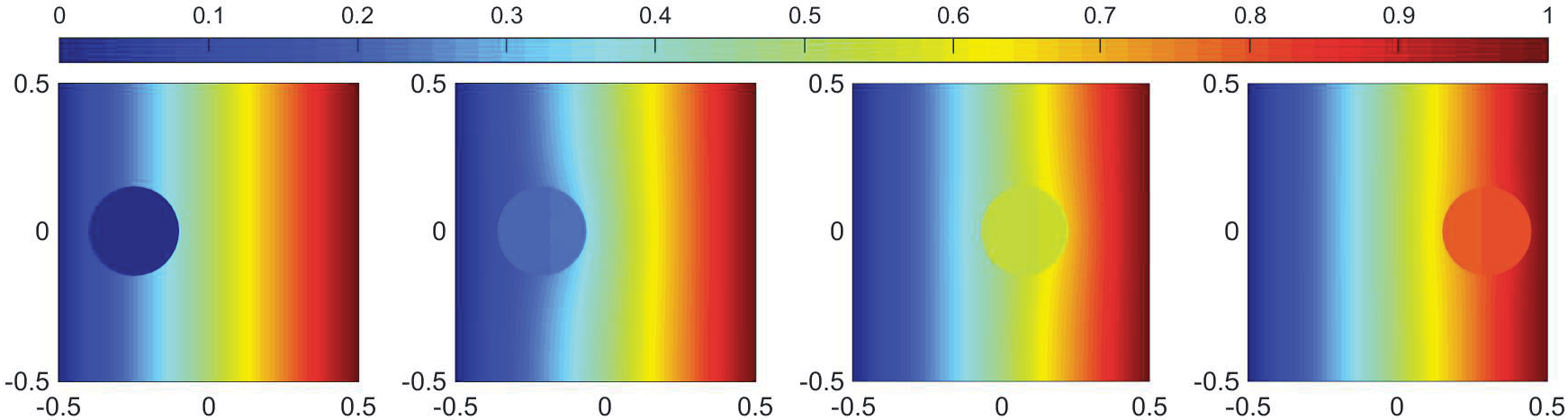}
    \caption{Snapshots of the bulk substance concentration $C$ at $t=0,10,50,100$.}
    \label{fig:var_weight_dirichlet_sb}
\end{figure}

The evolution of the bulk concentration \( C \) at times \( t = 0, 10, 50, 100 \) is shown in Fig.~\ref{fig:var_weight_dirichlet_sb}. Initially, the droplet contains no solute. As time progresses, solute gradually enters the droplet through the permeable interface. Meanwhile, interfacial chemical reactions lead to the emergence of a concentration jump of \( C \) across the interface.

\begin{figure}[htb!]
    \centering
    \includegraphics[scale=0.5]{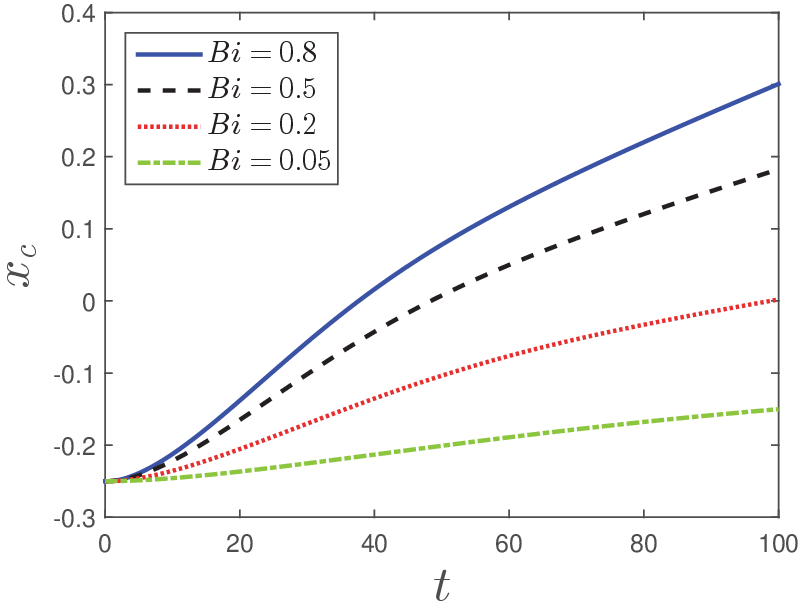}
    \includegraphics[scale=0.5]{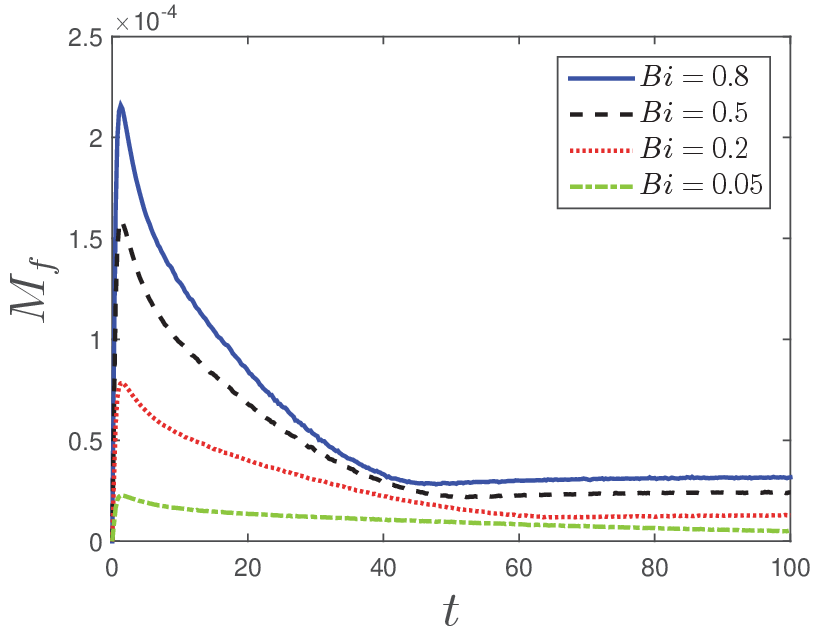}
    \caption{The evolution of the horizontal coordinate of the center of the drop (left panel) and the integral of the Marangoni force on the fluid interface (right panel) for different Biot numbers. }
    \label{fig:var_weight_dirichlet_biot_var}
\end{figure}

The influence of interface permeability is illustrated in Fig.~\ref{fig:var_weight_dirichlet_biot_var} a. For less permeable interfaces, the solute enters the droplet more slowly, which in turn slows down the interfacial reactions. This reduction in reaction rate diminishes the asymmetry of the Marangoni force along the x-axis $\int_{\Gamma}\nabla_s\gamma d\mathcal{H}^{d-1} \cdot (1,0)^T$ as shown in Fig.~\ref{fig:var_weight_dirichlet_biot_var} b.  As a result, the droplet experiences a weaker net driving force and exhibits slower motion.

%\subsection{Reaction-Induced Modulation of Permeability}
\subsection{Reaction-driven permeability gating: toward targeted delivery applications}
A compelling application of reactive, semi-permeable interfaces arises in the context of targeted drug delivery. In many therapeutic strategies, drug molecules must cross biological membranes, such as cellular membranes or vesicular encapsulations, to reach intracellular targets. These membranes typically exhibit selective permeability and can undergo chemical modifications in response to environmental cues, such as ligand binding or local enzymatic activity. To capture such processes, we consider a modified solute variable \( C_s \), governed in each bulk subdomain by the diffusion equation
\[
\partial_t^\bullet C_s^\pm = \frac{1}{Pe} \nabla \cdot \left( D_{C_s}^\pm \nabla C_s^\pm \right),
\]
subject to a reaction-modulated permeability condition on the interface \( \Gamma \):
\[
\frac{D_{C_s}^\pm}{Da \cdot Pe} \, {\bm n} \cdot \nabla C_s^\pm = Bi \cdot k_{C_s}(A_\Gamma) [\![ C_s ]\!].
\]

To represent biologically relevant gating behavior, we define the interfacial permeability function \( k_{C_s}(A_\Gamma) \) using a sigmoidal (logistic) form:
\[
k_{C_s}(A_\Gamma) = k_{\max} \cdot \frac{1}{1 + \exp(-\beta (A_\Gamma - A_0))},
\]
where \( k_{\max} \) denotes the maximum permeability, \( A_0 \) is the activation threshold, and \( \beta \) controls the sharpness of the transition. This formulation mimics switch-like responses in biological systems, such as ligand-gated ion channels or enzyme-sensitive drug carriers. Here, we use \( k_{\max} = 0.1 \), \( A_0 = 0.8 \), and \( \beta = 50 \), ensuring negligible leakage when \( A_\Gamma < 0.8 \), with sharp activation near \( A_\Gamma = 0.8 \).

Initially, \( C_s \) is confined to the exterior domain:
\[
C_s(\bm{x}, 0) = 
\begin{cases}
0, & \text{for } \bm{x} \in \Omega^+, \\
1, & \text{for } \bm{x} \in \Omega^-.
\end{cases}
\]

\begin{figure}[htb!]
    \centering
    \includegraphics[scale=0.52]{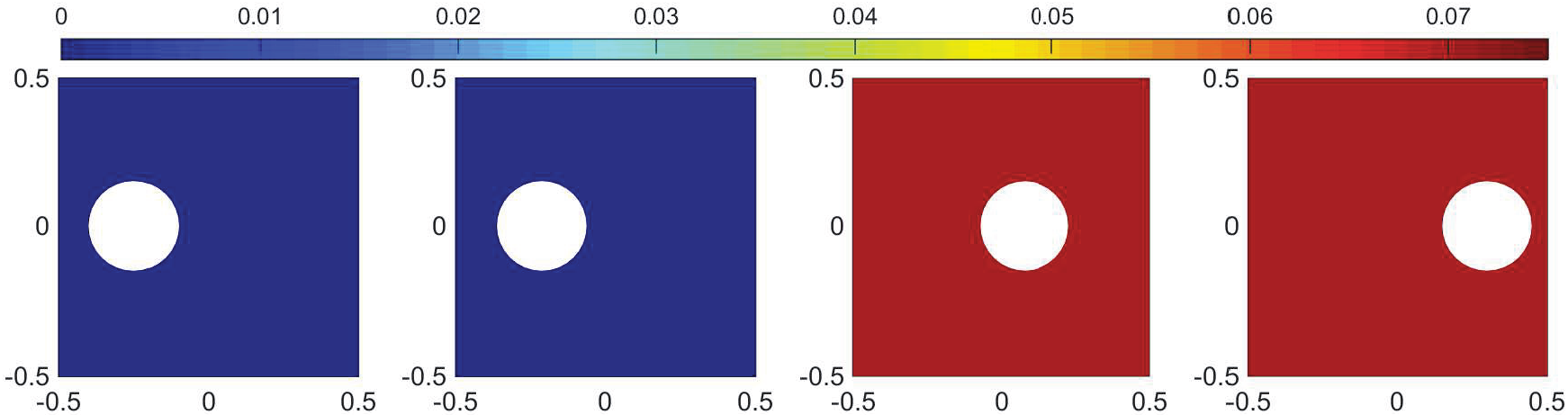}
    \caption{Snapshots of the bulk solute concentration \( C_s(\bm{x}, t) \) in \( \Omega^+ \) at $t=0,10,50,100$.}
    \label{fig:var_weight_dirichlet_sbs}
\end{figure}

Fig.~\ref{fig:var_weight_dirichlet_sbs} presents snapshots of the concentration distribution of \( C_s \) in \( \Omega^+ \) at times \( t = 0, 10, 50, 100 \). Initially, no leakage is observed due to low \( A_\Gamma \) and near-zero permeability. However, as the simulation progresses, the droplet self-propels toward regions with higher bulk solute concentration \( C \), driven by Marangoni stresses. These interfacial flows result from surface tension gradients induced by asymmetric production of \( A_\Gamma \), highlighting a chemically triggered, self-motile behavior.

\begin{figure}[htb!]
    \centering
    \includegraphics[scale=0.55]{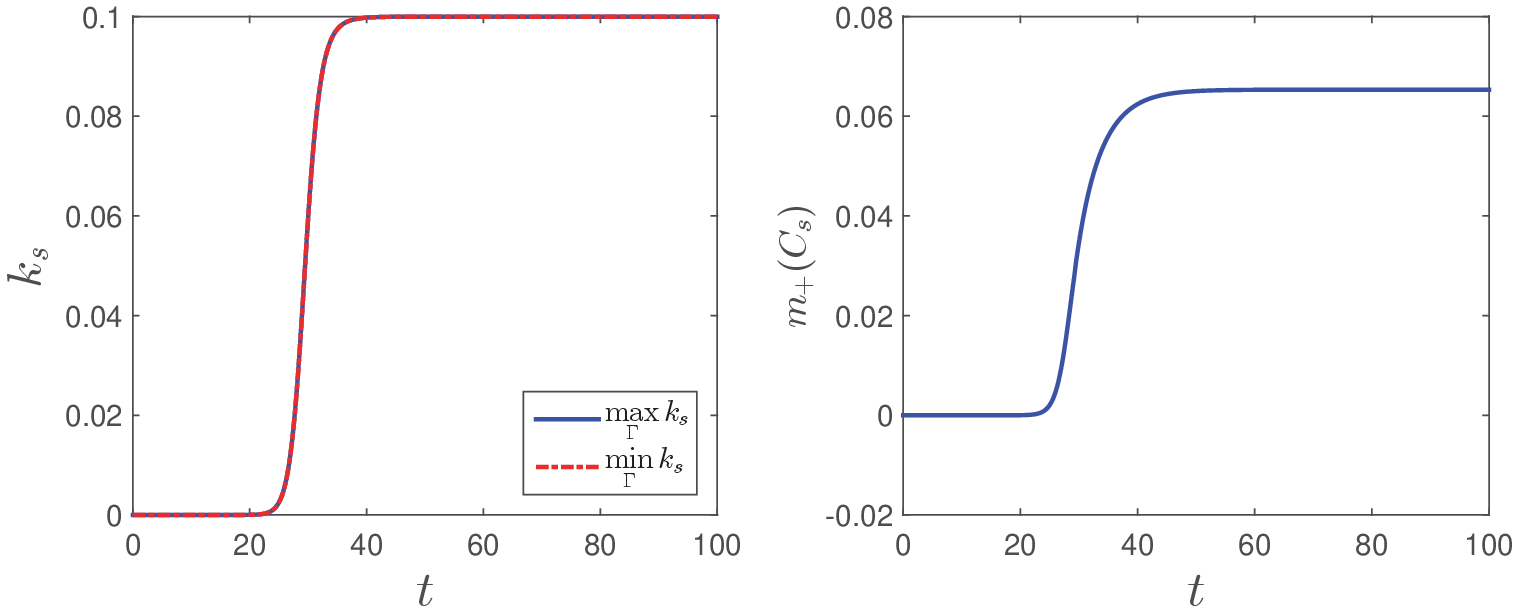}
    \caption{Left: Evolution of the maximum and minimum interfacial permeability \( k_{C_s} \). Right: Temporal evolution of solute mass \( m_+(C_s) \) in the droplet domain \( \Omega^+ \).}
    \label{fig:var_weight_dirichlet_ks_mass}
\end{figure}

Once \( A_\Gamma \) exceeds the activation threshold \( A_0 \) (around \( t = 20 \)), the permeability \( k_{C_s}(A_\Gamma) \) increases sharply (Fig.~\ref{fig:var_weight_dirichlet_ks_mass}, left), leading to a gating-like opening of the interface. The right panel of Fig.~\ref{fig:var_weight_dirichlet_ks_mass} shows the corresponding mass of \( C_s \) inside the droplet,
\[
m_+(C_s) = \int_{\Omega^+} C_s \, d\sL^{d},
\]
which rises monotonically until equilibrium is reached between the two sides. This progression captures the essence of targeted drug delivery: a vesicle remains sealed until reaching a reactive microenvironment, where chemical triggers induce a local permeability change, allowing the payload to be selectively released into the target tissue.

\begin{rmk} \label{rmk:leakyreaction}
The above model can be obtained by incorporating  by adding an additional term 
\[
\frac{1}{We \cdot Da} \int_{\Omega} f(C_s) d\sL^{d}
\]
to the total free energy, and a corresponding dissipation term
\[
\frac{1}{We} \int_{\Gamma} k_{C_s}(A_\Gamma) [\![f'(C_s) ]\!] [\![ C_s ]\!] \, d\mathcal{H}^{d-1}
\]
to the total dissipation functional.  
\end{rmk}
%For instance, surface-functionalized nanoparticles can release therapeutic agents upon encountering specific receptors or enzymatic activity at the interface, effectively triggering localized transmembrane transport. Moreover, certain drug carriers are designed with pH-sensitive or redox-reactive coatings that alter membrane permeability or promote membrane fusion upon activation. Such processes involve complex coupling between mass transport, interfacial reaction kinetics, and membrane mechanics. Accurate modeling of these mechanisms is essential for understanding delivery efficiency and optimizing design, especially in applications like cancer therapy, where microenvironmental heterogeneity strongly influences treatment outcomes.

\subsection{Reaction enhanced deformation under shear flow}\label{subsec:bubble_shear}

In this example, we investigate the influence of interfacial chemical reactions on droplet deformation under shear flow. The initial bulk concentration is set as \( C(\bm{x},0) = 0 \), and the initial interfacial concentrations are \( A_\Gamma(\bm{x},0) = 0.5 \), \( B_\Gamma(\bm{x},0) = 0 \), and \( C_\Gamma(\bm{x},0) = 0 \). The computational domain is defined as \(\Omega = [-0.75, 0.75] \times [-0.5, 0.5]\), with the initial velocity field \({\bm u}(\bm{x},0) = \bm{0} \).

Periodic boundary conditions are applied on the left and right boundaries. On the top and bottom boundaries, no-flux conditions are imposed for the chemical species, while Dirichlet boundary conditions are used for the velocity field to generate shear flow:
\[
{\bm n} \cdot \nabla C^+ = 0, \quad {\bm u} = (y, 0).
\]

The physical and numerical parameters used in this simulation are as follows:
\[
\begin{aligned}
&\rho^+ = 10, \quad \rho^- = 1, \quad \mu^+ = 10, \quad \mu^- = 1, \quad Re = 5, \quad Ca = 0.2, \quad Da = 1, \\
&Bi = 0.1,\quad  E = 0.4, \quad Pe = 1, \quad \omega_a =  \omega_b = \omega_c = 1, \quad \lambda_a = 1, \quad \lambda_c = 0.01,  \\
&k_r =  k_d^\pm = 1,\quad D_{A_\Gamma} = 100, \quad D_{B_\Gamma} = D_{C_\Gamma} = 1, \quad D^+ = 0.5, \quad D^- = 1.
\end{aligned}
\]

Additional parameter values are defined as needed within specific simulation contexts. The time step is set to \(\Delta t = 1 \times 10^{-3}\).

\begin{figure}[htb!]
    \centering
    \includegraphics[scale=0.6]{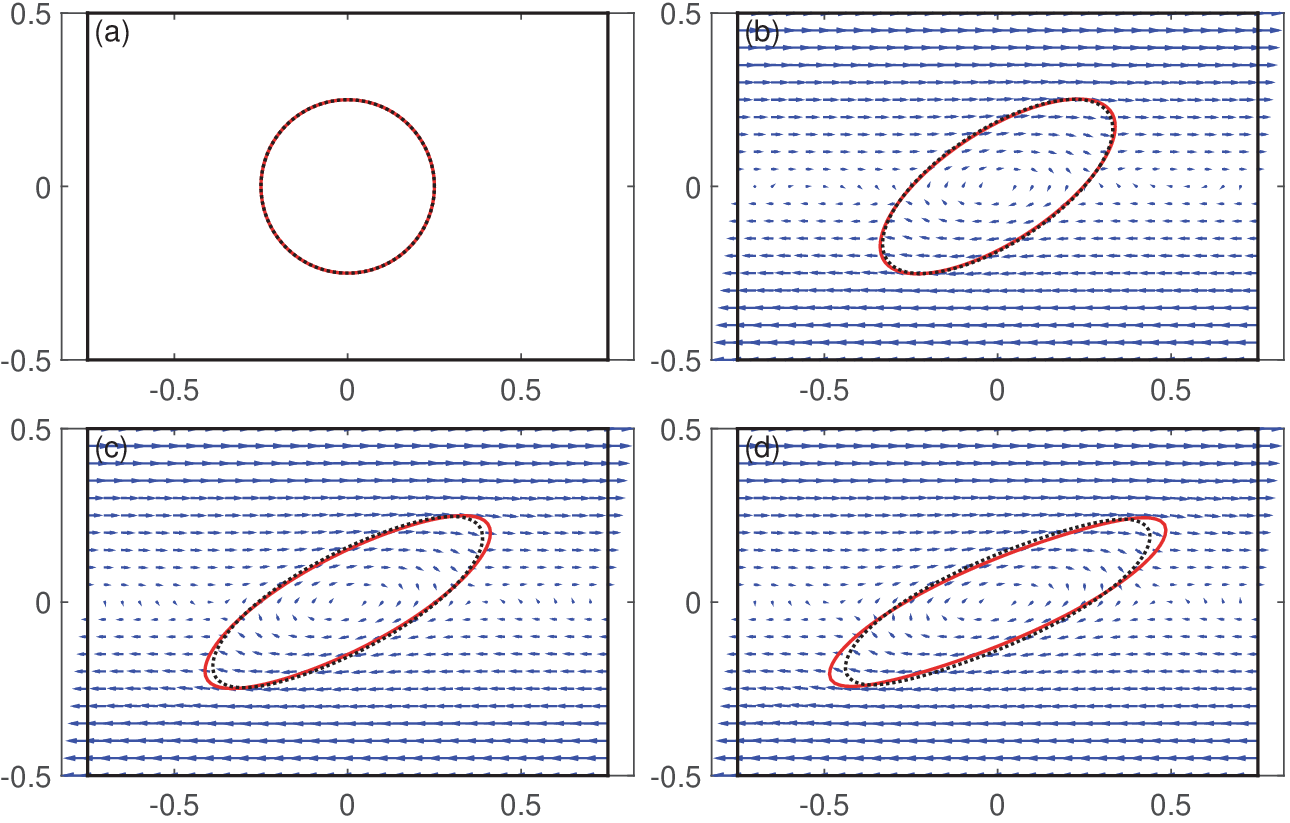}
    \caption{Snapshots of the fluid interface together with the velocity field at different times. Red solid line: \(k_r = 1\); black dotted line: \(k_r = 0\). (a) \(t = 0\); (b) \(t = 1\); (c) \(t = 1.5\); (d) \(t = 2\).}
    \label{fig:bubble1_shear_vel}
\end{figure}

Fig.~\ref{fig:bubble1_shear_vel} shows snapshots of the fluid interface together with the velocity at various time points. For comparison, we also show results for the case without chemical reaction (\(k_r = 0\)). The initial interface is a circle: \(\Gamma(0) = \{(x, y) \mid x^2 + y^2 = 0.25^2\}\), and \( Pe_\Gamma = 100 \). The interface is represented using \(J_\Gamma = 64\) markers. At \(t = 0\), the mesh contains \(K_\Omega = 2363\) vertices and \(J_\Omega = 4984\) triangles.

\begin{figure}[htb!]
    \centering
    \includegraphics[scale=0.6]{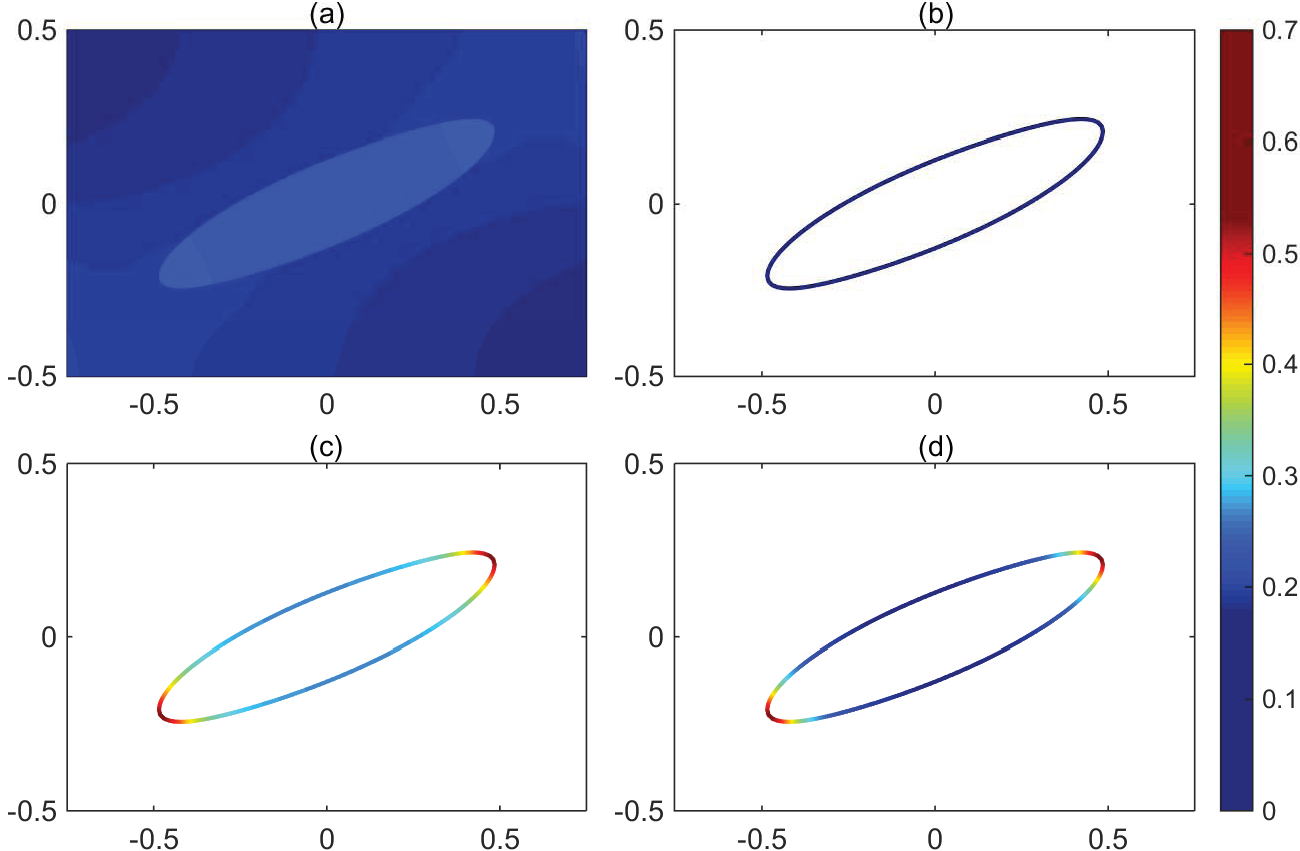}
    \caption{Concentration snapshots at \(t = 2\): (a) bulk \(C\); (b) interfacial \(A_\Gamma\); (c) interfacial \(B_\Gamma\); (d) interfacial \(C_\Gamma\).}
    \label{fig:bubble1_shear_sf_t2d0}
\end{figure}

Compared to the case without reaction, the deformation is significantly enhanced when \(k_r = 1\). This behavior can be attributed to two primary factors:
\begin{itemize}
\item \textbf{Increased total interfacial substance concentration due to chemical reaction.} Each molecule of \(C_{\Gamma}\) decomposes into one molecule each of \(A_{\Gamma}\) and \(B_{\Gamma}\). As a result, the total concentration of interfacial substances increases, which in turn reduces the surface tension and facilitates greater deformation.

\item \textbf{Convection-enhanced redistribution of interfacial substances.} With a Péclet number \(Pe_\Gamma = 10\), interfacial transport is convection-dominated. The imposed shear flow drives strong tangential convection along the interface (see Fig.~\ref{fig:bubble1_shear_vel} ), transporting interfacial substances toward the extremities of the elongated droplet (See   Fig.~\ref{fig:bubble1_shear_sf_t2d0}). This accumulation further reduces the local surface tension, thereby promoting elongation and reinforcing the deformation.
\end{itemize}

\subsection{The Rising Bubble} \label{subsec:bubble_rise}

In this subsection, we investigate the rising motion of a bubble under gravity, influenced by interfacial chemical reactions. The initial bulk concentration is set to \( C(\bm{x}, 0) = 0 \), while the interfacial concentrations are initialized as \( A_\Gamma(\bm{x}, 0) = 0.5 \), \( B_\Gamma(\bm{x}, 0) = 0 \), and \( C_\Gamma(\bm{x}, 0) = 0 \). The computational domain is \( \Omega = [-0.5, 0.5] \times [-0.5, 1.5] \), with the initial interface defined by \( \Gamma(0) = \{ (x, y) \mid x^2 + y^2 = 0.25^2 \} \). The initial velocity field is \( \bm{u}(\bm{x}, 0) = \bm{0} \).

On the boundary \( \partial\Omega \), we impose a no-flux condition \( \bm{n} \cdot \nabla C^+ = 0 \) for the bulk concentration. For the velocity field, Dirichlet conditions \( \bm{u} = \bm{0} \) are prescribed on \( \partial\Omega_1 \), and slip conditions \( \bm{u} \cdot \bm{n} = 0 \) on \( \partial\Omega_2 \).

The physical and numerical parameters are:
\[
\begin{aligned}
&\rho^+ = 10, \quad \rho^- = 1, \quad \mu^+ = 25, \quad \mu^- = 1, \quad Re = 100, \quad Ca = \frac{1}{24.5}, \quad Da = 1,  \\
&Bi = 0.1,\quad E = 0.4, \quad Pe_\Gamma = 10, \quad Pe = 1, \quad \omega_a = \omega_b =  \omega_c = 1, \quad \lambda_a = 1, \quad \lambda_c = 0.01, \\
&k_r = k_d^\pm = 1, \quad D_{A_\Gamma} = 1000, \quad D_{B_\Gamma} = D_{C_\Gamma} = 1, \quad D^+ = 0.5, \quad D^- = 1.
\end{aligned}
\]

The interface is represented using \( J_\Gamma = 64 \) markers, and a fixed time step \( \Delta t = 1 \times 10^{-3} \) is used. Unlike the previous simulations, gravitational forces are included in the momentum equations as \( \rho^\pm \bm{g} = -0.98 \rho^\pm \bm{e}_d \), where \( \bm{e}_d = (0, 1) \).

\begin{figure}[htb!]
    \centering
    \includegraphics[scale=0.54]{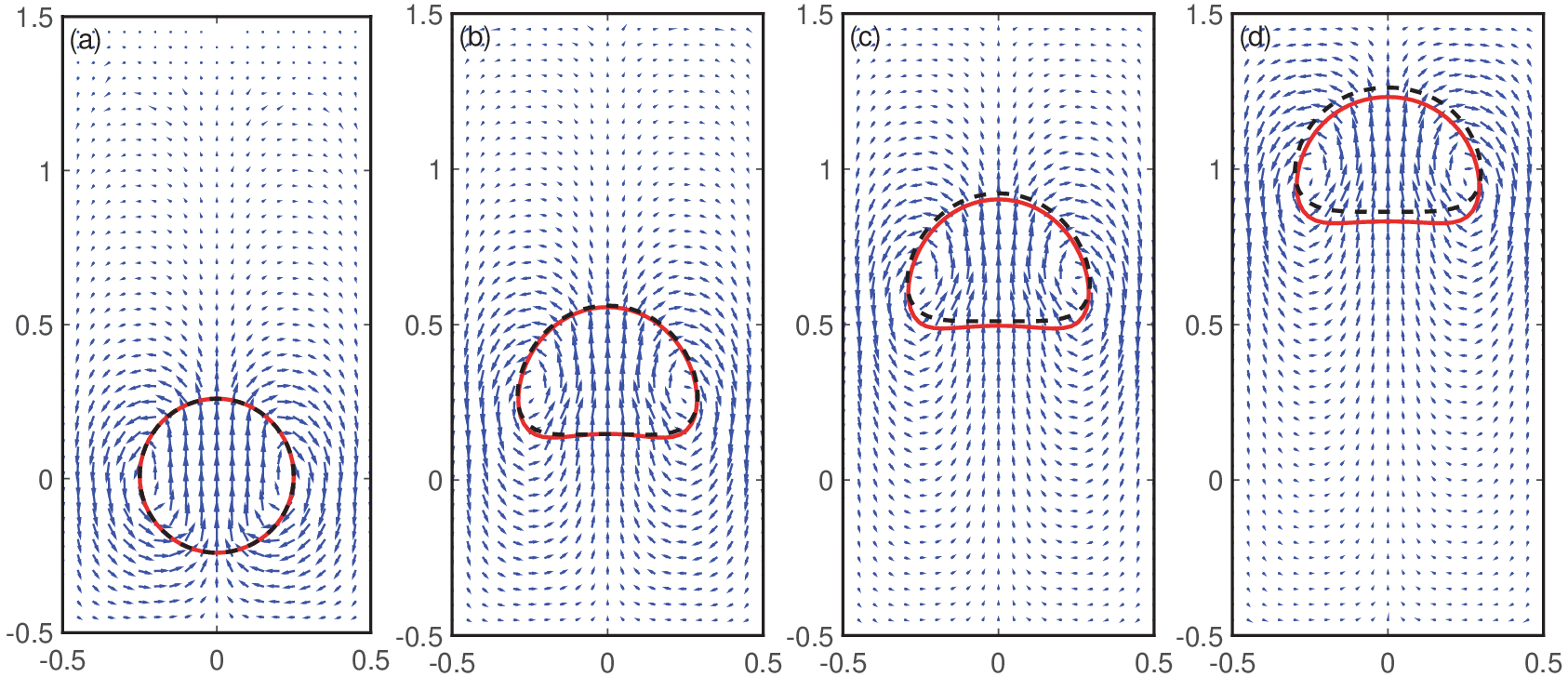}
    \caption{Snapshots of the fluid interface and velocity field at different times. Red solid line: with reaction; black dashed line: without reaction. (a) \( t=0.2 \), \( \max_{\bm{x} \in \Omega} |\bm{u}| = 0.1174 \), \( K_\Omega = 913 \), \( J_\Omega = 1956 \); (b) \( t=2.0 \), \( \max |\bm{u}| = 0.3354 \), \( K_\Omega = 927 \), \( J_\Omega = 1984 \); (c) \( t=4.0 \), \( \max |\bm{u}| = 0.3250 \), \( K_\Omega = 927 \), \( J_\Omega = 1984 \); (d) \( t=6.0 \), \( \max |\bm{u}| = 0.3042 \), \( K_\Omega = 873 \), \( J_\Omega = 1876 \).}
    \label{fig:bubble1_rising_velfluid}
\end{figure}

Fig.~\ref{fig:bubble1_rising_velfluid} shows snapshots of the fluid interface and velocity field at various times, comparing results with and without chemical reaction (red solid line: \( k_r = 1 \) vs. black dash line: \( k_r = 0 \)). The presence of the chemical reaction noticeably slows the upward motion of the bubble.   Mesh regeneration for the bulk is applied as needed throughout the simulation.

\begin{figure}[htb!]
    \centering
    \includegraphics[scale=0.52]{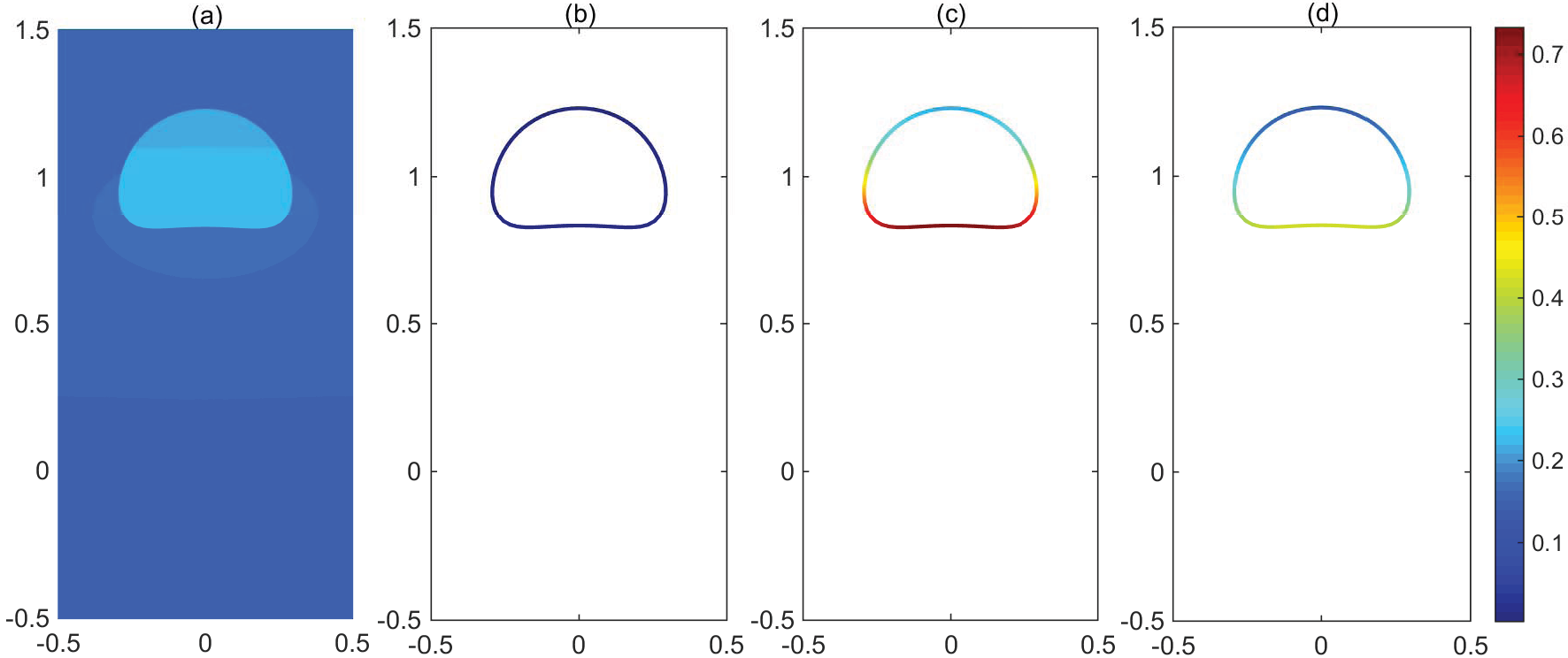}
    \caption{Snapshots of concentration fields at \( t = 6 \) for the interface with reaction. (a) Bulk concentration \( C \); (b) Interfacial concentration \( A_\Gamma \); (c) \( B_\Gamma \); (d) \( C_\Gamma \).}
    \label{fig:bubble1_rising_sb}
\end{figure}

Fig.~\ref{fig:bubble1_rising_sb} presents snapshots of the concentration fields at \( t=6 \). Compared to their initial values, \( B_\Gamma \), \( C_\Gamma \), and \( C \) have increased, while \( A_\Gamma \) has decreased, reflecting the progress of the interfacial reaction and the dissociation from the interface to the bulk region.

To quantitatively assess the effect of chemical reactions, we introduce the following discrete diagnostic quantities:
\begin{equation*}
 C_d|_{t_m} := \frac{\pi^{1/d} [2d~\text{vol}(\Omega^{-,m})]^{\frac{d-1}{d}}}{|\Gamma^m|}, \quad 
 y_c|_{t_m} := \frac{\int_{\Omega^{-,m}} (\mathbb{I} \cdot \bm{e}_d) \, d\sL^{d}}{\text{vol}(\Omega^{-,m})}, \quad 
 V_c|_{t_m} := \frac{\int_{\Omega^{-,m}} (\bm{u}^m \cdot \bm{e}_d) \, d\sL^{d}}{\text{vol}(\Omega^{-,m})},
\end{equation*}
where \( C_d \) measures bubble circularity, \( y_c \) denotes the center of mass in the vertical direction, and \( V_c \) indicates the average rise velocity.

\begin{figure}[htb!]
    \centering
    \includegraphics[scale=0.7]{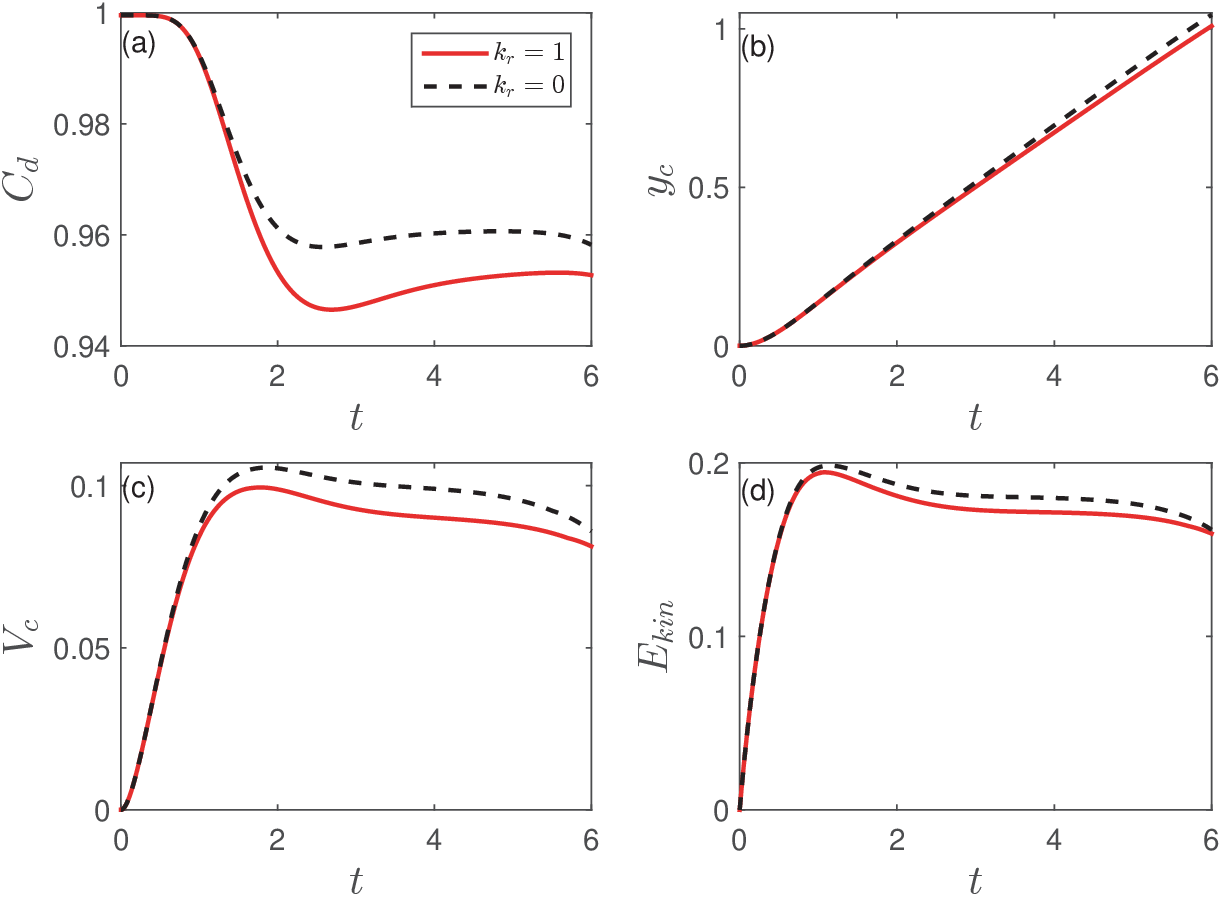}
    \caption{Time evolution of key quantities: (a) bubble circularity; (b) vertical center of mass; (c) rise velocity; (d) total kinetic energy. Black dashed line: without reaction; Red solid line: with reaction. }
    \label{fig:bubble1_rising_para}
\end{figure}

Fig.~\ref{fig:bubble1_rising_para} plots the time histories of these quantities along with the total kinetic energy of the system. In particular, panel (a) shows that interfacial reaction increases bubble deformation due to the decrease of surface tension. Panels (b) and (c) reveal that the center of mass and rise velocity are similar in both cases during the early stage, but diverge  later on, consistent with the trends observed in Fig.~\ref{fig:bubble1_rising_velfluid}. 
Our simulations reveal that interfacial chemical reactions affect the rising dynamics of the bubble. Specifically, in the presence of chemical reactions, the bubble undergoes more pronounced deformation and exhibits a reduced rising velocity compared to the non-reactive case. This phenomenon can be attributed to the production of surface-active species   through the interfacial reaction, which lowers the local surface tension. A reduced surface tension effectively softens the bubble interface, making it more susceptible to deformation under hydrodynamic stresses. The resulting increased deformation alters the flow field around the bubble, increasing the hydrodynamic drag. Consequently, the upward motion of the bubble is slowed down despite identical gravitational forcing. This mechanism is consistent with previous studies on surfactant-laden bubbles, where Marangoni effects and surface tension reduction have been shown to suppress bubble rise velocity through enhanced interfacial deformation and asymmetry \cite{sadhal1983stokes,magnaudet2000motion}.

\section{Conclusion}\label{sec:concluding}
In this work, we introduce a unified, thermodynamically consistent continuum model for fluid–structure interactions at semi-permeable interfaces with reaction dynamics. By systematically coupling interfacial deformation, fluid flow, surface and bulk transport, reaction kinetics, and selective permeability, the model captures the intricate multiphysics underlying a wide range of interfacial phenomena. The derivation, grounded in an energy variation approach, ensures the laws of mass conservation and energy dissipation.

 To address the computational challenges posed by evolving interfaces and tightly coupled surface–bulk dynamics, we develop a finite element scheme based on an ALE formulation. The fully coupled system is expressed in a weak form and is decomposed into two interacting subsystems. The first governs the fluid dynamics, assuming the solute concentrations are prescribed. The second describes the transport, adsorption, desorption, and chemical reactions of solute species within the bulk and along the interface. To ensure numerical stability and mesh integrity, the scheme incorporates the BGN framework, which preserves the geometric fidelity of the interface while maintaining the conservation properties intrinsic to the model. This approach enables efficient and accurate simulation of complex interfacial phenomena.

Through a series of numerical experiments, we validated the convergence and conservation properties of the scheme and explored complex interfacial phenomena. In particular, we demonstrated how interfacial reactions can create asymmetric surface tension distributions, giving rise to Marangoni-driven droplet motion. This chemically induced self-motility was shown to be highly sensitive to interfacial reaction kinetics and local concentration gradients. Furthermore, we modeled a switch-like increase in membrane permeability triggered by the accumulation of interfacial species, enabling gated solute transport across the droplet interface.
This mechanism was illustrated through a biologically inspired simulation, in which a droplet carrying solute molecules self-propels toward a reactive region and selectively releases its contents upon reaching a biochemical threshold. The interplay between Marangoni flows and reaction-activated permeability provides a mechanistic analog for targeted drug delivery systems, such as ligand- or enzyme-sensitive vesicles, which remain sealed during transport and only release their cargo after sensing specific molecular cues. Another biologically motivated application, we extended the model to simulate the ABCG1-mediated cholesterol efflux pathway. This process involves a cascade of reactions and interfacial transports where cholesterol is exported from the cell, transferred across the extracellular space, and ultimately absorbed by high-density lipoproteins (HDL). Our model captures each stage of this pathway: interfacial binding, transmembrane release, extracellular transport, and HDL uptake, demonstrating its capability to represent physiologically relevant, multiscale transport dynamics. The results reproduce key features of ABCG1-regulated cholesterol trafficking and illustrate the potential of this framework for modeling lipid transport in vascular and metabolic biology.

 In summary, our results underscore the importance of coupling chemical kinetics, transport mechanisms, and mechanical interface dynamics in a thermodynamically consistent way. The model serves as a versatile and extensible computational tool for investigating reactive interface problems across biological and industrial contexts, with direct applications in drug delivery, lipid metabolism, smart emulsions, and bioinspired microfluidic systems. 

%In general, our results highlight the critical role of a thermodynamically consistent coupling between chemical kinetics, transport phenomena, and interfacial dynamics in the modeling of reactive fluid-interface systems. The framework provides a robust and extensible foundation for studying biological membranes, smart emulsions, and microfluidic devices. 
%For example, the model can be applied to investigate the role of reactive oxygen species (ROS) produced in response to mechanical or chemical stimuli. ROS participate in redox reactions that alter protein conformations and lipid organization, thus modifying membrane tension, stiffness,  and permeability of ions across the membrane~\cite{zhou2010reaction,zeng2023mathematical}. 

\section*{Acknowledgment} 
 
The research of S.~Xu is partially supported by the National Natural Science Foundation of China (Grant No.~12271492). The work of W.~Shi is partially supported by the Scientific and Technological Innovation Programs of Higher Education Institutions in Shanxi (Grant No.~2019L0474). Z.~Zhang acknowledges partial support from the National Key R\&D Program of China (Grant No.~2023YFA1011403), the National Natural Science Foundation of China (Grant No.~92470112), the Shenzhen Sci-Tech Innovation Commission (Grant No.~20231120102244002), and the Guangdong Provincial Key Laboratory of Computational Science and Material Design (Grant No.~2019B030301001). The research of Q.~Zhao is partially supported by the National Natural Science Foundation of China (Grant No.~12401572).

%%%%%%%%%%%%%%%%%%%%% References %%%%%%%%%%%%%%%%%%%%%%%%%
\footnotesize
\bibliographystyle{abbrv}
\bibliography{myref}

\newpage
\begin{appendices}

\section{Integral Calculus}\label{sec:Appendix_Integral_Calculus}

% \begin{lma}(Reynolds Transport Formula) Let $C$ be any differentiable function defined in an evolving domain $\Omega(t)$, then
% \begin{equation}
%   \label{eq:Reynolds_Transport_Formula}
%   \mathrm{  \mathrm{\frac{d}{dt}}}\int_{\Omega(t)} C({\bm x},t) d\sL^{d} = \int_{\Omega(t)}\big(\partial^\bullet_t C+\nabla\cdot{\bm u}C\big) d\sL^{d}. 
% \end{equation}
% \end{lma}

% \begin{lma} (Bulk Divergence Formula) Let ${\bm F}$ be any differentiable vector function defined in a domain $\Omega$, then
% \begin{equation}
%   \label{eq:Bulk_Divergence_Formula1}
%   \int_{\Omega} \nabla\cdot {\bm F} d\sL^{d} = \int_{\partial\Omega} {\bm F}\cdot{\bm n} d\mathcal{H}^{d-1},
% \end{equation}
% where ${\bm n}$ is the unit outward normal vector of $\partial{\Omega}$.
% \end{lma}

% \begin{cor}
%   Using the identity $\nabla\cdot(C {\bm F}) = \nabla C \cdot {\bm F} + C \nabla\cdot {\bm F}$ and Eq. \eqref{eq:Bulk_Divergence_Formula1}, we have 
%   \begin{equation}
%     \label{eq:Bulk_Divergence_Formula2}
%     \int_{\Omega} C \nabla\cdot{\bm F} d\sL^{d} = \int_{\partial\Omega} C {\bm n}\cdot{\bm F} d\sL^{d} - \int_{\Omega} {\bm F}\cdot\nabla C d\sL^{d} . 
%   \end{equation}
% \end{cor}

\begin{lma} (Conservation of Bulk Mass) Let $C$ be the concentration function of a certain substance and defined in a domain $\Omega$, then
\begin{equation}
  \label{eq:Bulk_Mass_Conservation_Integral}
  \mathrm{  \mathrm{\frac{d}{dt}}}\int_{\hat{\Omega}} C d\sL^{d} + \int_{\partial\hat{\Omega}} {\bm J}_C\cdot\hat{\bm n} d\mathcal{H}^{d-1} = 0,
\end{equation}
where $\hat{\Omega}$ is any sub-domain of $\Omega$, and ${\bm J}_C$ is the diffusion flux and $\hat{\bm n}$ is the unit outward normal vector of $\partial\hat{\Omega}$. Applying the Reynolds transport formula and the   divergence formula, we can reformulate the left-hand side of Eq. \eqref{eq:Bulk_Mass_Conservation_Integral} as 
\begin{equation}
  \label{eq:Bulk_Mass_Conservation_Integra2}
  \mathrm{  \mathrm{\frac{d}{dt}}}\int_{\hat{\Omega}} C d\sL^{d} + \int_{\partial\hat{\Omega}} {\bm J}_C\cdot\hat{\bm n} d\mathcal{H}^{d-1}=\int_{\hat{\Omega}}\bigg( \frac{\partial C}{\partial t} + \nabla\cdot(C{\bm u}) + \nabla\cdot{\bm J}_C \bigg) d\sL^{d}.
\end{equation}
Since the above equation holds in any $\hat{\Omega}\subset\Omega$, we directly obtain the differential form of the conservation law for the substance in the bulk domain $\Omega$,
\begin{equation}
 \label{eq:Bulk_Mass_Conservation_Diff1}
  \frac{\partial C}{\partial t} + \nabla\cdot(C{\bm u}) + \nabla\cdot{\bm J}_C = 0.
\end{equation}
When ${\bm u}$ is a divergence-free vector, the above equation can be simplified as 
\begin{equation}
 \label{eq:Bulk_Mass_Conservation_Diff2}
 \partial^\bullet_t C + \nabla\cdot{\bm J}_C = 0.
\end{equation}
\end{lma}

\begin{lma}
  Let $C$ be a function defined in an evolving domain $\Omega(t)$ and satisfy the conservation law \eqref{eq:Bulk_Mass_Conservation_Diff1}, and $f(C)$ be an energy density function, then 
  \begin{align}
    \mathrm{  \mathrm{\frac{d}{dt}}}\int_{\Omega} f(C) d\sL^{d} &= \int_{\Omega}f^{\prime\prime}(C)\nabla C\cdot{\bm J}_C d\sL^{d} -\int_{\partial\Omega} f^{\prime}(C){\bm n}\cdot{\bm J}_C d\mathcal{H}^{d-1} + \int_{\Omega}(f(C)- Cf^{\prime}(C))\nabla\cdot{\bm u} d\sL^{d}. \label{eq:Bulk_Energy_Derivative1}
  \end{align}
  Specifically, when ${\bm u}$ is a divergence-free vector, Eq. \eqref{eq:Bulk_Energy_Derivative1} reduces to 
  \begin{equation}
    \mathrm{  \mathrm{\frac{d}{dt}}}\int_{\Omega} f(C) d\sL^{d} = \int_{\Omega} f^{\prime\prime}(C)\nabla C\cdot{\bm J}_C d\sL^{d} -\int_{\partial\Omega} f^{\prime}(C){\bm n}\cdot{\bm J}_C d\mathcal{H}^{d-1}. \label{eq:Bulk_Energy_Derivative2}
  \end{equation}
\end{lma}
\begin{proof} 
  The details of the derivation of Eq. \eqref{eq:Bulk_Energy_Derivative1} are as follows
  \begin{align}
  \mathrm{  \mathrm{\frac{d}{dt}}}\int_{\Omega} f(C) d\sL^{d} &= \int_{\Omega}\big[ f^{\prime}(C)\partial^\bullet_t C + f(C)\nabla\cdot{\bm u}\big] d\sL^{d}  = \int_{\Omega}\big[(f(C)-f^{\prime}(C)C)\nabla\cdot{\bm u} - f^{\prime}(C)\nabla\cdot{\bm J}_C \big] d\sL^{d} \nonumber
  \\
  & = \int_{\Omega} f^{\prime\prime}(C)\nabla C\cdot{\bm J}_C d\sL^{d} -\int_{\partial\Omega} f^{\prime}(C){\bm n}\cdot{\bm J}_C d\mathcal{H}^{d-1} + \int_{\Omega}(f(C)- Cf^{\prime}(C))\nabla\cdot{\bm u} d\sL^{d}, \label{eq:Bulk_Energy_Derivative3}
  \end{align}
  where we have used  the Reynolds transport formula and Eq.\eqref{eq:Bulk_Mass_Conservation_Diff1}.
\end{proof}

% \begin{lma} (Surface Transport Formula) Let $C_\Gamma$ be any differentiable function defined on an evolving surface $\Gamma(t)$, then
% \begin{equation}
% \label{eq:Surface_Transport_Formula}
%    \mathrm{\frac{d}{dt}}\int_{\Gamma(t)} C_\Gamma({\bm x},t) d\mathcal{H}^{d-1} = \int_{\Gamma(t)}\big(\partial^\bullet_t C_\Gamma+\nabla_s\cdot{\bm u}C_\Gamma\big) d\mathcal{H}^{d-1}. 
% \end{equation}
% \end{lma}

% \begin{lma} (Surface Divergence Formula) Let ${\bm F}_\Gamma$ be any differentiable vector function defined on a surface $\Gamma$, then
% \begin{equation}
%   \label{eq:Surface_Divergence_Formula1} 
%   \int_{\Gamma} \nabla_s\cdot {\bm F}_\Gamma d\mathcal{H}^{d-1} = \int_{\Gamma} {\bm F}_\Gamma\cdot{\bm n} \kappa~d\mathcal{H}^{d-1} + \int_{\partial\Gamma} {\bm F}_\Gamma\cdot{\bm m} d\mathcal{H}^{d-2},
% \end{equation}
% where ${\bm n}$ is the unit outward normal vector of $\Gamma$, and ${\bm m}$ is the unit outward conormal vector of ${\Gamma}$ at $\partial{\Gamma}$.
% \end{lma}

\begin{cor}
  Using the identity $\nabla_s\cdot(C_\Gamma {\bm F}_{\Gamma}) = \nabla_s C_\Gamma \cdot {\bm F}_{\Gamma} + C_\Gamma \nabla_s\cdot {\bm F}_{\Gamma}$ and the surface divergence formula, we have 
  \begin{equation}
  \label{eq:Surface_Divergence_Formula2}
    \int_{\Gamma} C_\Gamma \nabla_s\cdot{\bm F}_{\Gamma} d\mathcal{H}^{d-1} = \int_{\Gamma} C_\Gamma \kappa {\bm n}\cdot{\bm F}_{\Gamma} d\mathcal{H}^{d-1} -\int_{\Gamma} {\bm F}_{\Gamma}\cdot\nabla_s C_\Gamma d\mathcal{H}^{d-1} + \int_{\partial\Gamma} C_\Gamma {\bm F}_{\Gamma}\cdot{\bm m} d\mathcal{H}^{d-2}. 
  \end{equation}
\end{cor}

\begin{lma} (Conservation of Surface Mass)  Let $C_\Gamma$ be the concentration function of a certain substance and defined on a surface $\Gamma$, then
  \begin{equation}
  \label{eq:Surface_Mass_Conservation_Integral}
  \mathrm{  \mathrm{\frac{d}{dt}}}\int_{\hat{\Gamma}} C_\Gamma d\mathcal{H}^{d-1} + \int_{\partial\hat{\Gamma}} {\bm J}_{C_\Gamma}\cdot\hat{\bm m} d\mathcal{H}^{d-2} = \int_{\hat{\Gamma}} S d\mathcal{H}^{d-1},
  \end{equation}
where $\hat{\Gamma}$ is any subset of $\Gamma$, ${\bm J}_{C_\Gamma}$ is the surface diffusion flux tangential to $\hat{\Gamma}$ and $\hat{\bm m}$ is the unit outward conormal vector of $\hat{\Gamma}$ at $\partial\hat{\Gamma}$. Here, $S$ is a source for the surface concentration. Similarly, applying the surface transport formula and the surface divergence formula, we can reformulate the left-hand side of Eq. \eqref{eq:Surface_Mass_Conservation_Integral} as 
\begin{align}
    &  \mathrm{\frac{d}{dt}}\int_{\hat{\Gamma}} C_\Gamma d\mathcal{H}^{d-1} + \int_{\partial\hat{\Gamma}} {\bm J}_{C_\Gamma}\cdot\hat{\bm m} d\mathcal{H}^{d-2}\nonumber
    \\
   =&\int_{\hat{\Gamma}}\big( \partial^\bullet_t C_\Gamma +\nabla_s\cdot{\bm u}C_\Gamma + \nabla\cdot{\bm J}_{C_\Gamma} \big) d\mathcal{H}^{d-1}-\int_{\hat{\Gamma}}\kappa {\bm n}\cdot{\bm J}_{C_\Gamma}d\mathcal{H}^{d-1}.\label{eq:Surface_Mass_Conservation_Integra2}
\end{align}
Since ${\bm n}\cdot{\bm J}_{C_\Gamma}=0$ and $\hat{\Gamma}$ is arbitrary, we immediately obtain the differential form of the conservation law for the substance on the surface $\Gamma$,
\begin{equation}
\partial^\bullet_t C_\Gamma +\nabla_s\cdot{\bm u}C_\Gamma + \nabla\cdot{\bm J}_{C_\Gamma} = S.
\label{eq:Surface_Mass_Conservation_Diff1}
\end{equation}
 
\end{lma}

\begin{lma}
  Let $C_\Gamma$ be a function defined on an evolving surface $\Gamma$ and satisfy the conservation law \eqref{eq:Surface_Mass_Conservation_Diff1}, and $g(C_\Gamma)$ be an energy density function. Then 
  \begin{align}
    \mathrm{  \mathrm{\frac{d}{dt}}}\int_{\Gamma} g(C_\Gamma) d\mathcal{H}^{d-1} &= \int_{\Gamma} \big\{{\bm u}\cdot(\tilde{\gamma}(C_\Gamma)\kappa{\bm n}-\nabla_s\tilde{\gamma}(C_\Gamma)) + g^{\prime\prime}(C_\Gamma){\bm J}_{C_\Gamma}\cdot\nabla_s C_\Gamma + g^{\prime}(C_\Gamma)S \big\}d\mathcal{H}^{d-1}  \nonumber
    \\
    &+ \int_{\partial\Gamma}(\tilde{\gamma}(C_\Gamma){\bm u}-g^{\prime}(C_\Gamma){\bm J}_{C_\Gamma})\cdot{\bm m} d\mathcal{H}^{d-2}, \label{eq:Surface_Energy_Derivative1}
  \end{align}
  where $\tilde{\gamma}(C_\Gamma) := g(C_\Gamma)-g^{\prime}(C_\Gamma)C_\Gamma$ is the Legendre transform of $g(C_\Gamma)$.
  Specially, when $\Gamma$ is close, Eq. \eqref{eq:Surface_Energy_Derivative1} reduces to 
  \begin{equation}
    \mathrm{  \mathrm{\frac{d}{dt}}}\int_{\Gamma} g(C_\Gamma) d\mathcal{H}^{d-1} = \int_{\Gamma} \big\{{\bm u}\cdot(\tilde{\gamma}(C_\Gamma)\kappa{\bm n}-\nabla_s\tilde{\gamma}(C_\Gamma)) + g^{\prime\prime}(C_\Gamma){\bm J}_{C_\Gamma}\cdot\nabla_s C_\Gamma + g^{\prime}(C_\Gamma)S \big\}d\mathcal{H}^{d-1}. \label{eq:Surface_Energy_Derivative2}
  \end{equation}
\end{lma}
  \begin{proof}
   The details of the derivation of Eq. \eqref{eq:Surface_Energy_Derivative1} are as follows
    \begin{align}
    \mathrm{  \mathrm{\frac{d}{dt}}}\int_{\Gamma} g(C_\Gamma) d\mathcal{H}^{d-1} &= \int_{\Gamma}\big[g^{\prime}(C_\Gamma)\partial^\bullet_t C_\Gamma +g(C_\Gamma)\nabla_s\cdot{\bm u} \big] d\mathcal{H}^{d-1}\nonumber
    \\
    &=\int_{\Gamma} \tilde{\gamma}(C_\Gamma) \nabla_s\cdot{\bm u} d\mathcal{H}^{d-1} - \int_{\Gamma} g^{\prime}(C_\Gamma) \nabla_s\cdot{\bm J}_{C_\Gamma} d\mathcal{H}^{d-1} +  \int_{\Gamma} g^{\prime}(C_\Gamma)S d\mathcal{H}^{d-1}\nonumber
    \\
    &= \int_{\Gamma} \big\{{\bm u}\cdot(\tilde{\gamma}(C_\Gamma)\kappa{\bm n}-\nabla_s\tilde{\gamma}(C_\Gamma)) + g^{\prime\prime}(C_\Gamma){\bm J}_{C_\Gamma}\cdot\nabla_s C_\Gamma + g^{\prime}(C_\Gamma)S \big\}d\mathcal{H}^{d-1}  \nonumber
    \\
    &+ \int_{\partial\Gamma}(\tilde{\gamma}(C_\Gamma){\bm u}-g^{\prime}(C_\Gamma){\bm J}_{C_\Gamma})\cdot{\bm m} d\mathcal{H}^{d-2}, \label{eq:Surface_Energy_Derivative3}
    \end{align}
   where we have used the surface transport formula, Eq. \eqref{eq:Surface_Mass_Conservation_Diff1}, Eq. \eqref{eq:Surface_Divergence_Formula2} and ${\bm n}\cdot{\bm J}_{C_\Gamma}=0$.
  \end{proof}

\section{Some inequalities related with the energy dissipation}

\begin{lma}
\label{lma:Adsorption_Desorption_Source}
  For the given concentration functions $C^\pm$ defined in $\Omega^\pm$ and $C_\Gamma$ defined on $\Gamma$, there exists a function defined on $\Gamma$ 
  \begin{equation}
    S^\pm = k_{ad}^\pm\frac{C^\pm}{C_\infty} - k_d^\pm\left(\frac{C_\Gamma}{C_{\Gamma,\infty}}\right)^{\omega_c},
  \end{equation}
  so that the following inequality holds
  \begin{equation}
    \int_{\Gamma} (g^{\prime}(C_\Gamma)-f^{\prime}(C^\pm))S^\pm d\mathcal{H}^{d-1} \leq 0,
  \end{equation}
  where $\lambda_a = \exp\left(\frac{U_C-U_{C_\Gamma}}{RT}\right)$, $k_{ad}^\pm=\lambda_a k_d^\pm$ and $k_d^\pm$ are the adsorption and desorption coefficients respectively.
\end{lma}
\begin{proof}
  It is obvious that 
  \begin{align}
    [g^{\prime}(C_\Gamma)-f^{\prime}(C^\pm)]S^\pm &= \bigg[U_{C_\Gamma}-U_C + RT\ln\bigg(\frac{(C_\Gamma/C_{\Gamma,\infty})^{\omega_c}}{C^\pm/C_\infty}\bigg)\bigg]S^\pm =  RT\ln\bigg(\frac{(C_\Gamma/C_{\Gamma,\infty})^{\omega_c}}{\lambda_a C^\pm/C_\infty}\bigg)S^\pm.
  \end{align}
  By the simple inequality $(a-b)\ln(a/b) \geq 0$ for any $a > 0$ and $b > 0$, we can set up
\begin{equation}
  a = \lambda_a C^\pm/C_\infty,\quad b = (C_\Gamma/C_{\Gamma,\infty})^{\omega_c},
\end{equation}
such that Lemma \ref{lma:Adsorption_Desorption_Source} holds. 
\end{proof}

\begin{lma}
  \label{lma:Bulk_Concentration_Jump}
  Similarly, for the given concentration functions $C^\pm$ defined in $\Omega^\pm$, there exists a function defined on $\Gamma$
  \begin{equation}
    J_s = k_c\frac{[\![C]\!]}{C_\infty},
  \end{equation}
   so that the following inequality holds
  \begin{equation}
    \int_{\Gamma} [\![f^{\prime}(C)]\!] J_s d\mathcal{H}^{d-1} \geq 0,
  \end{equation}
  where $k_c$ is the mass transfer coefficient and $[\![f^{\prime}(C)]\!] = f^{\prime}(C^+)-f^{\prime}(C^-)$.
\end{lma}

\begin{lma}
  \label{lma:Chemical_Reaction_Source}
  For the given concentration functions $A_{\Gamma}$, $B_{\Gamma}$ and $C_{\Gamma}$ defined on $\Gamma$, there exists a function defined on $\Gamma$ 
  \begin{equation}
    \mathcal{R} = k_f\left(\frac{B_{\Gamma}}{B_{\Gamma,\infty}}\right)^{\omega_b}\left(\frac{C_{\Gamma}}{C_{\Gamma,\infty}}\right)^{\omega_c} - k_r\left(\frac{A_{\Gamma}}{A_{\Gamma,\infty}}\right)^{\omega_a},
  \end{equation}
  so that the following inequality holds
  \begin{equation}
    \int_{\Gamma} (g^{\prime}(A_\Gamma)-g^{\prime}(B_\Gamma)-g^{\prime}(C_\Gamma))\mathcal{R} d\mathcal{H}^{d-1} \leq 0,
  \end{equation}
  where $\lambda_c = \exp\left(\frac{U_{B_\Gamma}+U_{C_\Gamma}-U_{A_\Gamma}}{RT}\right)$ and $k_f = \lambda_c k_r$.
\end{lma}
\begin{proof}
It is easy to compute that 
\begin{align}
  g^{\prime}(A_\Gamma) - g^{\prime}(B_\Gamma)-g^{\prime}(C_\Gamma) & = (U_{A_\Gamma}-U_{B_\Gamma}-U_{C_\Gamma}) + RT\ln\left(\frac{\left(A_\Gamma/A_{\Gamma,\infty}\right)^{\omega_a}}{\left(B_\Gamma/B_{\Gamma,\infty}\right)^{\omega_b}\left(C_\Gamma/C_{\Gamma,\infty}\right)^{\omega_c}}\right) \nonumber
  \\
  & = RT\ln\bigg(\frac{(A_\Gamma/A_{\Gamma,\infty})^{\omega_a}}{\lambda_c(B_\Gamma/B_{\Gamma,\infty})^{\omega_b}(C_\Gamma/C_{\Gamma,\infty})^{\omega_c}}\bigg).\nonumber
\end{align}
Further, we have 
\begin{equation}
  [g^{\prime}(A_\Gamma) - g^{\prime}(B_\Gamma)-g^{\prime}(C_\Gamma)]\mathcal{R} = RT\ln\bigg(\frac{(A_\Gamma/A_{\Gamma,\infty})^{\omega_a}}{\lambda_c(B_\Gamma/B_{\Gamma,\infty})^{\omega_b}(C_\Gamma/C_{\Gamma,\infty})^{\omega_c}}\bigg)\mathcal{R}.
\end{equation}
By the simple inequality $(a-b)\ln(a/b) \geq 0$ for any $a > 0$ and $b > 0$, we can set up
\begin{equation}
  a = \lambda_c(B_\Gamma/B_{\Gamma,\infty})^{\omega_b}(C_\Gamma/C_{\Gamma,\infty})^{\omega_c},\quad b=(A_\Gamma/A_{\Gamma,\infty})^{\omega_a},
\end{equation}
so that Lemma \ref{lma:Chemical_Reaction_Source} holds. 
\end{proof}

\end{appendices}

\end{document}